\documentstyle[amsfonts,amssymb]{article}

\newcommand{\mb}{\mathbb}

\newcommand{\be}{\begin{equation}}
\newcommand{\ee}{\end{equation}}
\newcommand{\beq}{\begin{eqnarray}}
\newcommand{\eeq}{\end{eqnarray}}
\newcommand{\beqst}{\begin{eqnarray*}}
\newcommand{\eeqst}{\end{eqnarray*}}

\newcommand{\wt}{\widetilde}

\newcommand{\iy}{\infty}

\newcommand{\ve}{\varepsilon}
\newcommand{\pa}{\partial}

\def\R{{\mathbb R}}

\newcommand{\dsp}{\displaystyle}

\newtheorem{theorem}{Theorem}[section]
\newtheorem{lemma}[theorem]{Lemma} 
 
\newtheorem{corollary}[theorem]{Corollary} 
\newtheorem{remark}[theorem]{Remark} 
\newtheorem{proposition}[theorem]{Proposition}

\title{\bf Fundamental Solutions for  Wave Equation in \\  de~Sitter Model of Universe\\
}

\author{{\bf Karen Yagdjian\thanks{Corresponding author. 
\newline{\hspace*{0.5cm}{\it E-mail addresses}: yagdjian@utpa.edu (K.Yagdjian), agalstyan@utpa.edu (A.Galstian).}} \,\,  and  Anahit Galstian} }

\begin{document}

\date{}
\maketitle
\thispagestyle{empty} 
\vspace{-0.3cm}

\begin{center}  {\it Department of Mathematics, 
University of Texas-Pan American,  
1201 W.~University Drive, \\
Edinburg, TX 78541-2999,  
USA }
\end{center}

\addtocounter{section}{-1}
\renewcommand{\theequation}{\thesection.\arabic{equation}}
\setcounter{equation}{0}
\pagenumbering{arabic}
\setcounter{page}{1}
\thispagestyle{empty}

\begin{abstract}
\begin{small}
 In this article we construct the fundamental solutions  for the wave equation arising in the de~Sitter 
model of the universe.  We use the fundamental solutions  to  represent     solutions of the Cauchy problem and to prove  the $L^p-L^q$-decay estimates 
for the solutions of the equation with and without a source term. 
\medskip

\medskip


\end{small}
\end{abstract}

\section{Introduction and Statement of Results}
\label{S1}

\setcounter{equation}{0}
\renewcommand{\theequation}{\thesection.\arabic{equation}}

In this paper we construct the fundamental solutions  for the wave equation arising in the de~Sitter model of the
universe and use the fundamental solutions  to find representations of the solutions to the Cauchy problem as well as the decay rates for them.
\smallskip

After averaging on a suitable scale, our universe is homogeneous 
and isotropic; therefore, the properties of the universe can be properly described by treating the matter as a perfect homogeneous fluid. 
In the models of the universe proposed by Einstein \cite{Einstein} and de~Sitter \cite{Sitter} the line
element is  connected with the proper 
mass density  and the proper pressure in the universe by the field equations
for a perfect fluid. There are two alternatives,  
which lead to the solutions of Einstein    and de~Sitter, respectively \cite[Sec.132]{Moller}. 
\smallskip

The homogeneous and isotropic cosmological models possess highest symmetry that makes them  more amenable to rigorous study.
Among them we mention FLRW (Friedmann-Lematre-Robertson-Walker) models. The simplest class of cosmological models 
can be obtained if we assume additionally that the metric of the slices of constant time is flat and that 
the spacetime metric can be written in the form
\[
ds^2= -dt^2+ a^2(t)( d{x}^2   + d{y} ^2 +  d{z}^2 )
\]
with an appropriate scale factor $a(t)$. Although on the made  assumptions,   the  spatially flat FLRW models appear to give a good 
explanation of our universe. The assumption that the universe is expanding leads to the positivity of  
  the time derivative $\frac{d }{dt}a (t)$. A further assumption that the universe obeys the accelerated expansion suggests 
that the second derivative  $\frac{d^2 }{dt^2}a (t)$  is positive.  
A substantial amount of the observational material can be satisfactorily interpreted in terms of the models, 
which take into account existing acceleration of   the recession of distant galaxies.
\smallskip

The time dependence of the function $a(t)$ is  determined by  the Einstein field equations for gravity. 
The Einstein equations with the {\it cosmological constant} $\Lambda $ have form
\[
R_{\mu \nu }-\frac{1}{2} g_{\mu \nu }R = -8\pi GT_{\mu \nu }-\Lambda g_{\mu \nu },
\] 
where term $ \Lambda g_{\mu \nu }$ can be interpreted as an energy-momentum of the vacuum. 
 Even a small value of 
$\Lambda $ could have drastic effects on the evolution of the universe. Under the assumption of FLRW symmetry the equation of motion in the case of positive
cosmological constant \, $\Lambda $ \, leads to solution
\[
a(t)=a(0)e^{\sqrt{\frac{\Lambda }{3}}t},
\]
which produces models with exponentially accelerated expansion. The model described by the last equation is usually called the {\it de Sitter
model}.
\smallskip

The unknown of principal importance in the Einstein equations is 
a metric $g$. It comprises the basic geometrical feature of the gravitational field, 
and consequently explains the phenomenon of the mutual gravitational attraction of substance. 
In the presence of matter these equations contain a non-vanishing right hand side $-8\pi GT_{\mu \nu }$. 
In general, the matter fields described by the function  $\phi $  must satisfy some equations of motion, and 
in the case of the scalar field, the equation of motion is that $\phi $ 
should satisfy the  wave equation generated by the metric $g$. 
In the  de~Sitter universe  the equation for the scalar field with mass \,  $m$\,   and potential function \, $V$   \, 
is (See, e.g. \cite{Friedrich-Rendall,Rendall_2004}.)
\be
\label{1.3}
  \phi_{tt} +   n  \phi_t - e^{-2t} \bigtriangleup \phi = - m^2\phi  - V'(\phi ),
\ee
while for the massless scalar field   the equation is
\be
\label{Sitter_nl}
u_{tt} +   n  u_t - e^{-2t} \bigtriangleup u  = -  V'( u ).
\ee
Here $\bigtriangleup $ is the Laplacian on the flat metric. The time inversion 
 transformation $t \to -t$ reduces the last equation to the mathematically equivalent  equation
\be
\label{Gal_nl}
u_{tt} - n  u_t- e^{2t} \bigtriangleup u  = - V'(u ).
\ee
Thus, written out explicitly in coordinates the wave equation on de~Sitter spacetime takes the form 
\be
\label{1.10}
 u_{tt} + n H   u_t - e^{-2Ht}\bigtriangleup u =0\,.
\ee
In \cite{Rendall_2004} the following ansatz for the formal solutions of the last equation is suggested
\[
\sum_{m=0}^\infty \Big( A_m(x)e^{-mHt} +  B_m(x)t e^{- mHt} \Big) .
\]
It is shown that such solutions can be parametrized by $A_0$ and $A_n$. It is also claimed in \cite{Rendall_2004}  that any solution has an asymptotic expansion of the type derived on a formal level. 
\smallskip

In the case of de~Sitter universe  the line element   may be written \cite[Sec.134]{Moller}
\[
ds^2= -   c^2\, d{t }^2+ e^{2ct /R}( d{x }^2   + d{y } ^2 +  d{z  }^2)\,.
\]
The  coordinates \, $t $, $x $, $y  $, $z $ \, can take all values from $-\iy$ to $\iy$.
Here $R$ is the ``radius'' of the universe. 
The de~Sitter model allows us to get an explanation of the actual red shift 
of spectral lines observed by Hubble and Humanson \cite{Moller}. 
In a certain sense all  solutions  look like the de~Sitter solution at late times \cite{Heinzle-Rendall}.  
We  write the de~Sitter line element in the 
form 
\[
ds^2= -dt^2+ e^{2Ht}( d{x}^2   + d{y} ^2 +  d{z}^2 )\,,
\]
where $H = \sqrt{\Lambda /3}$ is Hubble constant. 
The spacetime metric in the higher dimensional analogue of de~Sitter space is
\[
ds^2 =-dt^2+e^{2Ht}\big( (dx^1)^2+ \ldots + (dx^n)^2\big)\,.
\]
It is a simplified version of the multidimensional cosmological models with the metric tensor given by 
\[
g =-e^{2\gamma (t)}dt^2+  e^{2\phi ^1(t)}g_1+ \ldots +  e^{2\phi ^n(t)}g_n \,,
\]
and can be chosen as a starting point for the study.  
The multidimensional cosmological models have attracted a lot of attention during recent years in constructing 
mathematical models of an anisotropic universe (see, e.g. \cite{Brozos-Vázquez,Heinzle-Rendall} and references therein). 
\smallskip

We take a principal part of the   equation (\ref{1.10}) as an initial model that can be treated
first:
\be
\label{1.11}
\pa_t^2 u  - e^{-2Ht}\bigtriangleup u =0\,.
\ee
For simplicity, we set $H=1$. The time inversion 
 transformation $t \to -t$ reduces the last equation to the mathematically equivalent  equation
\be
\label{G-S}
\pa_t^2 u  - e^{2t}\bigtriangleup u =0\,.
\ee
Hence, if we can find the fundamental solution for the linear equation (\ref{G-S})
associated with (\ref{Gal_nl}), then it generates  the fundamental solution for the linear equation (\ref{1.11}) associated with (\ref{Sitter_nl}).
\smallskip

The  equation (\ref{G-S}) is strictly hyperbolic. That implies the well-posedness of the Cauchy problem for (\ref{G-S})  
in the different functional spaces. The coefficient of the equation is an analytic function and 
Holmgren's theorem  implies a local uniqueness in the space of distributions. 
Moreover, the speed of propagation is finite, namely, 
it is equal to $e^{t} $ for every $ t \in {\mb R}$. 
The second-order strictly hyperbolic  equation (\ref{G-S}) possesses two fundamental solutions 
resolving the Cauchy problem. They can be written microlocally in terms of the Fourier integral operators \cite{Horm}, which
give a complete description of the wave front sets of the solutions.
The distance between two characteristic roots $\lambda_1 (t,\xi ) $ and $\lambda_2 (t,\xi ) $ of the equation (\ref{G-S}) is  
\[
|\lambda_1 (t,\xi ) - \lambda_2 (t,\xi )| = e^{t}|\xi |, \qquad t \in {\mb R}, \,\, \xi \in {\mb R}^n\,.
\]
It tends to zero as $t $ approaches $-\infty$. Thus, the operator is not uniformly (that is for all $t \in {\mb R} $) strictly hyperbolic. 
Moreover, the finite integrability
of the characteristic roots, $\int_{-\infty}^0 |\lambda_i (t,\xi )| dt <  \infty $,  leads to the 
existence of so-called ``horizon'' for that equation. More precisely,  any  signal emitted from the spatial point $x_0 \in {\mb R}^n$
at time $t_0 \in {\mb R} $ remains inside the ball $ |x -x_0 | <e^{t_0} $ for all time $t \in (-\infty, t_0) $. The  equation (\ref{G-S})
is neither Lorentz  invariant nor  invariant with respect to usual scaling and that brings additional difficulties. 
In particular, it can cause a nonexistence of the $L^p-L^q$ decay 
for the solutions in the backward direction of time. In \cite{yagdjian_birk} it is mentioned the model equation with permanently bounded domain of influence,   power decay of characteristic roots,  and without
$L^p-L^q$ decay 
for the solutions that illustrates that phenomenon. The above mentioned $L^p-L^q$ decay 
estimates are some of the important tools for studying nonlinear equations
(see, e.g. \cite{Racke,Shatah}).
\smallskip

The  equation  (\ref{G-S})
 was investigated in \cite{Galstian2001,Galstian} by the second author. More precisely,  in \cite{Galstian2001,Galstian} the resolving operator for the Cauchy problem
\be 
\label{0.5new_intr}
\pa_t^2 u  - e^{2t}\bigtriangleup u =0, \qquad u(x,0)= \varphi_0 (x), \quad u_t(x,0)= \varphi_1 (x)\,,
\ee
is written as a sum of the Fourier integral operators with the amplitudes given in terms of the 
Bessel functions and in terms of confluent hypergeometric functions. In particular,  it is proved  in \cite{Galstian2001,Galstian} that for $t>0$ 
the solution of the Cauchy problem  (\ref{0.5new_intr})   
is given by 
\beqst
\hspace*{-0.5cm}  u(x,t) 
& = &
- i\frac{2}{(2\pi)^{n}} \int_{\R^n}
\Big\{ 
e^{i[x\cdot \xi +( e^t-1) |\xi|] } H_{+}\big(\frac{1}{2} ;1;2ie^t |\xi| \big) 
H_{-}\big(\frac{3}{2} ; 3;2i|\xi| \big) \nonumber \\
&  &
\hspace{1.8cm}-\,
e^{i[x\cdot \xi  -(e^t-1) |\xi| ]} H_{-}\big(\frac{1}{2} ;1;2ie^t |\xi| \big)
H_{+}\big(\frac{3}{2} ;3;2i|\xi|  \big)\Big\}
|\xi|^2 {\mathcal F}(\varphi _0)(\xi) d\xi  \nonumber \\
&   &
- i\frac{1}{(2\pi)^{n}} \int_{\R^n}
 \Big\{ 
e^{i[x\cdot \xi +( e^t-1) |\xi|] } H_{+}\big(\frac{1}{2} ;1;2ie^t |\xi| \big) 
H_{-}\big(\frac{1}{2} ; 1;2i|\xi| \big) \nonumber \\
&  &
\hspace{1.8cm}
-\,e^{i[x\cdot \xi  -(e^t-1) |\xi| ]} H_{-}(\frac{1}{2} ;1;2ie^t |\xi| )
H_{+}\big(\frac{1}{2} ;1;2i|\xi|  \big)\Big\}
{\mathcal F}(\varphi _1)(\xi) d\xi \,.
\eeqst
In the notations of \cite{B-E} the last functions are   $H_{-}(\alpha ;\gamma ;z) = e^{i \alpha  \pi } \Psi (\alpha ;\gamma ;z) $
and $H_{+}(\alpha ;\gamma ;z) = e^{i \alpha  \pi } \Psi (\gamma -\alpha ;\gamma ;-z) $, where function $ \Psi ( a ;c ; z)  $ 
is defined in \cite[Sec.6.5]{B-E}. Here ${\mathcal F}(\varphi  )(\xi) $ is a Fourier transform of $ \varphi   (x) $. 
\smallskip

The typical $L^p-L^q$ decay estimates 
obtained in \cite{Galstian2001,Galstian} by dyadic decomposition of the phase space  
 contain some loss of regularity.
More precisely, it is proved that for the solution $u = u(x,t)$ to the Cauchy problem (\ref{0.5new_intr}) with $n \ge 2$, 
$\varphi _0(x) \in C_0^\iy ({\mb R}^n)$ and $\varphi _1(x) =0$ 
for all large $t\ge T >0$,  the following estimate is satisfied
\be
\label{1.16}
\|u(x,t)\|_{L^q({\mb R}^n)} 
  \le  
C (1+ e^t )^{ 
- \frac{1}{2} (n-1) (\frac{1}{p}-\frac{1}{q})}
\|\varphi _0\|_{W^N_p( {\mb R}^n)},
\ee
 where $1<p\le 2$, $\frac{1}{p}+\frac{1}{q}=1$, and 
$\frac
{1}{2} (n+1) (\frac{1}{p}-\frac{1}{q})\le N < \frac
{1}{2} (n+1) (\frac{1}{p}-\frac{1}{q})+1$ and $W^N_p( {\mb R}^n) $ is the Sobolev space.
In particular, the loss of regularity, $N$, is positive, unless $p=q=2$. 
This loss of regularity phenomenon exists for the classical wave equation as well.
 Indeed, it is well-known  (see, e.g., \cite{Littman_1963,Littman,Peral}) that for the Cauchy problem
$ u_{tt}-\bigtriangleup u=0\,, \quad u(x,0)=\varphi (x), \quad u_t(x,0)= 0 $,
the estimate
$ \| u(x,t)  \| _{ L^q  ({\mathbb R}_x^n)  }  \le C  \| \varphi  (x) \|_{L^q({\mathbb R}_x^n)}   $
fails to fulfill even for small positive  $t$ unless $q=2$. The obstacle  is created by the 
distinguishing feature of the (different from translation) Fourier integral operators of order zero, which   compose a resolving operator.
\smallskip

According to Theorem~1~\cite{Galstian}, for the solution $u = u(x,t)$ to the Cauchy problem (\ref{0.5new_intr}) with  $n \ge 2$, 
$\varphi _0(x) =0$ and $\varphi _1(x) \in C_0^\iy ({\mb R}^n)$ 
for all large $t\ge T >0$ and for any small $\ve>0$, the following estimate is satisfied
\[
\|u(x,t)\|_{L^q({\mb R}^n)} 
 \le 
C _\ve(1+t)(1+ e^t )^{ r_0- n(\frac{1}{p}-\frac{1}{q})}
\|\varphi _1\|_{W^N_p( {\mb R}^n)},\,
\]
where $1<p\le 2$, $\frac{1}{p}+\frac{1}{q}=1$, 
$r_0 = \max \{\ve; \frac{(n+1)}{2} (\frac{1}{p}-\frac{1}{q})-\frac{1}{q}\}$,  
$\frac{n+1}{2}(\frac{1}{p}-\frac{1}{q})-\frac{1}{q}\le N < \frac{n+1}{2}(\frac{1}{p}-\frac{1}{q})+\frac{1}{p}$.
\smallskip

The nonlinear equations (\ref{1.3}) and (\ref{Sitter_nl})
are those we would like to solve, but the linear problem is a natural first step.
Exceptionally efficient tool for the studying nonlinear equations is a   fundamental 
solution of the associate linear operator.  
\smallskip

In the construction of the fundamental solutions for the operator (\ref{G-S}) 
 we follow the approach proposed in \cite{YagTricomi} that allows us to represent  
the fundamental solutions as some integral of the family of the fundamental solutions of the Cauchy problem for the wave equation without source term.
The kernel of that integral contains Gauss's hypergeometric function. In that way, many properties of the wave equation can be extended to the
hyperbolic equations with the time dependent speed of propagation. That approach  was successfully applied in \cite{YagTricomi_GE,YagTricomi_JMAA} by the first author to
investigate the semilinear Tricomi-type equations.
\smallskip

The   operator of the equation (\ref{G-S}) is 
\[
{\mathcal S}:= \pa_t^2  - e^{2t}\bigtriangleup   \,,
\] 
where $x \in {\mb R}^n$, $t \in {\mb R}$, and $ \bigtriangleup $ is the Laplace operator, $\bigtriangleup := \sum_{j=1}^n \frac{\pa^2 }{\pa x_j ^2} $. We look for the  fundamental solution (Green's function, propagator in the literature on Physics) $E=E(x,t;x_0,t_0)$,
\[
E_{tt} - e^{2t}\Delta E = \delta (x-x_0,t-t_0),
\]
with a support in the ``forward light cone'' $D_+ (x_0,t_0) $, $x_0 \in {\mb R}^n$, $t_0 \in {\mb R}$,
and for the  fundamental solution with a support in the ``backward light cone'' $D_- (x_0,t_0) $, $x_0 \in {\mb R}^n$, $t_0 \in {\mb R}$,
defined as follows
\beq
\label{D+}
D_+ (x_0,t_0) 
& := &
\Big\{ (x,t)  \in {\mb R}^{n+1}  \, ; \, 
|x -x_0 | \leq e^t - e^{t_0} 
\,\Big\} \,,\\
\label{D-}
D_- (x_0,t_0) 
& := &
\Big\{ (x,t)  \in {\mb R}^{n+1}  \, ; \, 
|x -x_0 | \leq -(e^t - e^{t_0}) 
\,\Big\} \,.
\eeq
In fact,  any intersection of  $ D_- (x_0,t_0) $ with the hyperplane $t=const <t_0$ determines the so-called dependence domain
for the point $(x_0,t_0) $, while the  intersection of  $ D_+ (x_0,t_0) $ 
with the hyperplane $t=const >t_0$ is the so-called  domain of influence of the point $(x_0,t_0) $. The equation (\ref{G-S}) is
non-invariant    with respect to time inversion. Moreover, the domain of influence is wider than any given ball if time $const>t_0$ is sufficiently large, while 
the dependence domain is permanently, for all time  $const< t_0$, in the   ball of the radius $e^{t_0} $.
\smallskip

Define for $t_0 \in {\mb R}$ in the domain  $D_+ (x_0,t_0)\cup D_- (x_0,t_0) $  the function   
\beq
\label{E}
\hspace*{-0.5cm} E(x,t;x_0,t_0)
& := &
 \frac{1}{\sqrt{  (e^{t_0}  + e^t)^2  -(x-x_0)^2 }}F \Big(\frac{1}{2}, \frac{1}{2};1; \frac{(e^t  - e^{t_0})^2 - (x-x_0)^2}
{(e^t  + e^{t_0})^2 - (x-x_0)^2}  \Big) \,,
\eeq
where $F\big(a, b;c; \zeta \big) $ is the hypergeometric function (See, e.g. \cite{B-E}.).
Let $E(x,t;0,b)$ be a function  (\ref{E}),
and set
\[ 
E_+(x,t;0,t_0) 
 := 
\cases{ E(x,t;0,t_0) \quad \mbox{\rm in}  \,\, D_+ (0,t_0), \cr
0 \hspace*{2.0cm} \mbox{\rm elsewhere}} \,, \qquad
E_{-}(x,t;0,t_0) 
 :=  
\cases{ E(x,t;0,t_0) \quad \mbox{\rm in}  \,\, D_- (0,t_0), \cr
0 \hspace*{2.0cm} \mbox{\rm elsewhere}} \,.
\]
Since  function $E=E(x,t;0,t_0)$ is smooth in $D_{\pm} (0,t_0) $, 
it follows that $E_+(x,t;0,$ $t_0) $ and $E_-(x,t;0,$ $t_0) $ are   locally integrable functions and they define
 distributions whose supports are in   $D_{+}(0,t_0)$  and 
 $D_{-} (0,t_0) $, respectively. The next theorem  gives our first result.

\begin{theorem}
\label{T1}\label{T2.1_intr}
Suppose that $n=1$. The distributions $E_+(x,t;0,t_0) $ and  $E_-(x,t;0,t_0) $  are the fundamental solutions for the operator
${\mathcal S}:= \pa_t^2-e^{2t}\pa_x^2$ relative to point $(0,t_0)$, that is
\[
{\mathcal S} E_{\pm}(x,t;0,t_0) = \delta (x,t-t_0) \qquad \mbox{or} \qquad
\frac{\pa^2}{\pa t^2}E_{\pm}(x,t;0,t_0) - e^{2t} \frac{\pa^2}{\pa x^2}E_{\pm}(x,t;0,t_0)= \delta (x,t-t_0).
\]
\end{theorem}
\smallskip

To motivate one construction for the higher dimensional case $n>1$ we follow the approach suggested in \cite{YagTricomi} and represent fundamental solution  
$E_+(x,t;0,t_0) $  as follows 
\[
E_+(x,t;0,t_0) 
  =  
 \int_{ e^{t_0}- e^{t} }^  { e^t- e^{t_0}} E^{string} (x,r )   \frac{1}{\sqrt{(e^t  + e^{t_0} )^2  -r^2}}  
F\left(\frac{1}{2},\frac{1}{2};1; 
\frac{  (e^t  - e^{t_0} )^2    -r^2}
{(e^t  + e^{t_0} )^2    -r^2} \right)\, dr , \quad t>t_0,
\]
where the distribution $E^{string} (x,t )  $ is the fundamental solution of the Cauchy problem for the string equation:
\[
\frac{\pa^2 }{\pa t^2}  E^{string}   -   \frac{\pa^2 }{\pa x^2}  E^{string} =0, \qquad  E^{string}(x,0 )= \delta (x)
, \,\,\, E^{string}_t  (x,0 )=0\,.
\] 
Hence, $E^{string} (x,t ) = \frac{1}{2}\{ \delta (x+t )+ \delta (x-t )\}  $. The kernel (\ref{E}) is the even 
function of $x$ while $E^{string} (x,t ) $ is even with respect to $t$. The integral makes sense in the topology of the space of distributions.
The fundamental solution $E_-(x,t;0,t_0) $ for $t<t_0$ admits a similar representation.
\smallskip

We appeal to
the wave equation in Minkowski spacetime to obtain in the next theorem very similar representations of the fundamental solutions of the 
higher dimensional equation in  de~Sitter spacetime with $n>1$.
\begin{theorem}
If $x \in {\mb R}^n$, $n>1$, and $ \bigtriangleup $ is the Laplace operator, then for the operator
\[
{\mathcal S}:= \frac{\pa^2 }{\pa t ^2} - e^{2t}\bigtriangleup 
\] 
the  fundamental solution $E_{+,n}(x,t;x_0,t_0) $ $(= E_{+,n}(x-x_0,t;0,t_0))$, with a 
support in the  forward cone  $D_+ (x_0,t_0) $, $x_0 \in {\mb R}^n$, $t_0 \in {\mb R}$, \, 
supp$\,E_{+,n} \subseteq D_+ (x_0,t_0)$, is given by the following integral ($t>t_0$)
\be
\label{E+}
E_{+,n}(x-x_0,t;0,t_0) 
 = 
2   
  \int_{ 0}^{ e^t- e^{t_0} }   E^w (x-x_0,r )   \frac{1}{\sqrt{(e^t  + e^{t_0} )^2  -r^2}}  
F\left(\frac{1}{2},\frac{1}{2};1; 
\frac{  (e^t  - e^{t_0} )^2    -r^2}
{(e^t  + e^{t_0} )^2    -r^2} \right)\, dr .
\ee
Here the function 
$E^w(x,t;b)$   
is a fundamental solution to the Cauchy problem for the  wave equation
\[
E^w_{ tt} -   \bigtriangleup E^w  =  0 \,, \quad E^w(x,0)=\delta (x)\,, \quad E^w_{t}(x,0)= 0\,.
\]
The  fundamental solution $E_{-,n}(x,t;x_0,t_0) $ ($= E_{-,n}(x-x_0,t;0,t_0)$) with a 
support in the  backward cone  $D_- (x_0,t_0) $, $x_0 \in {\mb R}^n$, $t_0 \in {\mb R}$, \, 
supp$\,E_{-,n} \subseteq  D_- (x_0,t_0)$, is given by the following integral ($t<t_0$)
\be
\label{E-}
E_{-,n}(x-x_0,t;0,t_0) 
 = 
- 2   
  \int_{ e^t- e^{t_0} }^  { 0} E^w (x-x_0,r )   \frac{1}{\sqrt{(e^t  + e^{t_0} )^2  -r^2}}  
F\left(\frac{1}{2},\frac{1}{2};1; 
\frac{  (e^t  - e^{t_0} )^2    -r^2}
{(e^t  + e^{t_0} )^2    -r^2} \right)\, dr .
\ee
\end{theorem}

In particular, the formula (\ref{E+}) shows that Huygens's Principle is not valid for the waves propagating in the de~Sitter model of the universe. Fields satisfying a wave equation in the de~Sitter model of universe can be accompanied by {\it tails} propagating inside the light cone. This phenomenon
will be discussed in the spirit of \cite{Sonego-Faraoni}  in the forthcoming paper.
\smallskip

Next we use Theorem~\ref{T1}  to solve the Cauchy problem for the one-dimensional  equation 
\be
\label{Int_3} \label{TricEq}
u_{tt} - e^{2t}u_{xx} = f(x,t)\,, \qquad t > 0\,, \quad x \in {\mb R} \,, 
\ee
with vanishing initial data,
\be
\label{Int_4} \label{4.2}
u(x,0)= u_t(x,0)= 0\,.
\ee
\begin{theorem}
\label{T1.1}\label{T4.1}
Assume that the function $f$ is continuous along with its all second order derivatives, 
and that for every fixed $t$ it has a compact support, supp$f(\cdot ,t) \subset {\mb R}$.
Then the function $u=u(x,t)$ defined by  
\[
u(x,t)   = 
 \int_{ 0}^{t} db \int_{ x - (e^t- e^b)}^{x+ e^t- e^b}  f(y,b)  \frac{1}{\sqrt{(e^t  + e^b)^2  - (x-y)^2}}  
F\left(\frac{1}{2},\frac{1}{2};1; 
\frac{ (e^t  - e^b    )^2  - (x-y)^2}
{(e^t  + e^b    )^2  - (x-y)^2} \right) dy \, \nonumber 
\]
is a $C^2$-solution to the Cauchy problem for the equation (\ref{Int_3}) with  vanishing initial data,
(\ref{Int_4}).
\end{theorem}
\smallskip

The representation of the solution of the Cauchy problem for the one-dimensional case ($n=1$) 
of the equation (\ref{G-S}) without source term is given by the next theorem. 
\begin{theorem}
\label{T1.3}\label{TInt_2}
The solution $u=u (x,t)$ of the Cauchy problem 
\be
\label{oneDphy01}
u_{tt} - e^{2t}u_{xx} =0\, ,\qquad u(x,0) = \varphi_0  (x) \,, \qquad u_t(x,0) =\varphi_1  (x)\,  ,
\ee
with \,$\varphi_0    ,  \varphi_1  \in C_0^\infty ({\mb R})$ can be represented as follows
\beqst
u(x,t)   
&  =   &
 \frac{1}{2} e ^{-\frac{t}{2}}  \Big[ 
\varphi_0   (x+ e^t- 1)  
+     \varphi_0   (x - e^t+ 1)  \Big]  
+ \int_{ 0}^{e^t-1} \big[ 
\varphi_0   (x - z)  
+     \varphi_0   (x  + z)  \big] K_0(z,t)\,  dz \\
&  &    
+ \,\,\int_{0}^{  e^t-1} \,\Big[      \varphi_1    (x- z)  +   \varphi _1   (x + z)    \Big] K_1(z,t) dz \,,
\eeqst
where the kernels  $K_0(z,t)   $    and $K_1(z,t)   $ are defined by
\beqst
K_0(z,t) 
&: = &  
- \Big( \frac{\partial}{\partial t_0}E(z,t;0,t_0)\Big)\Big|_{t_0=0}  \\
& = &
 -  \frac{1}{2  ((e^t-1)^2-z^2 ) \sqrt{(e^t+1)^2-z^2}} \\
&  &
\times \left[ (1-e^{2 t}+z^2 ) 
F\Big(-\frac{1}{2},\frac{1}{2};1; \frac{ (e^t-1)^2 -z^2   }{ (e^t+1)^2 -z^2 }  \Big) 
+   2  (e^t -1) F\Big(\frac{1}{2},\frac{1}{2};1; \frac{ (e^t-1)^2 -z^2   }{ (e^t+1)^2 -z^2 }  \Big)  \right] \\
K_1(z,t)
& :=  &
E(z,t;0,0)=\frac{1}{  \sqrt{(1+e^t)^2-z^2 } }    F\Big(\frac{1}{2},\frac{1}{2};1; 
\frac{ (e^t-1)^2 -z^2   }{ (e^t+1)^2 -z^2 }  \Big)  
 \,, \quad 0\leq z\leq  e^t-1\,.
 \eeqst 
\end{theorem}

The kernel  $K_0(z,t)   $ has singularity at $z = e^t-1$. The kernels $K_0(z,t)  $ and $K_1(z,t)  $ play leading roles in the derivation of
decay estimates. Their main properties are listed and proved in Section~\ref{S11}.

\smallskip

Next we turn to the higher-dimensional equation with $n>1$.

\begin{theorem}
\label{T1.5}
If $n$ is odd, $n=2m+1$,  $m  \in {\mathbb N}$,  then the solution $u= u(x,t)$ to the Cauchy problem   
\be
\label{1.26}
u_{tt} - e^{2t}\Delta u = f ,\quad u(x,0)= 0  , \quad u_t(x,0)=0,
\ee
with \, $ f \in C^\infty ({\mb R}^{n+1})$\, and with the vanishing
initial data is given by the next expression 
\beq
\label{1.27}
u(x,t) 
& = &
2   \int_{ 0}^{t} db
  \int_{ 0}^{ e^t- e^b} dr_1 \,  \Bigg( \!\!\frac{\partial }{\partial r} 
\Big(  \frac{1}{r} \frac{\partial }{\partial r}\Big)^{\frac{n-3}{2} } 
\frac{r^{n-2}}{\omega_{n-1} c_0^{(n)} }  \!\!\int_{S^{n-1} } f(x+ry,b)\, dS_y  \!\!
\Bigg)_{r=r_1}   \nonumber \\
\label{7.1ndim_intr}
  &  & 
\hspace*{2cm}\times \frac{1}{\sqrt{(e^t  + e^b)^2  -r_1^2}}  
F\left(\frac{1}{2},\frac{1}{2};1; 
\frac{  (e^t  - e^b)^2    -r_1^2}
{(e^t  + e^b)^2    -r_1^2} \right) ,
\eeq
where $c_0^{(n)} =1\cdot 3\cdot \ldots \cdot (n-2 )$.
Constant $\omega_{n-1} $ is the area of the unit sphere $S^{n-1} \subset {\mathbb R}^n$. 

If $n$ is even, $n=2m$,  $m  \in {\mathbb N}$,  then the solution $u= u(x,t)$  is given by the next expression
\beq
\label{1.28} 
u(x,t) 
& = &
2 \int_{ 0}^{t} db
  \int_{ 0}^{ e^t- e^b} dr_1 \,  \Bigg(  \!\!\frac{\partial }{\partial r} 
\Big( \frac{1}{r} \frac{\partial }{\partial r}\Big)^{\frac{n-2}{2} } 
\frac{2r^{n-1}}{\omega_{n-1} c_0^{(n)} }  \!\!\int_{B_1^{n}(0)} \frac{f(x+ry,b) }{\sqrt{1-|y|^2}} \, dV_y 
 \!\!\Bigg)_{r=r_1}    \nonumber  \\
  &  & 
\hspace*{2cm}\times \frac{1}{\sqrt{(e^t  + e^b)^2  -r_1^2}}  
F\left(\frac{1}{2},\frac{1}{2};1; 
\frac{  (e^t  - e^b)^2    -r_1^2}
{(e^t  + e^b)^2    -r_1^2} \right) .
\eeq  
Here $B_1^{n}(0) :=\{|y|\leq 1\} $ is the unit ball in ${\mathbb R}^n$, while $c_0^{(n)} =1\cdot 3\cdot \ldots \cdot (n-1)$.
\end{theorem}
Thus, in both cases, of even and odd $n$, one can write 
\be
\label{1.29} 
u(x,t) 
  =  
2   \int_{ 0}^{t} db
  \int_{ 0}^{ e^t- e^b} dr  \,  v(x,r ;b)   \frac{1}{\sqrt{(e^t  + e^b)^2  -r ^2}}  
F\left(\frac{1}{2},\frac{1}{2};1; 
\frac{  (e^t  - e^b)^2    -r ^2}
{(e^t  + e^b)^2    -r ^2} \right) ,
\ee
where the function 
$v(x,t;b)$   
is a solution to the Cauchy problem for the  wave equation
\[
v_{tt} -   \bigtriangleup v  =  0 \,, \quad v(x,0;b)=f(x,b)\,, \quad v_t(x,0;b)= 0\,.
\]

The next theorem gives representation of the solutions of  equation (\ref{G-S}) with the initial data prescribed at $t=0$.
\begin{theorem}
\label{T1.6}
The solution $u=u (x,t)$ to the Cauchy problem  
\be
\label{1.30CP}
u_{tt}-  e^{2t} \bigtriangleup u =0, \quad u(x,0)= \varphi_0 (x) , \quad u_t(x,0)=\varphi_1 (x) 
\ee
with \, $ \varphi_0 $,  $ \varphi_1 \in C_0^\infty ({\mb R}^n) $, $n>1$, can be represented as follows:
\beq
\label{1.30}
u(x,t) 
& = &
 e ^{-\frac{t}{2}} v_{\varphi_0}  (x, \phi (t))
+ \, 2\int_{ 0}^{1} v_{\varphi_0}  (x, \phi (t)s) K_0(\phi (t)s,t)\phi (t)\,  ds  \nonumber \\
& &
+\, 2\int_{0}^1   v_{\varphi _1 } (x, \phi (t) s) 
  K_1(\phi (t)s,t) \phi (t)\, ds 
, \quad x \in {\mathbb R}^n, \,\, t>0,\,\, \phi (t):= e^t-1\,,  
\eeq
 by means of the kernels  $K_0$ and $K_1$ are defined in Theorem~\ref{T1.3}. 
 Here for $\varphi \in C_0^\infty ({\mathbb R}^n)$ and for $x \in {\mathbb R}^n$, $n=2m+1$, $m \in {\mathbb N}$,  
\begin{eqnarray*}
 v_\varphi  (x, \phi (t) s) :=  
\Bigg( \frac{\partial}{\partial r} \Big( \frac{1}{r} \frac{\partial }{\partial r}\Big)^{\frac{n-3}{2} } 
\frac{r^{n-2}}{\omega_{n-1} c_0^{(n)} } \int_{S^{n-1}  } 
\varphi (x+ry)\, dS_y   \Bigg)_{r=\phi (t) s}   
\end{eqnarray*}
while for $x \in {\mathbb R}^n$, $n=2m$,  $m \in {\mathbb N}$ ,   
\[
v_\varphi  (x, \phi (t) s) :=  \Bigg( \frac{\partial }{\partial r}  
\Big( \frac{1}{r} \frac{\partial }{\partial r}\Big)^{\frac{n-2}{2} } 
\frac{2r^{n-1}}{\omega_{n-1} c_0^{(n)}} \int_{B_1^{n}(0)}  \frac{1}{\sqrt{1-|y|^2}}\varphi (x+ry)\, dV_y 
 \Bigg)_{r=s \phi (t)}  \,.
\]
The function $v_\varphi  (x, \phi (t) s)$  coincides with the value $v(x, \phi (t) s) $ 
of the solution $v(x,t)$ of the Cauchy problem
\[
v_{tt}-  \bigtriangleup v =0, \quad v(x,0)= \varphi (x), \quad v_t(x,0)=0\,.
\]
\end{theorem}

\smallskip

As a consequence  of the theorems above  we obtain in Sections~\ref{S12}-\ref{S13}
for $n>1$  the following decay estimate
\beq
\label{1.32}
\hspace{-0.7cm}
\| (-\bigtriangleup )^{-s} u(x,t) \|_{ { L}^{  q} ({\mathbb R}^n)  } \!\! & \!\!\le &  \!\!\!\!
C   e^{t (  2s-n(\frac{1}{p}-\frac{1}{q})) } \int_{ 0}^{t} (1+  t -b  )\|  f(x, b)  \|_{ { L}^{p} ({\mathbb R}^n)  } \,db \nonumber \\
&  &
+ 
 C(e^ t-1)^{2s-n(\frac{1}{p}-\frac{1}{q}) } \left\{ \| \varphi_0  (x) \|_{ { L}^{p} ({\mathbb R}^n)  } + \|\varphi_1  (x)  \|_{ { L}^{p} ({\mathbb R}^n)  }(1+ t ) (1-e^{-t}) \right\} 
\eeq
provided that $s \ge 0$, $1<p \leq 2$, $\frac{1}{p}+ \frac{1}{q}=1$, $\frac{1}{2} (n+1)\left( \frac{1}{p} - \frac{1}{q}\right) \leq 
2s \leq n \left( \frac{1}{p} - \frac{1}{q}\right) < 2s+1$. Moreover, according to Theorem~\ref{T10.1}
the estimate (\ref{1.32}) is valid   for $n=1$ and $s=0$ as well as if $\varphi_0  (x)=0 $ and $\varphi_1  (x)=0 $. Case of $n=1$,  $f (x,t) =0 $, and non-vanishing $\varphi_1  (x)  $ and $\varphi_1  (x)  $ is discussed in Section~\ref{S11}.
\smallskip

The paper is organized as follows. In Section~\ref{S2} we construct the fundamental solutions of the operator (\ref{G-S}) for the case of $n=1$.
Then in Section~\ref{S3} we apply the fundamental solutions to solve the 
Cauchy problem with the source term and with the vanishing initial data  given at $t=0$. 
More precisely, we give a representation formula for the solutions. In Section~\ref{S4} we prove several basic properties
of the function $E(x,t;y,b)$.  In Sections~\ref{S5}-\ref{S6} we use formulas of  
Section~\ref{S4} to derive and to complete the list of   representation formulas 
for the solutions of the Cauchy problem for the case of one-dimensional spatial variable. The higher-dimensional
equation with the source term is considered in  Section~\ref{S7}, where we derive a representation formula for the solutions of 
the Cauchy problem with the source term and with the vanishing initial data  given at $t=0$. In same section this formula is used to 
derive the fundamental solutions of the operator and to complete the proof of Theorem~\ref{T1.6}.
Then in Sections~\ref{S10}-\ref{S13} we establish the $L^p-L^q$ decay estimates. 
Applications of all these results to the nonlinear equations will be done in the forthcoming paper.

\section{Fundamental Solutions. Proof of Theorem \ref{T1}}
\label{S2}

\setcounter{equation}{0}

In the characteristic coordinates $l$ and $m$, 
\be
\label{Gel1.8}
l =  x+ e^{t}, \qquad
m =  x- e^{t}
\ee
the operator 
\[
{\mathcal S} := \frac{\pa^2 }{\pa t ^2} - e^{2t}\frac{\pa^2 }{\pa x ^2}
\]
 reads
\beqst
\frac{\pa^2 }{\pa t ^2} -  e^{2t}  \frac{\pa^2  }{\pa x^2 } 
 & = &  
- (l-m)^{2 } \Bigg\{ \frac{\pa^2 }{\pa l \, \pa m} - \frac{1}{2(l-m)} \Big( \frac{\pa }{\pa l} -
\frac{\pa }{\pa m} \Big) \Bigg\} \,.
\eeqst
Consider point $(x,t)=(0,b)$, then two backward characteristics meet the $x$
line at the points 
$x= a$  and $x= -a$,
$a:= \phi (b) $. 
Note that the point $(l,m)=(\phi (b) ,-\phi (b) )$ represents point $(0,b)$ in characteristic
coordinates. The following lemma is an analog of   (2.2)\cite{GelfandII},  
where the Tricomi equation 
is considered. 
\begin{lemma}
The function
\[
E(l,m;a,b)  = 
(l-b)^{-1/2} (a-m)^{-1/2} 
F\Big(\frac{1}{2},\frac{1}{2};1; \frac{(l-a)(m-b)}{(l-b)(m-a)} \Big)
\]
solves the equation
\be
\label{f5.1}
\Bigg\{ \frac{\pa^2 }{\pa l \, \pa m} - \frac{1}{2(l-m)}  \Big( \frac{\pa }{\pa l} -
\frac{\pa }{\pa m} \Big) \Bigg\} E(l,m;a,b) = 0\,.
\ee
\end{lemma}

\noindent
{\bf Proof.} Indeed, after simple calculations, taking into account  (23) of  \cite[v.1, Sec.2.8]{B-E} 
\beq
\label{Fderivative}
  \frac{d}{d z}   F \left(\frac{1}{2},\frac{1}{2};1;z\right)   
& = &
\frac{1 }{2z (1 -z) }F \left(-\frac{1}{2},\frac{1}{2};1;z\right)
- \frac{1 }{2z}F \left(\frac{1}{2},\frac{1}{2};1;z\right)\,,
\eeq
we obtain
\beqst
&  &
   \partial _l\left((l-b)^{-\frac{1}{2}}(a-m)^{-\frac{1}{2}}
  F \left(\frac{1}{2},\frac{1}{2};1;\frac{(l-a)(m-b)}{(l-b)(m-a)}\right)\right)  \\
 & = &
\frac{(a-m) F \left(-\frac{1}{2},\frac{1}{2};1;\frac{(l-a)(m-b)}{(l-b)(m-a)}\right)-(l-m) F \left(\frac{1}{2},\frac{1}{2};1;\frac{(l-a)(m-b)}{(l-b)(m-a)}\right)}{2(l-a)
\sqrt{l-b} \sqrt{a-m} (l-m)}\,,
\eeqst
while
\beqst
&  &
  \partial _m\left((l-b)^{-\frac{1}{2}}(a-m)^{-\frac{1}{2}}
  F  \left(\frac{1}{2},\frac{1}{2};1;\frac{(l-a)(m-b)}{(l-b)(m-a)}\right)\right)  \\
 & = &
\frac{(b-l)F \left(-\frac{1}{2},\frac{1}{2};1;\frac{(l-a)(m-b)}{(l-b)(m-a)}\right)+(l-m) F \left(\frac{1}{2},\frac{1}{2};1;\frac{(l-a)(m-b)}{(l-b)(m-a)}\right)}{2\sqrt{l-b}
\sqrt{a-m} (b-m) (l-m)  }\,,
\eeqst
where, for the hypergeometric functions $F \left(\frac{1}{2},\frac{1}{2};1;z\right)$ and $F \left(-\frac{1}{2},\frac{1}{2};1;z\right)$ according to 
Sec.~2.1.3~\cite[v.1]{B-E} we have from the  Euler's formula
\[
F \left(\frac{1}{2},\frac{1}{2};1;z\right)= \frac{2}{\pi }\int_0^1 (1-t^2)^{-1/2} (1-zt^2)^{-1/2}dt  , \qquad 
F \left(-\frac{1}{2},\frac{1}{2};1;z\right)=  \frac{2}{\pi }\int_0^1 (1-t^2)^{-1/2} (1-zt^2)^{1/2}dt.
\]
These functions coincide with  the complete elliptic integrals of the first and second kind, $ {\bf K} (z) $ and ${\bf E} (z) $, respectively,
\beqst
&  &
{\bf K} (z) = \frac{\pi }{2}F \left(\frac{1}{2},\frac{1}{2},1,z^2\right) ,\qquad  {\bf E} (z) = \frac{\pi }{2}F \left(-\frac{1}{2},\frac{1}{2};1;z^2\right) \,.
\eeqst
(See (10) of \cite[v.1, Sec.~4.8, page 196]{B-E}   and   Sec.~13.8~\cite[v.2, page 317]{B-E}.) Then to calculate the second derivative we use    (21) of Sec.~2.8~\cite[v.1]{B-E} 
\[
 \frac{d}{d z}   F \left(-\frac{1}{2},\frac{1}{2};1;z\right)   
  =  
\frac{1 }{2z}F \left(-\frac{1}{2},\frac{1}{2};1;z\right)
- \frac{1 }{2z}F \left(\frac{1}{2},\frac{1}{2};1;z\right)\,,
\]
and obtain 
\beqst
&  &
\partial _l   \partial _m\left((l-b)^{-\frac{1}{2}}(a-m)^{-\frac{1}{2}}
  F  \left(\frac{1}{2},\frac{1}{2};1;\frac{(l-a)(m-b)}{(l-b)(m-a)}\right)\right)  \\
  &  = &
   \partial _l  \left(  \frac{(b-l)F \left(-\frac{1}{2},\frac{1}{2},1,\frac{(l-a)(m-b)}{(l-b)(m-a)}\right)+(l-m) F \left(\frac{1}{2},\frac{1}{2};1;\frac{(l-a)(m-b)}{(l-b)(m-a)}\right)}{2\sqrt{l-b}
\sqrt{a-m} (b-m) (l-m)  } \right) \\
 &  = &
 \frac{1}{4 (l-a) \sqrt{l-b} \sqrt{a-m} (b-m) (l-m)^2  } \\
 &  &
\times \Bigg[  \left(2 a b-a l-b l+l^2-(a+b) m+m^2\right)F \left(-\frac{1}{2},\frac{1}{2};1;\frac{(l-a)(m-b)}{(l-b)(m-a)}\right)\\
&  &
\hspace{3cm} +\,(l-m) (a-b-l+m) F \left(\frac{1}{2},\frac{1}{2};1;\frac{(l-a)(m-b)}{(l-b)(m-a)}\right)\Bigg]
\eeqst  
as well as
\beqst
&  &
  - \frac{1}{2(l-m)} \Big( \frac{\pa }{\pa l} -
\frac{\pa }{\pa m} \Big)   E(l,m;a,b) \\
& = & 
-\frac{1}{4 (a-l) \sqrt{l-b} \sqrt{a-m} (b-m) (l-m)^2   }  \\
&  &
\times \Bigg[\left((b-l) l+b m-m^2+a (-2 b+l+m)\right) F \left(-\frac{1}{2},\frac{1}{2};1;\frac{(l-a)(m-b)}{(l-b)(m-a)}\right)\\
&  &
\hspace{3cm} +(l-m) (-a+b+l-m) 
F \left(\frac{1}{2},\frac{1}{2};1;\frac{(l-a)(m-b)}{(l-b)(m-a)}\right) \Bigg]\,.
\eeqst
Hence (\ref{f5.1}) holds.
The lemma is proved. \hfill $\square$
\medskip

Consider operator 
\[
{\mathcal S}_{ch}^*  := \frac{\pa^2 \,}{\pa l \, \pa m} + \frac{1}{2(l-m)} \Big( \frac{\pa \, }{\pa l} -
\frac{\pa \,}{\pa m} \Big)  - \frac{1}{(l-m)^2} 
\]
that is the formally adjoint to the operator
\[
{\mathcal S}_{ch}  := \frac{\pa^2 \,}{\pa l \, \pa m} - \frac{1}{2(l-m)} \Big( \frac{\pa \, }{\pa l} -
\frac{\pa \,}{\pa m} \Big)  
\]
of the equation (\ref{f5.1}). The following lemma,   is an analog  of (2.4)\cite{GelfandII}.
\begin{lemma}
\label{L2.2}
If $v$ is a solution of the equation $ {\mathcal S}_{ch} ^* v = 0$, then $u = (l-m)^{-1}v$  
is a solution to ${\mathcal S}_{ch}  u=0$, and vice versa. 
\end{lemma}

\noindent
{\bf Proof.} Indeed, direct calculations lead to
\[
\frac{\pa \, }{\pa l} v 
  =    
(l-m)\frac{\pa \, }{\pa l}u +    u,\,\,   
\frac{\pa \, }{\pa m} v 
  =   
 (l-m)\frac{\pa \, }{\pa m}u - u,    \,\,
\frac{\pa^2 \, }{\pa l \, \pa m} v 
  =    
 (l-m)\frac{\pa ^2}{\pa m \,\pa l}u    +  \frac{\pa \, }{\pa m}   u  - \frac{\pa \, }{\pa l}u  \,.
\]
Then
\beqst
{\mathcal S}_{ch} ^* v 
& = & 
  (l-m)\frac{\pa ^2}{\pa m \,\pa l}u    +  \frac{\pa \, }{\pa m}   u  - \frac{\pa \, }{\pa l}u \\
&   &
+ \, \frac{1}{2(l-m)} \Big[ (l-m)\frac{\pa \, }{\pa l}u +    u -
 (l-m)\frac{\pa \, }{\pa m}u +  u \Big]   - \frac{1}{ (l-m)^2} (l-m)u \\
& = &
  (l-m)\Bigg\{ \frac{\pa ^2}{\pa m \,\pa l}u    -  \frac{1}{2(l-m)} \Big( \frac{\pa \, }{\pa l}u -  \frac{\pa \, }{\pa m}   u \Big)  \Bigg\} =0\,.
\eeqst 
Lemma is proved. \hfill $\square$

In the next lemma the Riemann function is presented.
\begin{proposition}
\label{PR}
The function 
\[
R(l,m;a,b) 
 = 
(l-m)E(l,m;a,b) 
 = 
(l-m)(l-b)^{-1/2} (a-m)^{-1/2} 
F\Big(\frac{1}{2},\frac{1}{2};1; \frac{(l-a)(m-b)}{(l-b)(m-a)} \Big)
\]
is the unique solution of the equation ${\mathcal S}_{ch} ^*v=0$ that satisfies the following conditions: \\
$(\mbox{\rm i})$ $\dsp{ R_l =  \frac{1}{2 (l-m)}R \quad }$ along the line $m=b$;\\
$(\mbox{\rm ii})$ $\dsp{ R_m =  - \frac{1}{2 (l-m)}R \quad }$ along the line $l=a$; \\
$(\mbox{\rm iii})$ $R(a,b;a,b)=  1$.
\end{proposition}
\medskip

\noindent
{\bf Proof.} It can be easily proven by the direct calculations. \hfill $\square$
\medskip

Next we use Riemann function $R(l,m;a,b)$ and function $ E(x,t;x_0,t_0)$ defined by (\ref{E}) to 
complete the proof of Theorem~\ref{T1}, which gives the  fundamental solution with a support 
in the forward cone $D_+ (x_0,t_0)$, $x_0 \in {\mb R}^n$, $t_0 \in {\mb R}$,
and  the  fundamental solution with a support in the backward cone $D_- (x_0,t_0)$, $x_0 \in {\mb R}^n$, $t_0 \in {\mb R}$,
defined by (\ref{D+}) and (\ref{D-}), respectively. 
\medskip

\noindent
{\bf Proof of Theorem~\ref{T1}.} 
We present a proof   for   $E_+(x,t;0,b) $ since for  
$E_-(x,t;0,b) $ it is similar.  
First, we note that the operator ${\mathcal S}$ is formally self-adjoint, ${\mathcal S}={\mathcal S}^*$. We must show that
\[
<E_+, {\mathcal S} \varphi > = \varphi (0,b)\,, \qquad 
\mbox{\rm for every} \quad \varphi \in C_0^\iy ({\mb R}^2)\,.
\]
Since $E(x,t;0,b) $ is locally integrable in ${\mb R}^{2}$, 
this is equivalent to showing that
\be
\label{Gel3.2}
\int \!\! \int_{{\mb R}^2} E_+(x,t;0,b) {\mathcal S} \varphi (x,t)\, dx \, dt = \varphi (0,b), \quad 
\mbox{\rm for every} \quad \varphi \in C_0^\iy ({\mb R}^2).
\ee
In the mean time $ {D(x,t)}/{D(l,m)} = (l-m)^{-1 } $ is the Jacobian of the transformation (\ref{Gel1.8}).
Hence the integral in the left-hand side of (\ref{Gel3.2}) is equal to
\beqst
&  &
\int \!\! \int_{{\mb R}^2} E_+(x,t;0,b) {\mathcal S} \varphi (x,t)\, dx \, dt  
  =   
\int_b^\iy  dt \int_{-(e^t-e^b)}^{e^t-e^b} E (x,t;0,b) {\mathcal S} \varphi (x,t)\, dx \\
& =  &
-  \int_{-\iy}^{-e^b}  
\int_{ e^b }^\iy E (l,m;e^b,-e^b)   (l-m)^{-1} \, dl \, dm  
  (l-m)^2
 \Bigg\{ \frac{\pa^2 }{\pa l \, \pa m} - \frac{1}{2 (l-m)} \Big( \frac{\pa }{\pa l} -
\frac{\pa }{\pa m} \Big)  
\Bigg\}  \varphi  .
\eeqst 
We will write $\varphi (l,m)$ for the function $\varphi (x,t)$ in the characteristic variables $l$, $m$. Then using the Riemann function $R$ constructed in Proposition~\ref{PR} we have
\beqst
&   &
\int \!\! \int_{{\mb R}^2} E_+(x,t;0,b) {\mathcal S} \varphi (x,t)\, dx \, dt \\
& = & 
-\int_{-\iy}^{-e^b} \!\! \int_{ e^b }^\iy R(l,m;e^b,-e^b)  
\Bigg\{ \frac{\pa^2 }{\pa l \, \pa m} - \frac{1}{2 (l-m)} \Big( \frac{\pa }{\pa l} -
\frac{\pa }{\pa m} \Big) \Bigg\}  \varphi 
\, dl \, dm \,.
\eeqst
Integrating by parts several times and applying Proposition \ref{PR}, we obtain
(\ref{Gel3.2}). Indeed,
\beqst
& & 
\int_{-\iy}^{-e^b} \!\! \int_{ e^b }^\iy R(l,m;e^b,-e^b)  
\Bigg\{ \frac{\pa^2 }{\pa l \, \pa m} - \frac{1}{2 (l-m)} \Big( \frac{\pa }{\pa l} -
\frac{\pa }{\pa m} \Big) \Bigg\}  \varphi 
\, dl \, dm \\
& = &
\int_{-\iy}^{-e^b} \!\!  \left[ R(l,m;e^b,-e^b)  
\frac{\pa  \varphi}{\pa m}  
 \right]_{l=e^b}^{l=\iy} dm - \int_{-\iy}^{-e^b} \!\! \int_{ e^b }^\iy  \left(\frac{\pa }{\pa l} R(l,m;e^b,-e^b) \right) 
\frac{\pa  \varphi }{\pa m}  
\, dl \, dm \\
&  &
- \int_{-\iy}^{-e^b} \!\! \int_{ e^b }^\iy R(l,m;e^b,-e^b)  
\frac{1}{2 (l-m)} \Big( \frac{\pa }{\pa l} -
\frac{\pa }{\pa m} \Big)  \varphi 
\, dl \, dm \,.
\eeqst
On the other hand, using the properties of Riemann function $R$ we obtain
\beqst
&  &
\int_{-\iy}^{-e^b} \!\!  \left[ R(l,m;e^b,-e^b)  
\frac{\pa  \varphi}{\pa m}  
 \right]_{l=e^b}^{l=\iy} dm \\
& = &
- \int_{-\iy}^{-e^b} \!\!   R(e^b,m;e^b,-e^b)  
\frac{\pa  }{\pa m}\varphi (e^b,m) \, dm  \\
& = &
- R(e^b,-e^b;e^b,-e^b)  
\varphi  \left. \right|_{l=e^b,\, m= -e^b}  +  \int_{-\iy}^{-e^b} \!\!   \left( \frac{\pa  }{\pa m} R(e^b,m;e^b,-e^b)  \right)
\varphi  (e^b,m) \, dm \\
& = &
- \left. \varphi  \right|_{l=e^b,\, m= -e^b} -  
\int_{-\iy}^{-e^b} \!\!   \frac{1}{2 (e^b-m)}  R(e^b,m;e^b,-e^b) 
\varphi  (e^b,m) \, dm\,,
\eeqst
while
\beqst
&  &
\int_{-\iy}^{-e^b} \!\! \int_{ e^b }^\iy  \left(\frac{\pa }{\pa l} R(l,m;e^b,-e^b) \right) 
\frac{\pa  \varphi }{\pa m}  
\, dl \, dm \\
& = &
\int_{ e^b }^\iy  \left(\frac{\pa }{\pa l} R(l,-e^b;e^b,-e^b) \right) 
 \varphi (l,  -e^b) \, dl - \int_{-\iy}^{-e^b} \!\! \int_{ e^b }^\iy  \left(\frac{\pa^2 }{\pa l\, \pa m} R(l,m;e^b,-e^b) \right) 
\varphi 
\, dl \, dm \\
& = &
\int_{ e^b }^\iy \frac{1}{2 (l-(-e^b))}  R(l,-e^b;e^b,-e^b) 
 \varphi (l,  -e^b) \, dl - \int_{-\iy}^{-e^b} \!\! \int_{ e^b }^\iy  \left(\frac{\pa^2 }{\pa l\, \pa m} R(l,m;e^b,-e^b) \right) 
\varphi 
\, dl \, dm \,.
\eeqst
Then
\beqst
&  &
\int_{-\iy}^{-e^b} \!\! \int_{ e^b }^\iy R(l,m;e^b,-e^b)  
\frac{1}{2 (l-m)} \frac{\pa \varphi}{\pa l}   
\, dl \, dm \\
& = & 
- \int_{-\iy}^{-e^b}  R(e^b,m;e^b,-e^b)  
\frac{1}{2 (e^b-m)}   \varphi (e^b,m) 
 \, dm - \int_{-\iy}^{-e^b} \!\! \int_{ e^b }^\iy  \left( \frac{\pa }{\pa l}\left( R(l,m;e^b,-e^b)  \frac{1}{2 (l-m)} \right)  \right)
 \varphi 
\, dl \, dm
\eeqst
and
\beqst
&  &
\int_{-\iy}^{-e^b} \!\! \int_{ e^b }^\iy R(l,m;e^b,-e^b)  
\frac{1}{2 (l-m)}  
\frac{\pa  \varphi }{\pa m}   
\, dl \, dm \\
&= &
\int_{ e^b }^\iy  R(l,-e^b;e^b,-e^b)  
\frac{1}{2 (l-(-e^b))}  \varphi  (l,-e^b)
\, dl  - \int_{-\iy}^{-e^b} \!\! \int_{ e^b }^\iy  \left(  \frac{\pa }{\pa m}  \left( R(l,m;e^b,-e^b)  
\frac{1}{2 (l-m)} \right) \right)  \varphi  
\, dl \, dm .
\eeqst
One more application of Proposition \ref{PR}   completes the proof of Theorem \ref{T1}.
\hfill $\square$

\section{Application to the Cauchy Problem: Source Term and $n=1$}
\label{S3}
\setcounter{equation}{0}

Consider now the Cauchy problem for the equation (\ref{Int_3}) 
with vanishing initial data (\ref{Int_4}).
For every  $(x,t) \in D_+ (0,b)$ one has $  -(e^t- e^b) \le x \le  e^t- e^b$, so that 
\[
 E(x,t;0,b)   =  
\frac{1}{\sqrt{  (e^t  + e^b)^2  -x^2  }}   F \left(\frac{1}{2},\frac{1}{2};1; \frac{    (e^t  - e^b)^2  -x^2  }
{  (e^t  + e^b)^2  -x^2        }\right) .
\]
The coefficient of the  equation (\ref{G-S}) is independent of $x$, therefore $E_+ (x,t;$ $y,b)$ $= E_+ (x-y,t;0,b)$.
Using the fundamental solution from Theorem \ref{T1} one can write the convolution 
\[
\hspace*{-0.5cm} u(x,t)
 = 
\int_{ -\iy}^{\iy} \! \int_{ -\iy}^{\iy} E_+ (x,t;y,b)f(y,b)\, db\, dy 
 = 
\int_{ 0}^{t}\! db \!\int_{ -\iy}^{\iy} E_+ (x-y,t;0,b) f(y,b) \,dy , 
\]
which is well-defined since supp$f \subset \{ t \geq 0 \} $. 
 Then according to the definition of the function $E_+ $ 
we obtain the statement of the Theorem~\ref{T1.1}. Thus, Theorem~\ref{T1.1} is proven.

\begin{remark} 
The argument of the hypergeometric function is nonnegative and bounded, 
\[
0 \leq  \frac{    (e^t  - e^b)^2  -z^2  
}
{  (e^t  + e^b)^2  -z^2        } < 1 \quad \mbox{  for all} \quad b\in (0,t),\,\, z \in ( e^b- e^t ,e^t  - e^b)\,.
 \]
The hypergeometric function $F\left(\frac{1}{2},\frac{1}{2};1; 
\frac{  (e^t  - e^b)^2  -z^2 }
{ (e^t  + e^b)^2  -z^2   } \right) $ at $b=t$ has a logarithmic singularity. Indeed, this follows  for $c=a+b\pm m$, ($m=0,1,2,\ldots$) 
from   formula 
15.3.10 of \cite[Ch.15]{A-S}:
\[
F  \left( a,b;a+b;z  \right) 
 = 
\frac{\Gamma (a+b)}{\Gamma (a)\Gamma (b)}\sum _{n=0}^\infty \frac{(a)_n(b)_n}{(n!)^2}
\left[ 2\psi (n+1) - \psi (a+n)- \psi (b+n) - \ln (1-z) \right] (1-z)^{n}\,,
\]
where $|\arg (1-z) |<\pi$ , $|1-z|<1 $. 
\end{remark}
 The following corollary is a manifestation of the  time-speed transformation principle introduced in \cite{YagTricomi}. It implies  
the existence of an operator transforming the solutions of the  Cauchy problem for the string equation
to the solutions of the  Cauchy problem for the inhomogeneous equation with time-dependent speed of propagation.
One may think of this transformation as  a ``two-stage'' Duhamel's principal,
but unlike  the last one, it reduces the equation with the time-dependent speed of propagation to 
the one with the speed of propagation independent of time. 

\begin{corollary}
\label{C2}
The solution $u= u (x,t)$ of the Cauchy problem (\ref{TricEq})-(\ref{4.2})  can be represented as follows
\[  
u(x,t) 
  =  
2 \int_{ 0}^{t} db \int_{ 0  }^{ e^t- e^b} dz \,  v(x,z;b)
\frac{1}{\sqrt{  (e^t  + e^b)^2  -z^2  }}   
F\left(\frac{1}{2},\frac{1}{2};1; 
\frac{  (e^t  - e^b)^2 -z^2}
{(e^t  + e^b)^2 -z^2} \right),
\] 
where the functions $v (x,t;\tau ) := \frac{1}{2} ( f(x+t,\tau ) +   f(x - t,\tau ))$,  $\tau \in [0,\iy)$,
form a one-parameter family  
of  solutions to the Cauchy problem for the string equation, that is,
\[
v_{tt} - v_{xx} =0\,, \qquad v (x,0;\tau ) = f(x,\tau )\,, \qquad v_t (x,0;\tau ) = 0\,.
\]
\end{corollary}
\medskip

\noindent
{\bf Proof.} From the theorem we have
\beqst 
u(x,t)  
& = &
 \int_{ 0}^{t} db \int_{  - (e^t- e^b) }^{ e^t- e^b} dz \, f(z+x,b)  
\frac{1}{\sqrt{  (e^t  + e^b)^2  -z^2  }}  
F\left(\frac{1}{2},\frac{1}{2};1; 
\frac{  (e^t  - e^b)^2 -z^2}
{(e^t  + e^b)^2 -z^2} \right) \\
& = &
\int_{ 0}^{t} db \int_{ 0  }^{ e^t- e^b} dz \, f(z+x,b)  
\frac{1}{\sqrt{  (e^t  + e^b)^2  -z^2  }}   
F\left(\frac{1}{2},\frac{1}{2};1; 
\frac{  (e^t  - e^b)^2 -z^2}
{(e^t  + e^b)^2 -z^2} \right) \\
&   &
+ \int_{ 0}^{t} db \int_{     0}^{ e^t- e^b} dz \, f(-z+x,b)  
\frac{1}{\sqrt{  (e^t  + e^b)^2  -z^2  }}  
F\left(\frac{1}{2},\frac{1}{2};1; 
\frac{  (e^t  - e^b)^2 -z^2}
{(e^t  + e^b)^2 -z^2} \right) \\
& = &
2\int_{ 0}^{t} db \int_{ 0  }^{ e^t- e^b} dz \, \frac{1}{2} \{f(x+z ,b) +  f(x-z ,b)\}
\frac{1}{\sqrt{  (e^t  + e^b)^2  -z^2  }}   
F\left(\frac{1}{2},\frac{1}{2};1; 
\frac{  (e^t  - e^b)^2 -z^2}
{(e^t  + e^b)^2 -z^2} \right) \,.
\eeqst
The corollary is proven. \hfill $\square$

\section{Some Properties of the Function  $E(x,t;y,b)$}
\label{S4}

\setcounter{equation}{0}

For $b \in {\mb R}$   the function $E(x,t;y,b)$ in the domain  $D_+ (y,b) \cup D_- (y,b)$ is defined by  (\ref{E}),
where $F\big(a, b;c; \zeta \big) $ is the hypergeometric function. In this section we collect some elementary auxiliary formulas 
to make proofs of the main theorems more transparent. For the simplicity we consider case $n=1$ in detail. The case of $n>1$ is
very similar.
\begin{proposition}
\label{P_E}
One has
\beq
\label{E_1}
\label{E_2}
E(x,t;y,b) 
 =   
E(x-y,t;0,b)  & , & 
E(x,t;0,b) 
  =   
E(-x,t;0,b),\\
\label{E_2a}
E(x,t;0,\ln(e^t-x)) 
& =  &
\frac{1}{2}\frac{1}{\sqrt{  e^t }  \sqrt{e^t-x}}, \\
\label{E_3a}
\frac{\pa }{\pa b} \Big(e^{ b}E( e^b  -e^t,t;0,b)  \Big)
& =  & 
\frac{1}{4}e^{-t/2}e^{b/2}, \\
\hspace{-0.7cm} \frac{\pa }{\pa b} \Big( b e^{ b}E( e^b  -e^t,t;0,b)  \Big)
 =  
\frac{\pa }{\pa b} \Big( b e^{ b}E(e^t- e^b  ,t;0,b) \Big)  
\label{E_4}\label{E_3}
& =  & 
 \frac{\pa }{\pa b} \Big( b e^{ b} \frac{1}{2}e^{-t/2}e^{-b/2}   \Big) =   \frac{1}{4}e^{-t/2}e^{b/2}(2+b), \\
\label{E_5}
\lim_{y\to x+ e^t- e^b} \frac{\pa }{\pa x} E(x-y,t;0,b) 
& =  & 
\frac{1}{16} e^{-2(b+t)} e^{b/2}e^{t/2}( e^{b }- e^{t }),\\
\label{E_6}
\lim_{y\to x- e^t+ e^b} \frac{\pa }{\pa x} E(x-y,t;0,b) 
& =  & 
\frac{1}{16} e^{-2(b+t)} e^{b/2}e^{t/2}( -e^{b }+ e^{t }), \\
\label{E_7}
\left[ \frac{\pa }{\pa b} E(x,t;0,b) \right]_{b= \ln(e^t-x)}
& =  & 
\frac{e^{-2t}\sqrt{  e^t } (-4e^t+x)}
{16 \sqrt{e^t-x}} \,,
\eeq
\vspace{-0.5cm} 
\beq
\frac{\pa E}{\pa b}  (z,t;0,0) 
    =      
\frac{1}{ 2  ((e^t-1)^2-z^2 ) \sqrt{(1+e^t)^2-z^2}}  
 \left\{ (1-e^{2 t}+z^2 ) 
F \left(-\frac{1}{2},\frac{1}{2};1; \frac{  (e^t-1)^2-z^2 }{  (e^t+1)^2-z^2 }   \right) \right.  \nonumber \\ 
\label{E_bzt00}
\left. +   2  (e^t -1) F \left(\frac{1}{2},\frac{1}{2};1;  \frac{  (e^t-1)^2-z^2 }{  (e^t+1)^2-z^2 }   \right)  \right\}.
\eeq  
\end{proposition}
\medskip

\noindent
{\bf Proof.} The properties (\ref{E_1})   and (\ref{E_2a}) are evident. 
To prove (\ref{E_3a}) and (\ref{E_3}) we write
\be   
 \label{4.11}  E(e^b  -e^t,t;0,b)
  =  
 (2 e^b   )^{-\frac{1}{2}} 
 (2e^t)^{-\frac{1}{2}}   F \left(\frac{1}{2}, \frac{1}{2};1; 0 \right) = \frac{1}{2} e^{-\frac{b}{2}}e^{-\frac{t}{2}}\,,
\ee  
that implies (\ref{E_3a}) and (\ref{E_4}).
To prove (\ref{E_5}) we denote 
\beqst
z:= \frac{(e^t  - e^{b})^2 - (x-y )^2 }
{(e^t  + e^{b})^2 - (x-y )^2 } \, ,
\eeqst
and obtain
\beq  
\label{4.1Ex}
\frac{\pa }{\pa x} E(x-y,t;0,b) 
& = &
 -\frac{1 }{2} (x-y+ e^t+ e^{b}   )^{-\frac{3}{2}} 
 (-x+y     + e^t + e^{b} )^{-\frac{1}{2}}  F \left(\frac{1}{2}, \frac{1}{2};1; z \right) \nonumber \\
&  &
+ \frac{1 }{2} (x-y+ e^t+ e^{b}   )^{-\frac{1}{2}} 
 (-x+y     + e^t + e^{b} )^{-\frac{3}{2}}  F \left(\frac{1}{2}, \frac{1}{2};1;z\right) \nonumber \\
&   &
+ \, \big( (e^t+ e^{b})^2 - (x-y)^2 \big)^{-\frac{1}{2}} F'_z \left(\frac{1}{2}, \frac{1}{2};1; z\right) \frac{\pa }{\pa x} z \,.
\eeq 
It is easily seen that
\beqst
\frac{\pa }{\pa x} z 
& = &
-\frac{ 8(x-y) e^{t + b} }
{[(x-y)^2- (e^t  + e^{b})^2]^2}   \,.
\eeqst
Here
\[
\lim_{y\to x+ e^t- e^b} z= 0, \qquad  
 \lim_{y\to x+ e^t- e^b} \frac{\pa }{\pa x} z =   \frac{ 1}{ 2  }( e^{ -b} - e^{- t }),
\]
while according to (23)  \cite[Sec.2.8 v.1]{B-E} we have
\beqst
 \partial _z   F \left(\frac{1}{2},\frac{1}{2},1,z\right)   
& = &
\frac{1 }{2z (1 -z) }F \left(-\frac{1}{2},\frac{1}{2},1,z\right)
- \frac{1 }{2z}F \left(\frac{1}{2},\frac{1}{2},1,z\right)\,.
\eeqst
Consequently,
\beqst
  \lim_{y\to x+ e^t- e^b}\partial _z   F \left(\frac{1}{2},\frac{1}{2},1,z\right)  
& = &
 \lim_{z\to 0} \frac{1 }{2z}\left\{ \frac{1 }{1 -z}F \left(-\frac{1}{2},\frac{1}{2},1,z\right)
- F \left(\frac{1}{2},\frac{1}{2},1,z\right) \right\}\,.
\eeqst
In fact (See, e.g.\cite{B-E}.),
\be
\label{4.12}\label{4.13}
 F\Big(\frac{1}{2},\frac{1}{2};1; z  \Big) = 1 + \frac{1}{4}z + O(z^2) \quad \mbox{\rm and } \quad
 F\Big(-\frac{1}{2},\frac{1}{2};1; z  \Big) = 1- \frac{1}{4}z + O(z^2)   \quad \mbox{\rm as } \quad
z\rightarrow 0   
\ee
imply
\be  
\label{Fderiv}
  \lim_{y\to x+ e^t- e^b}\partial _z   F \left(\frac{1}{2},\frac{1}{2},1,z\right) 
  =   
 \lim_{z\to 0} \frac{1 }{2z}\left\{ \frac{1 }{1 -z}\left( 1- \frac{1}{4}z + O(z^2)  \right)
-  \left( 1 + \frac{1}{4}z + O(z^2) \right) \right\} 
      =  
\frac{1 }{4}\,.
\ee  
Thus, according to (\ref{4.1Ex}) we obtain
\beqst 
 \lim_{y\to x+ e^t- e^b}\frac{\pa }{\pa x} E(x-y,t;0,b) & = &
  \lim_{y\to x+ e^t- e^b} -\frac{1 }{2} (x-y+ e^t+ e^{b}   )^{-\frac{3}{2}} 
 (-x+y     + e^t + e^{b} )^{-\frac{1}{2}} \\
&  &
+   \lim_{y\to x+ e^t- e^b}\frac{1 }{2} (x-y+ e^t+ e^{b}   )^{-\frac{1}{2}} 
 (-x+y     + e^t + e^{b} )^{-\frac{3}{2}}   \\
&   &
+ \,  \lim_{y\to x+ e^t- e^b} (x-y+ e^t+ e^{b}   )^{-\frac{1}{2}} 
 (-x+y     + e^t + e^{b} )^{-\frac{1}{2}} \frac{1 }{4} \cdot \frac{ 1}{ 2  }( e^{ -b} - e^{- t }) \\
& = &
 -\frac{1 }{2} (-e^t+ e^{b}+ e^t+ e^{b}   )^{-\frac{3}{2}} 
 (e^t- e^{b}     + e^t + e^{b} )^{-\frac{1}{2}} \\
&  &
+  \frac{1 }{2} (-e^t+ e^{b}+ e^t+ e^{b}   )^{-\frac{1}{2}} 
 (e^t- e^{b}    + e^t + e^{b} )^{-\frac{3}{2}}   \\
&   &
+ (-e^t+ e^{b}+ e^t+ e^{b}   )^{-\frac{1}{2}} 
 (e^t- e^{b}    + e^t + e^{b} )^{-\frac{1}{2}} \frac{1 }{8} ( e^{ -b} - e^{- t })\\
& = &
 -\frac{1 }{2} ( 2e^{b})^{-\frac{3}{2}} 
 (2e^t )^{-\frac{1}{2}} 
+  \frac{1 }{2} ( 2e^{b} )^{-\frac{1}{2}} 
 (2e^t)^{-\frac{3}{2}}   
+ (2e^{b} )^{-\frac{1}{2}} 
 (2e^t)^{-\frac{1}{2}} \frac{1 }{8} ( e^{ -b} - e^{- t })\\
& = &
\frac{1 }{16} e^{-2(b+t)}  e^{ \frac{t}{2}}  e^{ \frac{b}{2}} ( e^{ b} - e^{t })\,.
\eeqst
To prove (\ref{E_7}) we write
\beqst 
\frac{\pa }{\pa b} E(x,t;0,b) 
& = &
 \left(  \frac{\pa }{\pa b}\big( (e^t+ e^{b})^2 -  x^2 \big)^{-\frac{1}{2}}   \right) F \left(\frac{1}{2}, \frac{1}{2};1; \frac{ x^2- (e^t  - e^{b} )^2  }
{ x^2- (e^t  + e^{b} )^2 } \right)\\
&   &
+\, \big( (e^t+ e^{b})^2 -  x^2 \big)^{-\frac{1}{2}} \frac{\pa }{\pa b} \left(   F \left(\frac{1}{2}, \frac{1}{2};1; \frac{ x^2- (e^t  - e^{b} )^2  }
{ x^2- (e^t  + e^{b} )^2 }  \right) \right)\\
& = &
- e^{b}(e^t+ e^{b})\big( (e^t+ e^{b})^2 -  x^2 \big)^{-\frac{3}{2}}  F \left(\frac{1}{2}, \frac{1}{2};1; \frac{ x^2- (e^t  - e^{b} )^2  }
{ x^2- (e^t  + e^{b} )^2 } \right)\\
&   &
+\, \big( (e^t+ e^{b})^2 -  x^2 \big)^{-\frac{1}{2}}  F'_z \left(\frac{1}{2}, \frac{1}{2};1; \frac{ x^2- (e^t  - e^{b} )^2  }
{ x^2- (e^t  + e^{b} )^2 }  \right) \frac{\pa }{\pa b} \frac{(e^t  - e^{b} )^2- x^2}
{(e^t  + e^{b} )^2- x^2}\\
& = &
- e^{b}(e^t+ e^{b})\big( (e^t+ e^{b})^2 -  x^2 \big)^{-\frac{3}{2}}  F \left(\frac{1}{2}, \frac{1}{2};1; \frac{ x^2- (e^t  - e^{b} )^2  }
{ x^2- (e^t  + e^{b} )^2 } \right)\\
&   &
+\, \big( (e^t+ e^{b})^2 -  x^2 \big)^{-\frac{1}{2}}  F'_z \left(\frac{1}{2}, \frac{1}{2};1; \frac{ x^2- (e^t  - e^{b} )^2  }
{ x^2- (e^t  + e^{b} )^2 } \right) \\
&  &
\times  \frac{-2e^{b}(e^t  - e^{b} )  [(e^t  + e^{b} )^2- x^2]- [(e^t  - e^{b} )^2- x^2]2e^{b}(e^t  + e^{b} )}
{[(e^t  + e^{b} )^2- x^2]^2}\\
& = &
- e^{b}(e^t+ e^{b})\big( (e^t+ e^{b})^2 -  x^2 \big)^{-\frac{3}{2}}  F \left(\frac{1}{2}, \frac{1}{2};1; \frac{ x^2- (e^t  - e^{b} )^2  }
{ x^2- (e^t  + e^{b} )^2 }  \right)\\
&   &
+\,   F'_z \left(\frac{1}{2}, \frac{1}{2};1; \frac{ x^2- (e^t  - e^{b} )^2  }
{ x^2- (e^t  + e^{b} )^2 } \right)  \frac{4e^{b} e^t x^2 -4e^{b} e^t  (e^{2t}  - e^{2b} ) }
{[(e^t  + e^{b} )^2- x^2]^2\sqrt{(e^t+ e^{b})^2 -  x^2}}\,.
\eeqst
On the other hand (\ref{Fderiv}) implies
\beq  
\label{3.13a}
\left[ \frac{\pa }{\pa b} E(x,t;0,b) \right]_{b= \ln(e^t-x)}
& = &
- \left[ e^{b}(e^t+ e^{b})\big( (e^t+ e^{b})^2 -  x^2 \big)^{-\frac{3}{2}}  \right]_{b= \ln(e^t-x)} F \left(\frac{1}{2}, \frac{1}{2};1; 0 \right)
\nonumber \\
&   &
+\,   F'_z \left(\frac{1}{2}, \frac{1}{2};1; 0 \right) \left[  \frac{4e^{b} e^t x^2 -4e^{b} e^t  (e^{2t}  - e^{2b} ) }
{[(e^t  + e^{b} )^2- x^2]^2\sqrt{(e^t+ e^{b})^2 -  x^2}} \right]_{b= \ln(e^t-x)} \nonumber \\
& = &
- \left[ e^{b}(e^t+ e^{b})\big( (e^t+ e^{b})^2 -  x^2 \big)^{-\frac{3}{2}}  \right]_{b= \ln(e^t-x)} \nonumber \\
&   &
+\,  \frac{1}{4} \left[  \frac{4e^{b} e^t x^2 -4e^{b} e^t  (e^{2t}  - e^{2b} ) }
{[(e^t  + e^{b} )^2- x^2]^2\sqrt{(e^t+ e^{b})^2 -  x^2}} \right]_{b= \ln(e^t-x)}\,.
\eeq 
Then
\be 
\label{3.13}
\left[ e^{b}(e^t+ e^{b})\big( (e^t+ e^{b})^2 -  x^2 \big)^{-\frac{3}{2}}  \right]_{b= \ln(e^t-x)}  
  =  
(e^t-x)(2e^t -x) \big( 4 e^t  (e^t-x)\big)^{-\frac{3}{2}}   
  =  
\frac{e^{-2t}\sqrt{  e^t }(2e^t -x) }{ 8  \sqrt{e^t-x}} 
\ee 
and
\beq
\left[  \frac{e^{b} e^t x^2 -e^{b} e^t  (e^{2t}  - e^{2b} ) }
{[(e^t  + e^{b} )^2- x^2]^2\sqrt{(e^t+ e^{b})^2 -  x^2}} \right]_{b= \ln(e^t-x)} & = &
\label{3.14}
- \frac{e^{-2t}\sqrt{  e^t } x}
{16 \sqrt{e^t-x}} \,.
\eeq
Hence (\ref{3.13a}), (\ref{3.13}), and  (\ref{3.14}) prove (\ref{E_7}).
\medskip

To prove (\ref{E_bzt00})   we use (\ref{E_7})    and  
\beqst
\hspace*{-0.5cm}\frac{\pa }{\pa b}  E(z,t;0,b)
& = &
- e^{b}(e^t + e^{b}) (   (e^t + e^{b})^2 -z^2 )^{-\frac{3}{2}}F \left(\frac{1}{2}, \frac{1}{2};1; \frac{(e^t  - e^{b})^2- z^2}
{(e^t  + e^{b})^2- z^2} \right) \\
&  &
+   
 (   (e^t + e^{b})^2 -z^2)^{-\frac{1}{2}} \frac{\pa }{\pa b}F \left(\frac{1}{2}, \frac{1}{2};1; \frac{(e^t  - e^{b})^2- z^2}
{(e^t  + e^{b})^2- z^2} \right) \,.
\eeqst
If we denote
\[
\zeta = \frac{(e^t  - e^{b})^2- z^2}
{(e^t  + e^{b})^2- z^2} \,, \qquad \zeta_0 = \frac{(e^t  - 1)^2- z^2}{(e^t  + 1)^2- z^2} \,,
\]
then
\[
 \frac{\pa }{\pa b}\zeta 
  =  
\frac{  4e^{ b}e^t(e^{2b}    - e^{2t}) +  4z^2e^{b}e^t }
{[(e^t  + e^{b})^2- z^2]^2} \,,\qquad \frac{\pa \zeta}{\pa b} \Big|_{b=0}  
  =    
\frac{  4 e^t(1    - e^{2t}) +  4z^2 e^t }
{[(e^t  + 1)^2- z^2]^2} \,.
\]
Hence
\beqst
\hspace*{-0.5cm}\frac{\pa E}{\pa b}  (z,t;0,0)
& = &
-  (e^t + 1) [   (e^t + 1)^2 -z^2 ]^{-\frac{3}{2}}F \left(\frac{1}{2}, \frac{1}{2};1;  \zeta_0 \right) \\
&  &
+   
 [   (e^t + 1)^2 -z^2]^{-\frac{1}{2}}  F_\zeta  \left(\frac{1}{2}, \frac{1}{2};1;  \zeta_0 \right) \frac{  4 e^t(1    - e^{2t}) +  4z^2 e^t }
{[(e^t  + 1)^2- z^2]^2} \,. 
\eeqst
According to (\ref{Fderivative}) we obtain
\beqst
&  &
\frac{\pa E}{\pa b}  (z,t;0,0) \\
& = &
-  (e^t + 1) [   (e^t + 1)^2 -z^2 ]^{-\frac{3}{2}}F \left(\frac{1}{2}, \frac{1}{2};1;  \zeta_0 \right) \\
&  &
+   
 [   (e^t + 1)^2 -z^2]^{-\frac{1}{2}} \left[  \frac{1 }{2\zeta_0   (1 -\zeta_0  ) }F \left(-\frac{1}{2},\frac{1}{2};1;\zeta_0 \right)
- \frac{1 }{2\zeta_0 }F \left(\frac{1}{2},\frac{1}{2};1;\zeta_0 \right)  \right] \frac{  4 e^t(1    - e^{2t}) +  4z^2 e^t }
{[(e^t  + 1)^2- z^2]^2}  \\
& = &
-  (e^t + 1) [   (e^t + 1)^2 -z^2 ]^{-\frac{3}{2}}F \left(\frac{1}{2}, \frac{1}{2};1;  \zeta_0 \right) \\
&  &
+   
 [   (e^t + 1)^2 -z^2]^{-\frac{3}{2}} \frac{2  [ e^t(1    - e^{2t}) +   z^2 e^t ]}{  (e^t  - 1)^2- z^2  } \left[  \frac{1 }{  1 -\zeta_0  }F \left(-\frac{1}{2},\frac{1}{2};1;\zeta_0 \right)
- F \left(\frac{1}{2},\frac{1}{2};1;\zeta_0 \right)  \right]\,.
\eeqst
The term with $F \left(\frac{1}{2},\frac{1}{2};1;\zeta_0 \right) $ contains a factor 
\beqst
&  &
-  (e^t + 1) [   (e^t + 1)^2 -z^2 ]^{-\frac{3}{2}}  
-  
 [   (e^t + 1)^2 -z^2]^{-\frac{3}{2}} \frac{2  [ e^t(1    - e^{2t}) +   z^2 e^t ]}{  (e^t  - 1)^2- z^2  }     \\
& = &
- \frac{[   (e^t + 1)^2 -z^2 ]^{-1}}{[ (e^t  - 1)^2- z^2]\sqrt{ (e^t + 1)^2 -z^2} }  \Bigg[ (e^t + 1)  [ (e^t  - 1)^2- z^2]  
+   2  [ e^t(1    - e^{2t}) +   z^2 e^t ]   \Bigg] \,,
\eeqst
where
\beqst 
 (e^t + 1)  [ (e^t  - 1)^2- z^2]  
+   2  [ e^t(1    - e^{2t}) +   z^2 e^t ]  
& = &
(-e^t +1) [  (e^t +1)^2         -  z^2 ]  \,.
\eeqst
The coefficient of $ F \left(-\frac{1}{2},\frac{1}{2};1;\zeta_0 \right)$ is 
\beqst
 [   (e^t + 1)^2 -z^2]^{-\frac{3}{2}} \frac{2  [ e^t(1    - e^{2t}) +   z^2 e^t ]}{  (e^t  - 1)^2- z^2  }  \frac{1 }{  1 -\frac{(e^t  - 1)^2- z^2}{(e^t  + 1)^2- z^2}   }  = 
 \frac{1 }{   2     }  \frac{1    - e^{2t} +   z^2 }{ [ (e^t  - 1)^2- z^2 ] \sqrt{(e^t + 1)^2 -z^2} } \,.
\eeqst
The formula (\ref{E_bzt00}) and, consequently, the proposition are proven. \hfill $\square$

\section{The Cauchy Problem:  Second Datum and $n=1$}
\label{S5}

\setcounter{equation}{0}

In this section we prove Theorem~\ref{T1.3} in the case of  $\varphi_0  (x)=0$. More precisely, we have to prove that 
the solution $u (x,t)$ of the Cauchy problem  (\ref{oneDphy01}) with  $\varphi_0  (x)=0$ and  $\varphi_1  (x)=\varphi (x)$
can be represented as follows
\be 
\label{3.1n}
u(x,t)  
 =   
   \int_{0}^{  e^t-1} \,\Big[      \varphi   (x+ z)  +   \varphi   (x - z)    \Big] K_1(z,t) dz
 =   
   \int_{0}^{  1} \,\Big[      \varphi   (x+ \phi (t) s )  +   \varphi   (x - \phi (t) s)    \Big] K_1(\phi (t) s,t) \phi (t) ds ,
\ee
where $\phi (t) =  e^t-1$. 
The proof of the theorem is splitted into several steps.

\begin{proposition}
\label{C3a}
The solution $u = u (x,t)$ of the Cauchy problem    (\ref{oneDphy01}) with  $\varphi_0  (x)=0$ and  $\varphi_1  (x)=\varphi (x)$
can be represented as follows
\beqst
u(x,t)  
& = &  
  \int_{ 0}^{t} \, db \Big[\frac{1}{4}e^{-t/2}e^{b/2}(2+b) +  \frac{1}{16}  b e^{-3t/2} e^{b/2} ( e^{b }- e^{t }) \Big] \Big[
     \varphi   (x+ e^t- e^b)     
+   \varphi   (x - e^t+ e^b)    \Big] \\
&  &
+ \int_{ 0}^{t} b e^{2b}\, db \int_{ x - (e^t- e^b)}^{x+ e^t- e^b} dy \, \varphi   (y)  
\Big( \frac{\pa }{\pa y} \Big)^2 E(x-y,t;0,b) \,.
\eeqst 
\end{proposition}
 \medskip

 \noindent
 {\bf Proof.} We look for the solution $u=u(x,t)$ of the form $
u(x,t)=w(x,t)+ t  \varphi (x)$. 
Then $ u_{tt} - e^{2t}u_{xx} =0$ implies 
\beqst
&  &
w_{tt} - e^{2t}w_{xx} = t e^{2t}\varphi^{(2)} (x), \qquad  w(x,0) = 0,\quad w_t(x,0)= 0 \,.
\eeqst
We set $f(x,t)=  t e^{2t}\varphi^{(2)} (x) $ and due to Theorem~\ref{T1.1} obtain  
\[
w(x,t)   
  =  
 \int_{ 0}^{t} b e^{2b}\, db \int_{ x - (e^t- e^b)}^{x+ e^t- e^b} dy \, \varphi^{(2)} (y)  
E(x-y,t;0,b)  \,.
\]
Then we integrate by parts:
\beqst
w(x,t) 
& = &
 \int_{ 0}^{t} b e^{2b}\, db  \Bigg(
 \varphi^{(1)} (x+ e^t- e^b)  E(  -e^t+ e^b ,t;0,b)   
-  \varphi^{(1)} (x - e^t+ e^b)  E(e^t- e^b  ,t;0,b)\Bigg) \\
&  &
- \int_{ 0}^{t} b e^{2b}\, db \int_{ x - (e^t- e^b)}^{x+ e^t- e^b} dy \, \varphi^{(1)} (y)  
\frac{\pa }{\pa y} E(x-y,t;0,b)\,.
 \eeqst
But
\[
 \varphi^{(1)} (x+ e^t- e^b) = -  e^{-b} \frac{\pa }{\pa b}  \varphi  (x+ e^t- e^b) \quad \mbox{\rm and} \quad
 \varphi^{(1)} (x - e^t+ e^b)=    e^{-b} \frac{\pa }{\pa b}  \varphi  (x - e^t+ e^b)
\]
by one more integration by parts imply
\beqst
&  &
w(x,t)  \\
& = &
\Bigg[b e^{b} \Bigg(
-      \varphi  (x+ e^t- e^b)  E(  -e^t+ e^b ,t;0,b)   
-       \varphi  (x - e^t+ e^b)  E(e^t- e^b  ,t;0,b)\Bigg)\Bigg]_{ b=0}^{b=t}  \\
&  &
- \int_{ 0}^{t} \, db  \Bigg(
-     \varphi  (x+ e^t- e^b)   \frac{\pa }{\pa b} \Big( b e^{ b}E(  -e^t+ e^b ,t;0,b)  \Big) 
-    \varphi  (x - e^t+ e^b)   \frac{\pa }{\pa b} \Big(b e^{ b}E(e^t- e^b  ,t;0,b)\Big) \Bigg) \\
&  &
- \int_{ 0}^{t} b e^{2b}\, db \int_{ x - (e^t- e^b)}^{x+ e^t- e^b} dy \, \varphi^{(1)} (y)  
\frac{\pa }{\pa y} E(x-y,t;0,b)\\
& = &  
-  2 t e^{t}  
   \varphi  (x )  E( 0 ,t;0,t)     \\
&  &
- \int_{ 0}^{t} \, db  \Bigg(
-     \varphi  (x+ e^t- e^b)   \frac{\pa }{\pa b} \Big( b e^{ b}E(  -e^t+ e^b ,t;0,b)  \Big) 
-    \varphi  (x - e^t+ e^b)   \frac{\pa }{\pa b} \Big(b e^{ b}E(e^t- e^b  ,t;0,b)\Big) \Bigg) \\
&  &
- \int_{ 0}^{t} b e^{2b}\, db \int_{ x - (e^t- e^b)}^{x+ e^t- e^b} dy \, \varphi^{(1)} (y)  
\frac{\pa }{\pa y} E(x-y,t;0,b)\,.
 \eeqst
Since  \, $ E( 0 ,t;0,t)=e^{-t}/2 $ \,  we obtain 
\beqst
&  &
u(x,t)   \\
& = &  
- \int_{ 0}^{t} \, db  \Bigg(
-     \varphi  (x+ e^t- e^b)   \frac{\pa }{\pa b} \Big( b e^{ b}E(  -e^t+ e^b ,t;0,b)  \Big) 
-    \varphi  (x - e^t+ e^b)   \frac{\pa }{\pa b} \Big(b e^{ b}E(e^t- e^b  ,t;0,b)\Big) \Bigg) \\
&  &
- \int_{ 0}^{t} b e^{2b}\, db \int_{ x - (e^t- e^b)}^{x+ e^t- e^b} dy \, \varphi^{(1)} (y)  
\frac{\pa }{\pa y} E(x-y,t;0,b)\,.
\eeqst
Then we apply (\ref{E_3})  of Proposition \ref{P_E} to derive the next representation
\beqst
u(x,t)   
& = &  
  \int_{ 0}^{t} \, db \frac{1}{4} e^{-t/2}e^{b/2}(2+b)  \Big(
     \varphi  (x+ e^t- e^b)     
+   \varphi  (x - e^t+ e^b)    \Big) \\
&  &
- \int_{ 0}^{t} b e^{2b}\, db \int_{ x - (e^t- e^b)}^{x+ e^t- e^b} dy \, \varphi^{(1)} (y)  
\frac{\pa }{\pa y} E(x-y,t;0,b)\,.
\eeqst
The integration by parts and \, $\frac{\pa }{\pa y} E(x-y,t;0,b)=-\frac{\pa }{\pa x} E(x-y,t;0,b)  $ \, imply
\beqst
u(x,t)  
& = &  
  \int_{ 0}^{t} \, db \frac{1}{4}e^{-t/2}e^{b/2}(2+b)  \Big(
     \varphi  (x+ e^t- e^b)     
+   \varphi  (x - e^t+ e^b)    \Big) \\
&  &
+ \int_{ 0}^{t} b e^{2b}\, db  \,  \varphi  (x+ e^t- e^b)  
 \Big[ \frac{\pa }{\pa x} E(x-y,t;0,b) \Big]_{y= x+ e^t- e^b}  \\
&  &
- \int_{ 0}^{t} b e^{2b}\, db  \,  \varphi  (x - (e^t- e^b))  
 \Big[ \frac{\pa }{\pa x} E(x-y,t;0,b) \Big]_{ y= x - (e^t- e^b)}  \\
&  &
+ \int_{ 0}^{t} b e^{2b}\, db \int_{ x - (e^t- e^b)}^{x+ e^t- e^b} dy \, \varphi  (y)  
\Big( \frac{\pa }{\pa y} \Big)^2 E(x-y,t;0,b)  \,.
\eeqst
The application of  (\ref{E_5}) and  (\ref{E_6}) from Proposition \ref{P_E} leads to
\beqst
u(x,t)  
& = &  
  \int_{ 0}^{t} \, db \frac{1}{4}e^{-t/2}e^{b/2}(2+b)  \Big[
     \varphi  (x+ e^t- e^b)     
+   \varphi  (x - e^t+ e^b)    \Big] \\
&  &
+ \int_{ 0}^{t} b e^{2b}\, db  \,  \frac{1}{16} e^{-2(b+t)} e^{b/2}e^{t/2}( e^{b }- e^{t }) \Big[ \varphi  (x+ e^t- e^b)  
+ \varphi  (x - (e^t- e^b))   \Big] \\
&  &
+ \int_{ 0}^{t} b e^{2b}\, db \int_{ x - (e^t- e^b)}^{x+ e^t- e^b} dy \, \varphi  (y)  
\Big( \frac{\pa }{\pa y} \Big)^2 E(x-y,t;0,b)  \\ 
& = &  
  \int_{ 0}^{t} \, db \frac{1}{4} e^{-t/2}e^{b/2}(2+b)  \Big[
     \varphi  (x+ e^t- e^b)     
+   \varphi  (x - e^t+ e^b)    \Big] \\
&  &
+ \int_{ 0}^{t}\, db  \,  \frac{1}{16}  b e^{-3t/2} e^{b/2} ( e^{b }- e^{t }) \Big[ \varphi  (x+ e^t- e^b)  
+ \varphi  (x -  e^t+ e^b)    \Big] \\
&  &
+ \int_{ 0}^{t} b e^{2b}\, db \int_{ x - (e^t- e^b)}^{x+ e^t- e^b} dy \, \varphi  (y)  
\Big( \frac{\pa }{\pa y} \Big)^2 E(x-y,t;0,b) \,.
\eeqst  
Finally,
\beq 
\label{5.4}
u(x,t)  
& = &  
  \int_{ 0}^{t} \, db \Big[\frac{1}{4}e^{-t/2}e^{b/2}(2+b) +  \frac{1}{16}  b e^{-3t/2} e^{b/2} ( e^{b }- e^{t }) \Big] \Big[
     \varphi  (x+ e^t- e^b)     
+   \varphi  (x - e^t+ e^b)    \Big] \nonumber \\
&  &
+ \int_{ 0}^{t} b e^{2b}\, db \int_{ x - (e^t- e^b)}^{x+ e^t- e^b} dy \, \varphi  (y)  
\Big( \frac{\pa }{\pa y} \Big)^2 E(x-y,t;0,b) \,.
\eeq   
To get last representation we have used (\ref{E_2}) and (\ref{4.11}).
The proposition is proven.
\hfill $\square$

\begin{corollary}
\label{C5.3}
The solution $u= u (x,t)$ of the Cauchy problem  (\ref{oneDphy01}) with  $\varphi_0  (x)=0$ and  $\varphi_1  (x)=\varphi (x)$
can be represented as follows
\beqst
u(x,t)  
& = &  
  \int_{ 0}^{t} \, db \Big[\frac{1}{4}e^{-t/2}e^{b/2}(2+b) +  \frac{1}{16}  b e^{-3t/2} e^{b/2} ( e^{b }- e^{t }) \Big] \Big[
     \varphi   (x+ e^t- e^b)     
+   \varphi   (x - e^t+ e^b)    \Big] \nonumber \\
&  &
+ \int_{ 0}^{t} b e^{2b}\, db  \int_{ 0}^{    e^t- e^b} dz \, \Big[\varphi  (x-z)  
 +\varphi   (x+z)  \Big]
\Big( \frac{\pa }{\pa z} \Big)^2 E(z ,t;0,b)
\eeqst
as well as by (\ref{3.1n}), 
where 
\be
\label{3.3}
K_1(z,t) 
 = 
\Big[\frac{1}{4}e^{-t/2}(2+\ln (e^{t }-z )) -  \frac{1}{16} e^{-3t/2} z  \ln (e^{t }-z)  \Big]  \frac{1}{\sqrt{e^{t }-z}} 
+  \int_{ 0}^{\ln (e^t-z)}  b e^{2b} 
\Big( \frac{\pa }{\pa z} \Big)^2 E(z ,t;0,b)  db. 
 \ee
\end{corollary}
\medskip

\noindent
{\bf Proof of corollary.}  In  this proof we drop subindex  of $\varphi _1$.
To   prove (\ref{3.1n}) with $K_1(z,t) $ defined by (\ref{3.3}) we apply   (\ref{5.4})  and write
\beqst
u(x,t)  
& = &  
  \int_{ 0}^{t} \, db \Big[\frac{1}{4}e^{-t/2}e^{b/2}(2+b) +  \frac{1}{16}  b e^{-3t/2} e^{b/2} ( e^{b }- e^{t }) \Big] \Big[
     \varphi (x+ e^t- e^b)     
+   \varphi (x - e^t+ e^b)    \Big] \\
&  &
+ \int_{ 0}^{t} b e^{2b}\, db \int_{ x - (e^t- e^b)}^{x+ e^t- e^b} dy \, \varphi (y)  
\Big( \frac{\pa }{\pa y} \Big)^2 E(y-x ,t;0,b) \\ 
& = &  
  \int_{ 0}^{t} \, db \Big[\frac{1}{4}e^{-t/2}e^{b/2}(2+b) +  \frac{1}{16}  b e^{-3t/2} e^{b/2} ( e^{b }- e^{t }) \Big] \Big[
     \varphi (x+ e^t- e^b)     
+   \varphi (x - e^t+ e^b)    \Big] \\
&  &
+ \int_{ 0}^{t} b e^{2b}\, db \int_{  - (e^t- e^b)}^{ e^t- e^b} dz \, \varphi (z+x)  
\Big( \frac{\pa }{\pa z} \Big)^2 E(z ,t;0,b) \\ 
& = &  
  \int_{ 0}^{t} \, db \Big[\frac{1}{4}e^{-t/2}e^{b/2}(2+b) +  \frac{1}{16}  b e^{-3t/2} e^{b/2} ( e^{b }- e^{t }) \Big] \Big[
     \varphi (x+ e^t- e^b)     
+   \varphi (x - e^t+ e^b)    \Big] \\
&  &
+ \int_{ 0}^{t} b e^{2b}\, db  \int_{ 0}^{    e^t- e^b} dz \, \Big[\varphi (x-z)  
 +\varphi (x+z)  \Big]
\Big( \frac{\pa }{\pa z} \Big)^2 E(z ,t;0,b)\,.
\eeqst 
Next we make change $z= e^{b }- e^{t }$, $dz= e^{b }db$, and $b=\ln (z+e^{t })$ in 
\beqst  
&   &  
  \int_{ 0}^{t} \, db \Big[\frac{1}{4}e^{-t/2}e^{b/2}(2+b) +  \frac{1}{16}  b e^{-3t/2} e^{b/2} ( e^{b }- e^{t }) \Big] \Big[
     \varphi (x+ e^t- e^b)     
+   \varphi (x - e^t+ e^b)    \Big] \\ 
& = &  
  \int_{0}^{  e^t-1} \,\Big[      \varphi (x+ z)  +   \varphi (x - z)    \Big] 
\Big[\frac{1}{4}e^{-t/2}(2+\ln (e^{t }-z )) -  \frac{1}{16} e^{-3t/2} z  \ln (e^{t }-z )  \Big]  \frac{1}{\sqrt{e^{t }-z}} dz \,.
\eeqst 
Then
\beqst 
u(x,t)  & = &  
   \int_{0}^{  e^t-1} \,\Big[      \varphi (x+ z)  +   \varphi (x - z)    \Big] 
\Big[\frac{1}{4}e^{-t/2}(2+\ln (e^{t }-z )) -  \frac{1}{16} e^{-3t/2} z  \ln (e^{t }-z )  \Big]  \frac{1}{\sqrt{e^{t }-z}} dz  \\
&  &
+  \int_{ 0}^{    e^t- 1} dz \Big[\varphi (x-z)  
 +\varphi   (x+z)  \Big]\int_{ 0}^{\ln (e^t-z)}\, db \, b e^{2b} 
\Big( \frac{\pa }{\pa z} \Big)^2 E(z ,t;0,b) \\
& = &  
   \int_{0}^{  e^t-1} \,\Big[      \varphi (x- z)  +   \varphi (x + z)    \Big] K_1(z,t)\, dz\,,
\eeqst
where  $K_1(z,t)   $ is defined by (\ref{3.3}). Corollary is proven. \hfill $\square$

The next lemma completes the proof of Theorem~\ref{T1.3}.
\begin{lemma} 
The kernel   $K_1(z,t)   $ defined by (\ref{3.3}) coincides with one given in Theorem~\ref{T1.3}. 
\end{lemma}
\medskip

\noindent
{\bf Proof.} We have by integration by parts
\beqst
\int_{ 0}^{\ln (e^t-z)} \!\!  b e^{2b} 
\Big( \frac{\pa }{\pa z} \Big)^2 E(z ,t;0,b)   db
\!\!& \!\!= \!\!&\!\!
\int_{ 0}^{\ln (e^t-z)} \!\! b \Big( \frac{\pa }{\pa b} \Big)^2 E(z ,t;0,b)   db\\
\!\! &\!\! =\!\!&\!\!
\ln (e^t-z)  \Bigg[ \frac{\pa }{\pa b}  E(z ,t;0,b) \Bigg]_ {b=\ln (e^t-z)}
\!\!-      E(z ,t;0,\ln (e^t-z)) +  E(z ,t;0,0) . 
 \eeqst 
On the other hand,  (\ref{E_2a}) and (\ref{E_7}) of Proposition \ref{P_E} imply 
 \beqst
\int_{ 0}^{\ln (e^t-z)}  b e^{2b} 
\Big( \frac{\pa }{\pa z} \Big)^2 E(z ,t;0,b)  db
& = &
\ln (e^t-z)  \frac{\pa }{\pa b}  E(z ,t;0,\ln (e^t-z)) 
-     \frac{1}{2}  e^{-\frac{t}{2}} 
( e^t   - z )^{-\frac{1}{2}}    +  E(z ,t;0,0) \\
& = &
\ln (e^t-z) \frac{e^{-2 t} \sqrt{e^t} \left(-4 e^t+z\right)}{16 \sqrt{e^t-z}} 
-     \frac{1}{2}  e^{-\frac{t}{2}} 
( e^t   - z )^{-\frac{1}{2}}    +  E(z ,t;0,0) .
 \eeqst
 Thus, for the kernel   $K_1(z,t)   $ defined by (\ref{3.3}) we have
 \beqst
K_1(z,t) 
& = &
\Big[\frac{1}{4}e^{-t/2}(2+\ln (e^{t }-z )) -  \frac{1}{16} e^{-3t/2} z  \ln (e^{t }-z )  \Big]  \frac{1}{\sqrt{e^{t }-z}} \\
&  &
 +  \ln (e^t-z) \frac{e^{-2 t} \sqrt{e^t} \left(-4 e^t+z\right)}{16 \sqrt{e^t-z}} 
-     \frac{1}{2}  e^{-\frac{t}{2}} 
 \frac{1}{\sqrt{e^{t }-z}}    +  E(z ,t;0,0) \\
 & = &
\Big[\frac{1}{4}e^{-t/2}\ln (e^{t }-z ) -  \frac{1}{16} e^{-3t/2} z  \ln (e^{t }-z )  \Big]  \frac{1}{\sqrt{e^{t }-z}} \\
&  &
 +  \ln (e^t-z) \frac{e^{-2 t} \sqrt{e^t} \left(-4 e^t+z\right)}{16 \sqrt{e^t-z}} 
  +  E(z ,t;0,0) \\
 & = &
  E(z ,t;0,0) \,.
 \eeqst
The last line can be easily transformed into $K_1(z,t)$ of Theorem~\ref{T1.3}. Lemma is proven. \hfill $\square$

\section{The Cauchy Problem:  First Datum and $n=1$}
\label{S6}
\setcounter{equation}{0}

In this section we prove Theorem~\ref{T1.3} in the case of  $\varphi_1  (x)=0$. Thus, we have to prove for
the solution $u=u (x,t)$ of the Cauchy problem (\ref{oneDphy01}) with $\varphi_1  (x)=0$
the representation given by   Theorem~\ref{T1.3} in the case of  $\varphi_1  (x)=0$, which is equivalent to 
\[
u(x,t)  
  =    
\frac{1}{2} e ^{-\frac{t}{2}}  \Big[ 
\varphi_0   (x+ e^t- 1)  
+     \varphi_0   (x - e^t+ 1)  \Big]  
+ \, \int_{ 0}^{1} \big[ 
\varphi_0   (x - \phi (t)s)  
+     \varphi_0   (x  + \phi (t)s)  \big] K_0(\phi (t)s,t)\phi (t)\,  ds ,
\]
where $\phi (t) = e^t- 1$. 
The proof of this case consists of the several steps.

\begin{proposition} 
The solution $u= u (x,t)$ of the Cauchy problem (\ref{oneDphy01}) 
can be represented as follows
\beqst
u(x,t)  
& = &  
\frac{1}{2} e ^{-\frac{t}{2}}  \Big[ 
\varphi_0   (x+ e^t- 1)  
+    \varphi_0   (x - e^t+ 1)  \Big]  +
\int_{ 0}^{t} \frac{1}{4}e^\frac{b}{2} e^{-\frac{t}{2} } \Big[
\varphi_0   (x+ e^t- e^b)   
+ \varphi_0   (x - e^t+ e^b)   \Big] \, db \\
&  &
+\int_{ 0}^{t}    \frac{1}{16} e^{-2 t } e^\frac{b}{2}e^{\frac{t}{2}}( e^{b }- e^{t })\Big[  \varphi_0  (x+ e^t- e^b)  
 + \varphi_0  (x -  e^t+ e^b )  
 \Big]  \, db  \\
&  &
+ \int_{ 0}^{t}   e^{2b}\, db \int_{ x - (e^t- e^b)}^{x+ e^t- e^b} dy \, \varphi_0   (y)  
\Big( \frac{\pa }{\pa y}\Big)^2 E(x-y,t;0,b)\,.
\eeqst 
\end{proposition}
\medskip

\noindent
{\bf Proof.} We set $u(x,t)= w(x,t)+ \varphi_0 (x)$, then 
\[
w_{tt} - e^{2t}w_{xx} =e^{2t}\varphi_{0,xx} \, ,\qquad w(x,0) = 0 \,, \qquad w_t(x,0) =0\,  .
\]
Next we plug $f(x,t)= e^{2t}\varphi_{0,xx} $ in the formula given by Theorem~\ref{T1.1}
 and obtain  
\[
w(x,t)   = 
 \int_{ 0}^{t}  e^{2b}\, db \int_{ x - (e^t- e^b)}^{x+ e^t- e^b} dy \, \varphi_0 ^{(2)} (y) 
E(x-y,t;0,b) \,. 
\]
Then we integrate by parts
\beqst
w(x,t) & = &
 \int_{ 0}^{t} e^{2b}\, db  \left(
 \varphi_0 ^{(1)} (x+ e^t- e^b)  E(  -e^t+ e^b ,t;0,b)   
-  \varphi_0 ^{(1)} (x - e^t+ e^b)  E(e^t- e^b  ,t;0,b)\right) \\
&  &
- \int_{ 0}^{t}  e^{2b}\, db \int_{ x - (e^t- e^b)}^{x+ e^t- e^b} dy \, \varphi_0 ^{(1)} (y)  
\frac{\pa }{\pa y} E(x-y,t;0,b)\,.
 \eeqst
On the other hand,
\[
 \varphi_0 ^{(1)} (x+ e^t- e^b) = -  e^{-b} \frac{\pa }{\pa b}  \varphi_0   (x+ e^t- e^b) ,\qquad
 \varphi_0 ^{(1)} (x - e^t+ e^b)=    e^{-b} \frac{\pa }{\pa b}  \varphi_0   (x - e^t+ e^b)
\]
imply
\beqst
w(x,t) & = &
\int_{ 0}^{t}  e^{ b}\, db  \left(
-    \frac{\pa }{\pa b}  \varphi_0   (x+ e^t- e^b)  E(  -e^t+ e^b ,t;0,b)   
-    \frac{\pa }{\pa b}  \varphi_0   (x - e^t+ e^b)  E(e^t- e^b  ,t;0,b)\right) \\
&  &
- \int_{ 0}^{t} e^{2b}\, db \int_{ x - (e^t- e^b)}^{x+ e^t- e^b} dy \, \varphi_0 ^{(1)} (y)  
\frac{\pa }{\pa y} E(x-y,t;0,b) \,.
\eeqst
One more integration by parts leads to
\beqst
&  &
w(x,t)  \\
& = &
-  2  e^{t}  
   \varphi_0  (x )  E( 0 ,t;0,t)     \\
   &  &
 -     \left(
-      \varphi_0   (x+ e^t- 1)  E(  -e^t+ 1 ,t;0,0)   
-       \varphi_0   (x - e^t+ 1)  E(e^t- 1  ,t;0,0)\right) \\ 
&  &
- \int_{ 0}^{t} \, db  \left(
-     \varphi_0   (x+ e^t- e^b)   \frac{\pa }{\pa b} \Big(   e^{ b}E(  -e^t+ e^b ,t;0,b)  \Big) 
-    \varphi_0   (x - e^t+ e^b)   \frac{\pa }{\pa b} \Big(  e^{ b}E(e^t- e^b  ,t;0,b)\Big) \right) \\
&  &
- \int_{ 0}^{t}   e^{2b}\, db \int_{ x - (e^t- e^b)}^{x+ e^t- e^b} dy \, \varphi_0 ^{(1)} (y)  
\frac{\pa }{\pa y} E(x-y,t;0,b)\\
& = &  
-    
   \varphi_0   (x )     +  \frac{1}{2} e ^{-\frac{t}{2}}  \Big( 
\varphi_0   (x+ e^t- 1)  
+    \varphi_0   (x - e^t+ 1)  \Big) \\ 
&  &
- \int_{ 0}^{t} \, db  \left(
-     \varphi_0   (x+ e^t- e^b)   \frac{\pa }{\pa b} \Big(   e^{ b}E(  -e^t+ e^b ,t;0,b)  \Big) 
-    \varphi_0   (x - e^t+ e^b)   \frac{\pa }{\pa b} \Big(  e^{ b}E(e^t- e^b  ,t;0,b)\Big) \right) \\
&  &
- \int_{ 0}^{t}   e^{2b}\, db \int_{ x - (e^t- e^b)}^{x+ e^t- e^b} dy \, \varphi_0 ^{(1)} (y)  
\frac{\pa }{\pa y} E(x-y,t;0,b)\,.
 \eeqst
 \newpage
 
 \noindent
We have used
\beqst 
E( 0 ,t;0,t)=\frac{1}{2} e^{-t} , \qquad E(e^t- 1,t;0,0) = E(1-e^t ,t;0,0) 
 =   
\frac{1}{2} e ^{-\frac{t}{2}} \,.
\eeqst
Hence
\beqst
u(x,t) 
\!\! & \!\! = \!\! &  \!\! 
\frac{1}{2} e ^{-\frac{t}{2}}  \Big( 
\varphi_0   (x+ e^t- 1)  
+    \varphi_0   (x - e^t+ 1)  \Big) \\ 
&  &
- \!  \int_{ 0}^{t} \, db   \Big(
-   \varphi_0   (x+ e^t- e^b)   \frac{\pa }{\pa b} \Big(   e^{ b}E(  -e^t+ e^b ,t;0,b)  \Big) 
-    \varphi_0   (x - e^t+ e^b)   \frac{\pa }{\pa b} \Big(  e^{ b}E(e^t- e^b  ,t;0,b)\Big) \Big) \\
&  &
- \!\int_{ 0}^{t}   e^{2b}\, db \int_{ x - (e^t- e^b)}^{x+ e^t- e^b} dy \, \varphi_0 ^{(1)} (y)  
\frac{\pa }{\pa y} E(x-y,t;0,b) \,.
 \eeqst
Next we apply (\ref{E_2})  and (\ref{E_3a})  of Proposition \ref{P_E} and the integration by parts to obtain
\beqst
\hspace{-0.3cm} &  &
u(x,t) \\
\hspace{-0.3cm} & = & 
\frac{1}{2} e ^{-\frac{t}{2}}  \Big( 
\varphi_0   (x+ e^t- 1)  
+    \varphi_0   (x - e^t+ 1)  \Big)  +
\int_{ 0}^{t} \, db \frac{1}{4}e^\frac{b}{2} e^{-\frac{t}{2} } \Big(
\varphi_0   (x+ e^t- e^b)   
+ \varphi_0   (x - e^t+ e^b)   \Big) \\
\hspace{-0.3cm} &  &
- \int_{ 0}^{t}   e^{2b}\, db  \Big[ \varphi_0   (y)  
\frac{\pa }{\pa y} E(x-y,t;0,b)\Big]_{ y= x - (e^t- e^b)}^{y= x+ e^t- e^b} 
+ \int_{ 0}^{t}   e^{2b}\, db \int_{ x - (e^t- e^b)}^{x+ e^t- e^b} dy \, \varphi_0   (y)  
\Big( \frac{\pa }{\pa y}\Big)^2 E(x-y,t;0,b) .
\eeqst
We have due to (\ref{E_5}) and (\ref{E_6}) of Proposition \ref{P_E} 
\beqst
u(x,t) 
& = &  
\frac{1}{2} e ^{-\frac{t}{2}}  \Big( 
\varphi_0   (x+ e^t- 1)  
+    \varphi_0   (x - e^t+ 1)  \Big)  +
\int_{ 0}^{t}  db \frac{1}{4}e^\frac{b}{2} e^{-\frac{t}{2} } \Big(
\varphi_0   (x+ e^t- e^b)   
+ \varphi_0   (x - e^t+ e^b)   \Big) \\
&  &
- \int_{ 0}^{t}   e^{2b}\, db  \Big[ -\varphi_0  (x+ e^t- e^b)  
\frac{1}{16} e^{-2(b+t)} e^{b/2}e^{t/2}( e^{b }  -e^{t })\\
&  &
\hspace{5cm} + \,\varphi_0  (x - (e^t- e^b))  
\frac{1}{16} e^{-2(b+t)} e^{b/2}e^{t/2}( -e^{b }+ e^{t }) \Big]  \\
&  &
+ \int_{ 0}^{t}   e^{2b}\, db \int_{ x - (e^t- e^b)}^{x+ e^t- e^b} dy \, \varphi_0   (y)  
\Big( \frac{\pa }{\pa y}\Big)^2 E(x-y,t;0,b) \,,
 \eeqst
which coincides with the desired   representation. The proposition is proven.
\hfill $\square$

\noindent
{\bf Completion of the proof of Theorem~\ref{T1.3}.}  We make change $z= e^{b }- e^{t }$, $dz= e^{b }db$, and $b=\ln (z+e^{t })$ in 
the second and third terms of the representation given by the previous proposition:
\beqst 
&   &  
\int_{ 0}^{t}  \frac{1}{4}e^\frac{b}{2} e^{-\frac{t}{2} } \Big[
\varphi_0   (x+ e^t- e^b)   
+ \varphi_0   (x - e^t+ e^b)   \Big] \, db\\
&  &
+\int_{ 0}^{t}    \frac{1}{16} e^{-2 t } e^\frac{b}{2}e^{\frac{t}{2}}( e^{b }- e^{t })\Big[  \varphi_0  (x+ e^t- e^b)  
 + \varphi_0  (x -  e^t+ e^b )  
 \Big]  \, db  \\
& = &  
\int_{0}^{e^t - 1}  \Big[ \frac{1}{4} e^{-\frac{t}{2} } -    \frac{1}{16} e^{-2 t } e^{\frac{t}{2}}z \Big] 
\frac{1}{\sqrt{e^{t }-z} } \Big[  \varphi_0  (x-z)  
 + \varphi_0  (x +z )  
 \Big]  \, dz  \,.
\eeqst
Next we consider the last term apply (\ref{E_2}), and change the order of integration:
\beqst
&  &
\int_{ 0}^{t}   e^{2b}\, db \int_{ x - (e^t- e^b)}^{x+ e^t- e^b} dy \, \varphi_0   (y)  
\Big( \frac{\pa }{\pa y}\Big)^2 E(x-y,t;0,b) \\
& = &
\int_{ 0}^{t}   e^{2b}\, db \int_{0}^{e^t- e^b} dz \, \Big[ \varphi_0   (x-z)  + \varphi_0   (x+z)\Big]
\Big( \frac{\pa }{\pa z}\Big)^2 E(z,t;0,b)  \\
& = &
\int_{ 0}^{e^t-1} dz  \Big[ \varphi_0   (x-z)  + \varphi_0   (x+z)\Big]  \int_{0}^{\ln(e^t- z)} \, e^{2b}\, db
\Big( \frac{\pa }{\pa z}\Big)^2 E(z,t;0,b) \,.
\eeqst
On the other hand, due to\,  $\left( \frac{\pa }{\pa z}\right)^2 E(z,t;0,b)= e^{-2b} \left( \frac{\pa }{\pa b}\right)^2 E(z,t;0,b)$ \,
the last integral is equal to
\beqst
&  &
\int_{ 0}^{e^t-1} dz  \Big[ \varphi_0   (x-z)  + \varphi_0   (x+z)\Big]  \int_{0}^{\ln(e^t- z)}  \Big( \frac{\pa }{\pa b}\Big)^2 E(z,t;0,b) \, db\\
& = &
\int_{ 0}^{e^t-1} dz  \Big[ \varphi_0   (x-z)  + \varphi_0   (x+z)\Big] \Big[  \frac{\pa }{\pa b}  E(z,t;0,\ln(e^t- z)) -  \frac{\pa }{\pa b}  E(z,t;0,0)   \Big]\,.
 \eeqst
According to (\ref{E_7}) and (\ref{E_bzt00}) we have 
\beqst
&   &
 \Big[ \frac{1}{4} e^{-\frac{t}{2} } -    \frac{1}{16} e^{-2 t } e^{\frac{t}{2}}z \Big] 
\frac{1}{\sqrt{e^{t }-z} } +   \frac{\pa E}{\pa b}  (z,t;0,\ln(e^t- z)) -  \frac{\pa  E}{\pa b} (z,t;0,0)  \\
& = &
- \frac{1}{ 2  ((e^t-1)^2-z^2 ) \sqrt{(1+e^t)^2-z^2}} \left\{ (1-e^{2 t}+z^2 ) 
F \left(-\frac{1}{2},\frac{1}{2};1; \frac{  (e^t-1)^2-z^2 }{  (e^t+1)^2-z^2 }   \right) \right.\\
&  &
\left. +   2  (e^t -1) F \left(\frac{1}{2},\frac{1}{2};1;  \frac{  (e^t-1)^2-z^2 }{  (e^t+1)^2-z^2 }   \right)  \right\}\,.
\eeqst
Theorem~\ref{T1.3} is proven.
\hfill $\square$

\section{n-Dimensional Case, $n>1$}
\label{S7}\label{S8} \label{S9}
\setcounter{equation}{0}

{\sl The proof of Theorem~\ref{T1.5}}.  
Let us consider the case $x \in {\mathbb R}^n$, where $n=2m+1$,\, $m \in  {\mathbb N}$.
First for the given function $u=u(x,t)$ we   define the spherical means of 
$u$ about point $x$:
\beqst
I_u(x,r,t) 
 & = &
\frac{1}{\omega_{n-1} } \int_{S^{n-1}  } u(x+ry,t)\, dS_y \,,   
\eeqst
where $\omega_{n-1} $ denotes the area of the unit sphere $S^{n-1} \subset {\mb R}^n$. Then we define
an operator $ \Omega _r $ by 
\[
\Omega _r ( u) (x,t) := \Big( \frac{1}{r} \frac{\partial }{\partial r}\Big)^{m-1} r^{2m-1}I_u(x,r,t) \,.
\]
One can show that there are constants $c_j^{(n)} $, $j=0,\ldots,m-1$, where $n=2m+1$, with
$c_0^{(n)} =1\cdot 3\cdot 5\cdots (n-2)$, \, 
such that
\[
\Big( \frac{1}{r} \frac{\partial }{\partial r}\Big)^{m-1} r^{2m-1} \varphi (r) 
= r \sum_{j=0}^{m-1} c_j^{(n)} r^{j} \frac{\partial^j }{\partial r^j} \varphi (r) \,.
\]
One can recover the functions according to 
\beq
\label{7.1}
u(x,t) 
& = &
\lim_{r \to 0}  I_u (x,r,t) = \lim_{r \to 0}  \frac{1}{c_0^{(n)}r} \Omega _r ( u) (x,t) \,,\\
\label{7.2}
u(x,0) 
& = &
\lim_{r \to 0}  \frac{1}{c_0^{(n)}r} \Omega _r ( u) (x,0) \,,
\quad u_t(x,0) = \lim_{r \to 0}  \frac{1}{c_0^{(n)}r} \Omega_ r ( \partial_t u) (x,0)  \,. 
\eeq
It is well known that $\Delta_x \Omega_r  h =  
\frac{\pa^2}{\pa\, r^2}  \Omega_r  h $  for every function $h \in C^2({\mb R}^n)$.
Therefore we arrive at the following mixed problem for the  function $v(x,r,t) := \Omega_r  (u )(x,r,t) $: 
\beqst
\cases{
v_{tt} (x,r,t) - e^{2t} v_{rr}  (x,r,t) = F (x,r,t)  
\quad  \mbox{\rm for all}\quad   t \ge 0\,, 
\,\, r \ge 0\,, \,\, x \in {\mb R}^n\,,\cr
v(x,0,t) = 0 
\quad \mbox{\rm for all}\quad t \ge 0\,, \quad x \in {\mb R}^n\cr
v(x,r,0) = 0 \,, \quad v_t(x,r,0) = 0
\quad \mbox{\rm for all}
\quad r \ge 0\,, \quad x \in {\mb R}^n\,,\cr
F (x,r,t)  := \Omega_r  (f ) (x,t) \,, \quad  F (x,0,t)  =0 \,, 
\quad \mbox{\rm for all}\quad x \in {\mb R}^n\,.}
\eeqst
It must be noted here that the spherical mean $I_u$ defined for $r>0$ has an 
extension as even function for $r<0 $ and hence $\Omega_r  (u ) $ has a natural extension as an
odd function. That allows replacing the mixed problem with the Cauchy problem. Namely, 
let functions $\wt v$ and $\wt F$ be the continuations of the functions $v$ and $F$, respectively, by 
\[
\hspace*{-0.7cm} 
\wt v(x,r,t)   = \cases{\,v (x,r,t) , \,\, if \,\,r \ge 0 \cr  - v (x,- r,t) , \,\,if \,\,r \le 0 }\,, \quad 
\wt F(x,r,t)   = \cases{\,F (x,r,t) , \,\,if \,\,r \ge 0 \cr  - F (x,- r,t) , \,\,if \,\,r \le 0 } \,.
\]
Then $\wt v$ solves the Cauchy problem
\beqst
&  &
\wt v_{tt} (x,r,t) - e^{2t} \wt v_{rr}  (x,r,t) = \wt F (x,r,t)  
\quad \mbox{\rm for all}\quad t \ge 0\,, 
\quad r \in {\mb R} \,, \quad x \in {\mb R}^n\,,\\
&  &
\wt v(x,r,0) = 0 \,, \quad \wt v_t(x,r,0) = 0
\quad \mbox{\rm for all}
\quad r \in {\mb R} \,, \quad x \in {\mb R}^n .
\eeqst 
Hence according to Theorem \ref{T4.1} one has the representation 
\[
 \wt v (x,r,t)   = 
\int_{ 0}^{t} db \int_{ r - (e^t- e^b)}^{r+ e^t- e^b} dr_1 \,  \wt F (x,r_1,b)   
\frac{1}{\sqrt{(e^t  + e^b )^2 - (r-r_1)^2 }} 
F\left(\frac{1}{2},\frac{1}{2};1; 
\frac{ (e^t  - e^b )^2 - (r-r_1)^2 }
{  (e^t  + e^b )^2 - (r-r_1)^2 } \right) . 
\]
Since $u(x,t) = \lim_{r \to 0} \big( \wt v (x,r,t)/(c_0^{(n)}r)\big) $,
we consider the case with $r < t$ in the above representation to obtain:
\beqst
u(x,t) 
& = &
\lim_{r \to 0}
\frac{1}{c_0^{(n)}r} \int_{ 0}^{t} db \int_{  - (e^t- e^b)}^{ e^t- e^b} dz \,  \wt F (x,z+r,b)   
\frac{1}{\sqrt{(e^t  + e^b)^2-z^2}}
F\left(\frac{1}{2},\frac{1}{2};1; 
\frac{ (e^t  - e^b)^2-z^2}
{  (e^t  + e^b)^2-z^2} \right) \\
&  = &
\frac{1}{c_0^{(n)} }   \int_{ 0}^{t} db
  \int_{ 0}^{ e^t- e^b} dr_1 \,  \lim_{r \to 0}
\frac{1}{ r}  \left\{\wt F (x,r+r_1,b) +   \wt F (x,r-r_1,b)\right\} \\
  &  & 
\hspace*{3cm}\times \frac{1}{\sqrt{(e^t  + e^b)^2-r_1^2}}
F\left(\frac{1}{2},\frac{1}{2};1; 
\frac{ (e^t  - e^b)^2-r_1^2}
{  (e^t  + e^b)^2-r_1^2} \right) \,.
\eeqst
Then by definition of the function $\wt F $ we replace 
$\lim_{r \to 0} \frac{1}{ r} \Big\{  \wt F (x,r- r_1,b)  + \wt F (x,r+ r_1,b) \Big\}  $
with \\
$2\Big( \frac{\partial }{\partial r} F (x,r ,b) \Big)_{r=r_1}  $ in the last formula. 
The definitions of $F (x,r ,t)$  and of the operator $\Omega _r $  yield:
\beqst
u(x,t) 
& = &
\frac{2}{c_0^{(n)} }   \int_{ 0}^{t} db
  \int_{ 0}^{ e^t- e^b} dr_1 \,  \left( \frac{\partial }{\partial r} 
\Big( \frac{1}{r} \frac{\partial }{\partial r}\Big)^{m-1} r^{2m-1}I_f(x,r,t) 
\right)_{r=r_1} \\
  &  & 
\hspace*{2cm}\times \frac{1}{\sqrt{(e^t  + e^b)^2-r_1^2}}
F\left(\frac{1}{2},\frac{1}{2};1; 
\frac{ (e^t  - e^b)^2-r_1^2}
{  (e^t  + e^b)^2-r_1^2} \right) \,,
\eeqst
where $x \in {\mb R}^n$, $n=2m+1$, $m \in {\mb N}$. Thus the solution to the Cauchy problem 
is given by (\ref{1.27}). We employ the method of descent to complete the proof for 
the case with even $n$, $n=2m$, $m \in {\mb N}$. 
Theorem~\ref{T1.5} is proven.
\hfill $\square$

\noindent
{\bf Proof of (\ref{E+}) and (\ref{E-}).} We set $f(x,b)=\delta (x)\delta (t-t_0) $ in  (\ref{1.27}) and (\ref{1.28}),
and we obtain  (\ref{E+}) and (\ref{E-}),  
where if $n$ is odd,   
\[
E^{w} (x,t) :=\frac{1}{\omega_{n-1} 1\cdot 3\cdot 5\ldots \cdot (n-2) } \frac{\partial }{\partial t} 
\Big( \frac{1}{t} \frac{\partial }{\partial t}\Big)^{\frac{n-3}{2} } 
\frac{1}{t} \delta (|x|-t) \, ,
\] 
while for $n$ even  we have
\[
 E^{w} (x,t) :=\frac{2}{\omega_{n-1} 1\cdot 3\cdot 5\ldots \cdot (n-1) } \frac{\partial }{\partial t} 
\Big( \frac{1}{t} \frac{\partial }{\partial t}\Big)^{\frac{n-2}{2} } 
\frac{1}{\sqrt{t^2-|x|^2}}\chi _{B_t(x)} \, .
\] 
Here $\chi _{B_t(x)} $ denotes  the characteristic function of the ball $B_t(x):= \{x \in {\mb R}^{n};\, |x| $ $\leq t \}  $.
Constant $\omega_{n-1} $ is the area of the unit sphere $S^{n-1} \subset {\mb R}^n$. The distribution $\delta (|x|-t) $ is defined  by 
\[
\label{delta}
<\delta (|\cdot |-t) , f (\cdot ) > = \int_{|x|=t}f(x)\, dx 
\quad \mbox{\rm for all} \quad f \in C^\iy_0({\mb R}^{n}) \,.
\]

{\sl The proof of Theorem~\ref{T1.6}}.  First we consider  case of  $\varphi_0  (x)=0$. More precisely, we have to prove that 
the solution $u (x,t)$ of the Cauchy problem (\ref{1.30CP}) with  $\varphi_0  (x)=0$
can be represented by (\ref{1.30}) with  $\varphi_0  (x)=0$. The next lemma will be used in both cases.

\begin{lemma}
\label{L8.2} \label{L9.2}
Consider the mixed problem
\begin{eqnarray*}
\cases{ 
v_{tt} - e^{2t}   v_{rr} =0 \quad \mbox{\rm for all}\quad t \ge 0\,, 
\quad r \ge 0\,,  \quad \cr  
v(r,0)= \tau _0 (r) \,, \quad v_t(r,0)= \tau _1 (r) \quad \mbox{\rm for all} 
\quad r \ge 0\,, \cr
v(0,t)=0 \quad \mbox{\rm for all}\quad t \ge 0\,,}
\end{eqnarray*}
and denote by  $\widetilde \tau _0 (r) $ and $\widetilde \tau _1 (r) $  the 
continuations of the functions $\tau _0 (r) $ and $\tau _1 (r) $  for negative $r$ as odd functions: 
\, $\widetilde \tau _0 (-r) =- \tau _0 (r)$ and  $\widetilde \tau _1 (-r) =- \tau _1 (r)$  for all $r \geq 0$, respectively. 
Then solution $v(r,t) $ to the mixed problem is given by the restriction of (\ref{3.1n}) to
$r \ge 0$:
\begin{eqnarray*} 
  v(r,t) 
& = &
\frac{1}{2} e ^{-\frac{t}{2}}  \Big[ 
 \widetilde \tau _0 (r+ e^t- 1)  
+  \widetilde  \tau _0  (r - e^t+ 1)  \Big]   
+ \, \int_{ 0}^{1} \big[ 
 \widetilde \tau _0  (r - \phi (t)s)  
+   \widetilde \tau _0  (r  + \phi (t)s)  \big] K_0(\phi (t)s,t)\phi (t)\,  ds \\
&  &
+ \int_0^1\Big[  \widetilde \tau_1 \Big( r +  \phi (t) s \Big) +  \widetilde \tau_1 \Big( r -  \phi (t) s \Big) \Big] K_1(\phi (t)s,t) \phi (t)\,ds \,, 
\end{eqnarray*}
where $K_0(z,t)$ and  $K_1(z,t)$ are defined in Theorem~\ref{T1.3} and $\phi (t)= e^t-1$.
\end{lemma}

\medskip

\noindent
{\bf Proof.} This lemma is a direct consequence of Theorem~\ref{T1.3}. \hfill $\square$

\medskip

Now let us  consider the case $x \in {\mathbb R}^n$, where $n=2m+1$.
First for the given function $u=u(x,t)$ we   define the spherical means of 
$u$ about  point $x$. 
One can recover the functions by means of (\ref{7.1}), (\ref{7.2}), and
\begin{eqnarray*}
\varphi_i(x) 
& = &
\lim_{r \to 0}  I_{\varphi_i} (x,r) = \lim_{r \to 0}  \frac{1}{c_0^{(n)}r} \Omega_r ( \varphi_i) (x) \,, \quad i=0,1   \,.
\end{eqnarray*}
Then 
we arrive at the following mixed problem 
\begin{eqnarray*}
\cases{
v_{tt} (x,r,t) - e^{2t} v_{rr}  (x,r,t) = 0  
\quad  \mbox{\rm for all}\quad   t \ge 0\,, 
\,\, r \ge 0\,, \,\, x \in {\mathbb R}^n\,,\cr
v(x,0,t) = 0 
\quad \mbox{\rm for all}\quad t \ge 0\,, \quad x \in {\mathbb R}^n\,,\cr
v(x,r,0) = 0 \,, \quad v_t(x,r,0) = \Phi _1(x,r)
\quad \mbox{\rm for all}
\quad r \ge 0\,, \quad x \in {\mathbb R}^n\,,}
\end{eqnarray*}
with the  unknown function $v(x,r,t) := \Omega_r  (u )(x,r,t) $, where
\begin{eqnarray}
\label{8.2}
&  &
\Phi _i (x,r)  := \Omega_r  (\varphi_i  ) (x) 
= \Big( \frac{1}{r} \frac{\partial }{\partial r}\Big)^{m-1} r^{2m-1}\frac{1}{\omega_{n-1} } \int_{S^{n-1}  } 
\varphi_i(x+ry)\, dS_y
\,, \\
\label{8.3}
&  &
\Phi _i (x,0)  =0 \,, \quad i=0,1, \quad
\quad \mbox{\rm for all}\quad x \in {\mathbb R}^n\,.
\end{eqnarray} 
Then, according to Lemma \ref{L8.2} and   
$u(x,t) = \lim_{r \to 0} \big(  v (x,r,t)/(c_0^{(n)}r)\big) $,
 we obtain:
\begin{eqnarray*}
u(x,t) 
 & = &
\frac{  1}{c_0^{(n)}} \lim_{r \to 0} \frac{  1}{r} \int_0^1\Big[ \wt \Phi _1 \big(x, r +  \phi (t) s \big) +   \wt \Phi _1 \big(x, r -  \phi (t) s \big) \Big] K_1(\phi (t)s,t) \phi (t)\,ds  \,.
\end{eqnarray*}
The last limit is equal to 
\begin{eqnarray*}
&   &
2 \int_0^1\left( \frac{ \partial }{\partial r} \Phi _1 (x,  r  )   
  \right)_{r= \phi (t)  s}  K_1(\phi (t)s,t) \phi (t)\,ds  \\
& = &
 2\int_0^1 \left( \frac{\partial}{\partial r} \Big( \frac{1}{r} \frac{\partial }{\partial r}\Big)^{\frac{n-3}{2} } 
\frac{r^{n-2}}{\omega_{n-1} } \int_{S^{n-1}  } 
\varphi_1 (x+ry)\, dS_y \right)_{r=\phi (t) s}    K_1(\phi (t)s,t) \phi (t)\,ds  \,.
\end{eqnarray*}
Thus, Theorem~\ref{T1.6} in the case of  $\varphi_0  (x)=0$ is proven.
\bigskip

Now we turn to  the case of  $\varphi_1  (x)=0$. 
Thus, we arrive at the following mixed problem 
\begin{eqnarray*}
\cases{
v_{tt} (x,r,t) - e^{2t} v_{rr}  (x,r,t) = 0  
\quad  \mbox{\rm for all}\quad   t \ge 0\,, 
\,\, r \ge 0\,, \,\, x \in {\mathbb R}^n\,,\cr
v(x,r,0) = \Phi _0(x,r) \,, \quad v_t(x,r,0) = 0
\quad \mbox{\rm for all}
\quad r \ge 0\,, \quad x \in {\mathbb R}^n\,,\cr
v(x,0,t) = 0 
\quad \mbox{\rm for all}\quad t \ge 0\,, \quad x \in {\mathbb R}^n\,,}
\end{eqnarray*}
with the  unknown function $v(x,r,t) := \Omega_r  (u )(x,r,t) $ defined by (\ref{8.2}), (\ref{8.3}). 
Then, according to Lemma \ref{L9.2} and   
$u(x,t) = \lim_{r \to 0} \big(  v (x,r,t)/(c_0^{(n)}r)\big) $,
 we obtain:
\begin{eqnarray*}
u(x,t) 
 & = &
 \frac{1}{c_0^{(n)}} e ^{-\frac{t}{2}} \lim_{r \to 0}\frac{1}{2r} \Big[ 
 \widetilde \Phi _0 (x,r+ e^t- 1)  
+  \widetilde  \Phi _0 (x,r - e^t+ 1)  \Big]   \\
&  &
+ \, \frac{2}{c_0^{(n)}} \int_{ 0}^{1} \lim_{r \to 0}\frac{1}{2r} \big[ 
 \widetilde \Phi _0 (x,r - \phi (t)s)  
+   \widetilde \Phi _0 (x,r  + \phi (t)s)  \big] K_0(\phi (t)s,t)\phi (t)\,  ds\,, \\
& = &
 \frac{1}{c_0^{(n)}} e ^{-\frac{t}{2}} \left( \frac{\pa }{\pa r} 
\Phi _0 (x,r)  \right)_{r=\phi (t)}   
+ \, \frac{2}{c_0^{(n)}} \int_{ 0}^{1} \left( \frac{\pa }{\pa r} 
\Phi _0 (x,r)  \right)_{r=\phi (t)s}  K_0(\phi (t)s,t)\phi (t)\,  ds\\
& = &
 e ^{-\frac{t}{2}} v_{\varphi_0}  (x, \phi (t))
+ \, 2\int_{ 0}^{1} v_{\varphi_0}  (x, \phi (t)s) K_0(\phi (t)s,t)\phi (t)\,  ds\,.
\end{eqnarray*}
Theorem~\ref{T1.6} is proven. \hfill $\square $

\section{$L^p-L^q$ Decay and $L^q-L^q$ Estimates for the Solutions of  One-dimensional Equation,   $n=1$}
\label{S10}
\setcounter{equation}{0}

Consider now the Cauchy problem for the equation (\ref{Int_3})  with the source term and with vanishing initial  data (\ref{Int_4}). 
\begin{theorem}
\label{T10.1}
For every function $f \in C^2 ({\mb R}\times [0,\infty))$ such that $f(\cdot ,t) \in C_0^\infty ({\mb R}_x)$ the solution 
$u = u(x,t)$ of the Cauchy problem (\ref{Int_3}), (\ref{Int_4})  satisfies inequality
\beqst
\| u(x,t) \| _{L^q({\mb R}_x)} & \leq  &
C  e^{\frac{t}{\rho}-t} \int_{ 0}^{t}   (1+  t-b  )        \|  f(x,b) \| _{L^p({\mb R}_x)} \, db 
\eeqst
for all $t >0$, where $1 <p<\rho '$, $ \frac{1}{q} = \frac{1}{p} - \frac{1}{\rho '} $, $ \rho <2$, $ \frac{1}{\rho } + \frac{1}{\rho '}= 1$. 
\end{theorem}
\medskip

\noindent
{\bf Proof.} 
Using the fundamental solution from Theorem \ref{T1} one can write the convolution 
\[
\hspace*{-0.5cm} u(x,t)
 = 
\int_{ -\iy}^{\iy} \! \int_{ -\iy}^{\iy} E_+ (x,t;y,b)f(y,b)\, db\, dy 
 = 
\int_{ 0}^{t}\! db \!\int_{ -\iy}^{\iy} E_+ (x-y,t;0,b) f(y,b) \,dy  \,.
\]
Due to Young's inequality we have
\beqst
\| u(x,t) \| _{L^q({\mb R}_x)} & \leq  &
c_k \int_{ 0}^{t} db \Bigg( \int_{  - (\phi (t) - \phi (b))}^{ \phi (t) - \phi (b)} 
|E (x,t;0,b)|^\rho  dx \Bigg)^{1/\rho } \|  f(x,b) \| _{L^p({\mb R}_x)} ,
\eeqst
where $1 <p<\rho '$, $ \frac{1}{q} = \frac{1}{p} - \frac{1}{\rho '} $, $ \frac{1}{\rho } + \frac{1}{\rho '}= 1$. The integral in parentheses can be transformed as follows 
\[
 \int_{  - (\phi (t) - \phi (b))}^{ \phi (t) - \phi (b)} 
|E (x,t;0,b)|^\rho  dx 
 = 
2 e^{b-b\rho }\int_{  0}^{ e^{t-b}  - 1} 
((e^{t-b}  + 1)^2  - r^2  )^{-\frac{\rho }{2}}  
F\left(\frac{1}{2},\frac{1}{2};1; 
\frac{ (e^{t-b}  - 1 )^2 -  r^2   }
{ (e^{t-b}  + 1 )^2  - r^2} \right)^\rho  d r .
\]
\begin{lemma}
\label{L10.2}
For all $z>1$ the following estimate 
\[
\int_{  0}^{ z  - 1} 
((z  + 1)^2  - r^2  )^{-\frac{\rho }{2}}  
F\left(\frac{1}{2},\frac{1}{2};1; 
\frac{ (z  - 1 )^2 -  r^2   }
{ (z  + 1 )^2  - r^2} \right)^\rho  d r  
  \leq  
C (1+ \ln   z )^\rho (z  - 1) (z  + 1)^{-\rho } F\Big(\frac{1}{2},\frac{\rho }{2};\frac{3}{2}; \frac{  (z  - 1)^2  }{ (z  + 1)^2  } \Big)
\]
is fulfilled, provided that $1 <p<\rho '$, $ \frac{1}{q} = \frac{1}{p} - \frac{1}{\rho '} $, $ \frac{1}{\rho } + \frac{1}{\rho '}= 1$. 
In particular, if $\rho < 2$, then 
\beqst
\int_{  0}^{ z  - 1} 
((z  + 1)^2  - r^2  )^{-\frac{\rho }{2}}  
F\left(\frac{1}{2},\frac{1}{2};1; 
\frac{ (z  - 1 )^2 -  r^2   }
{ (z  + 1 )^2  - r^2} \right)^\rho  d r 
& \leq &
C (1+ \ln   z )^\rho (z  - 1) (z  + 1)^{-\rho } \,.
\eeqst
\end{lemma}
\medskip

\noindent
{\bf Proof.} We rewrite the argument of the hypergeometric function as follows
\[
\frac{ (z  - 1 )^2 -  r^2   }
{ (z  + 1 )^2  - r^2} = 1- \frac{  4z  }{ (z  + 1)^2 -r ^2}\,.
\]
If
\be
\label{r}
r \ge \sqrt{(z  + 1)^2-8z} \,,
\ee
then
\be
\label{9.2}
\frac{  4z  } { (z  + 1)^2 -r ^2}  \ge \frac{1}{2}    \Longrightarrow 0< 1- \frac{  4z  } { (z  + 1)^2 -r ^2} \le \frac{1}{2} 
\ee
for such $r$ and $z$  implies
\be
\label{9.3}
\left| F\left(\frac{1}{2},\frac{1}{2};1; 1-
\frac{  4z  }
{ (z  + 1)^2 -r ^2} \right) \right| \le C  \,.         
\ee
Hence for $\rho >0$ we have
\[
\int_{  \sqrt{(z  + 1)^2-8z}}^{ z  - 1} 
((z  + 1)^2  - r^2  )^{-\frac{\rho }{2}}  
F\left(\frac{1}{2},\frac{1}{2};1; 
1- \frac{  4z  }{ (z  + 1)^2 -r ^2} \right)^\rho  d r 
 \leq 
C(z  - 1) (z  + 1)^{-\rho } F\left(\frac{1}{2},\frac{\rho }{2};\frac{3}{2}; \frac{  (z  - 1)^2  }{ (z  + 1)^2  } \right).
\]
If
\[
r \le \sqrt{(z  + 1)^2-8z} \quad {\rm and} \quad z\ge 6  \,,
\]
then $8 <8z \leq (z  + 1)^2 - r^2 \leq (z  + 1)^2 $,
implies
\[
\left| F\left(\frac{1}{2},\frac{1}{2};1; 1-
\frac{  4z  }
{ (z  + 1)^2 -r ^2} \right) \right| \le C \left| \ln \left( \frac{  4z  }
{ (z  + 1)^2 -r ^2} \right)  \right|  \le C   (1+ \ln   z )  \,.       
\]
Hence 
\beqst
&  &
\int_{0 }^{\sqrt{(z  + 1)^2-8z  }} 
((z  + 1)^2  - r^2  )^{-\frac{\rho }{2}}  
F\left(\frac{1}{2},\frac{1}{2};1; 
1- \frac{  4z  }{ (z  + 1)^2 -r ^2} \right)^\rho  d r \\
& \leq &
C   (1+ \ln   z )^\rho (z  - 1) (z  + 1)^{-\rho } F\Big(\frac{1}{2},\frac{\rho }{2};\frac{3}{2}; \frac{  (z  - 1)^2  }{ (z  + 1)^2  } \Big)\,.
\eeqst
The lemma is proven. \hfill $\square $
\medskip

\noindent
{\bf Completion of the proof of Theorem~\ref{T10.1}.} Thus for $\rho <2$ and $ z= e^{t-b}$ we have 
\beqst
\| u(x,t) \| _{L^q({\mb R}_x)} 
& \leq  &
c  \int_{ 0}^{t}  e^{\frac{b}{\rho}-b } (1+ \ln   z )  (z  - 1)^{1/\rho } (z  + 1)^{-1 }   \|  f(x,b) \| _{L^p({\mb R}_x)} \, db \\
& \leq  &
c  \int_{ 0}^{t}  e^{\frac{b}{\rho}-b } (1+  t-b  )  (e^{t-b}  - 1)^{1/\rho } (e^{t-b}  + 1)^{-1 }   \|  f(x,b) \| _{L^p({\mb R}_x)} \, db \\
& \leq  &
c  \int_{ 0}^{t}  e^{\frac{b}{\rho} -b} (1+  t-b  )  e^{\frac{t}{\rho}-\frac{b}{\rho}}   e^{-t+b}     \|  f(x,b) \| _{L^p({\mb R}_x)} \, db \,.
\eeqst
The last inequality implies the  estimate of the statement of theorem. Theorem~\ref{T10.1} is proven. \hfill $\square $
\bigskip

\begin{proposition}
\label{P9.3}
The solution $u=u(x,t) $ of the Cauchy problem 
\[
u_{tt} - e^{2t}u_{xx} =0\, ,\qquad u(x,0) = \varphi_0  (x) \,, \qquad u_t(x,0) =\varphi_1  (x)\,  ,
\]
with $\varphi_0  $, $\varphi_1  \in C_0^\infty ({\mb R})$ satisfies the following estimate
\be 
\label{9.1}   
\|  u(x,t) \|_{ { L}_{  q} ({\mathbb R})  } 
  \le     
 C \left(  \|  \varphi_0   (x)   \|_{ { L}_{  q} ({\mathbb R})  } +  (1+ t )\|  \varphi_1   (x)   \|_{ { L}_{  q} ({\mathbb R})  } \right)
\, \quad \mbox{for all} \,\quad t \in (0,\infty) .
\ee 
\end{proposition}
\medskip  

\noindent
{\bf Proof.} First we consider the  equation without source term but with the second datum that is the case of $\varphi _0=0$. 
For the convenience we drop subindex of
$\varphi _1$.
Then we apply the representation given by Theorem~\ref{TInt_2} for the solution $u=u (x,t)$ of the Cauchy problem with $\varphi _0=0$,
and obtain
\[  
\|  u(x,t) \|_{ { L}^{  q} ({\mathbb R})  } 
 \le    
2 \|  \varphi   (x)   \|_{ { L}^{  q} ({\mathbb R})  }  \int_{0}^{  e^t-1} \, |K_1(r,t)| dr 
\,.
\]
To estimate the last integral we write
\beq   
\label{10.1} 
 \int_{0}^{  e^t-1} \, |K_1(r,t)| dr 
& \leq  &   
I_1(e^t)  
\,,
\eeq 
where for $z= e^t>1$ we denote
\beq 
\label{10.2}   
I_1 (z)
& :=  &   
 \int_{0}^{  z-1} \, \frac{1}{  \sqrt{(1+z)^2-r^2 } }    F\left(\frac{1}{2},\frac{1}{2};1; 
\frac{  r^2 - (z  - 1 )^2  }
{   r^2 - (z  + 1 )^2 } \right) dr
\,.
\eeq 
Then, according to Lemma \ref{L10.2} (where $\rho =1$) we have for that integral the following estimate
\be 
\label{10.4I1}  
I_1 (e^t)
  \leq  
C (1+ t )\,.
\ee 
\newpage

\noindent
Finally, (\ref{10.1}) to  (\ref{10.4I1})  imply the $L^q- L^q$ estimate (\ref{9.1})  for the case of $\varphi _0=0$.
\medskip

Next we consider the equation without source but with the first datum, that is,   the  case of $\varphi _1=0$. 
We apply the representation given by Theorem~\ref{TInt_2} for the solution $u=u (x,t)$ of the Cauchy problem with $\varphi _1=0$,
and obtain
\beqst
\|  u(x,t) \|_{ { L}^{  q} ({\mathbb R})  } 
& \le &   
e ^{-\frac{t}{2}} 
\|  \varphi_0   (x) \|_{ { L}^{  q} ({\mathbb R})  }  
+ 2\|  \varphi_0   (x) \|_{ { L}^{  q} ({\mathbb R})  }  \int_{ 0}^{e^t-1} \left|   K_0(z,t)  \right|  dz \,.
\eeqst
Thus, we have to estimate the integral 
$
\int_{ 0}^{e^t-1} \left|   K_0(r,t)  \right|  dr \,.
$
The following lemma completes the proof of proposition.
\begin{lemma} 
\label{L10.3}
The kernel $K_0(r,t)$ has an integrable singularity at $r=e^t-1$, more precisely, one has 
\beqst
\int_{ 0}^{e^t-1} \left|   K_0(r,t)  \right|  dr \leq C \quad \mbox{for all} \quad    t \in[0,\infty)\,. 
\eeqst 
\end{lemma}
\medskip

\noindent
{\bf Proof.} Consider the argument $\frac{ (e^t-1)^2 -r^2   }{ (e^t+1)^2 -r^2 }  $ of the hypergeometric function and its derivative.
Denote $z= e^t$, then $0 \leq \frac{ (e^t-1)^2 -r^2   }{ (e^t+1)^2 -r^2 } =  \frac{ (z-1)^2 -r^2   }{ (z+1)^2 -r^2 } \leq 1$.
The formula (\ref{4.12})   describes the behavior of those functions at the neighbourhood of zero.
Hence, if $\varepsilon  >0$ is small, then for all $z$ and $r$ such that 
\be
\label{9.7}
\frac{ (z-1)^2 -r^2   }{ (z+1)^2 -r^2 } \leq \varepsilon 
\ee
one has
\be
\label{9.8}
F\Big(\pm \frac{1}{2},\frac{1}{2};1; \frac{ (z-1)^2 -r^2   }{ (z+1)^2 -r^2 }   \Big) 
=  1 \pm \frac{1 }{4}\frac{ (z-1)^2 -r^2   }{ (z+1)^2 -r^2 }   + O\left(\left( \frac{ (z-1)^2 -r^2   }{ (z+1)^2 -r^2 }\right)^2\right)   \,. 
\ee
Consider therefore two zones, 
\beq
\label{9.10}
Z_1(\varepsilon ,z) 
& := &
\left\{ (z,r) \,\Big|\, \frac{ (z-1)^2 -r^2   }{ (z+1)^2 -r^2 } \leq \varepsilon,\,\, 0 \leq r \leq z-1 \right\} ,\\ 
\label{9.11}
Z_2(\varepsilon ,z) 
& := &
\left\{ (z,r) \,\Big|\, \varepsilon \leq  \frac{ (z-1)^2 -r^2   }{ (z+1)^2 -r^2 },\,\, 0 \leq r \leq z-1  \right\}.
\eeq
We  split integral into two parts:
\beqst
\int_{ 0}^{e^t-1} \left|   K_0(r,t)  \right|  dr 
& = &
\int_{ (z,r) \in Z_1(\varepsilon ,z)   }   \left|   K_0(r,t)  \right|  dr 
+ \int_{ (z,r) \in Z_2(\varepsilon ,z)  }  \left|   K_0(r,t)  \right|  dr \,.
\eeqst 
In the first zone we have
\beq
\label{9.15}
 &  &
\left| (1-z^{2}+r^2 ) 
F\Big(-\frac{1}{2},\frac{1}{2};1; \frac{ (z-1)^2 -r^2   }{ (z+1)^2 -r^2 }   \Big) 
+   2  (z -1) F\Big(\frac{1}{2},\frac{1}{2};1; \frac{ (z-1)^2 -r^2   }{ (z+1)^2 -r^2 }   \Big) \right| \nonumber \\
& = &
\big( (z-1)^{2}-r^2 \big)\left|1  + \frac{1 }{4}\frac{ 3-z^{2} -2   z  +r^2     }{ (z+1)^2 -r^2 }   \right|  
+ \big(z^{2}+2z -3-r^2  \big)O\left(\left( \frac{ (z-1)^2 -r^2   }{ (z+1)^2 -r^2 }\right)^2\right) \,.
\eeq
Consider therefore,
\beqst
A_1
&  :=  &
\int_{ (z,r) \in Z_1(\varepsilon ,z)  }  \frac{1}{     \sqrt{(z+1)^2-r^2}    }   dr   \leq   
\int_{ 0}^{z-1}  \frac{1}{   \sqrt{(z+1)^2-r^2}    } dr  
  \leq   
\frac{\pi}{2} \qquad \mbox{\rm for all}  \qquad z \in [1,\infty) \,, \\
A_2
& := &
\int_{  (z,r) \in Z_1(\varepsilon ,z) }  \frac{1}{ ((z+1)^2 -r^2 )\sqrt{(z+1)^2-r^2}    }  \big|3-z^{2} -2   z  +r^2 \big|  dr \\
& \leq   &
C\int_{ 0}^{z-1}  \frac{1}{ \sqrt{(z+1)^2-r^2}    }    dr \\
& \leq   &
\frac{\pi}{2} \qquad \mbox{\rm for all}  \qquad z \in [1,\infty)  \,,
\eeqst 
and
\beqst
A_3
& := &
\int_{  (z,r) \in Z_1(\varepsilon ,z) }  \frac{z^2+2z-3-r^2}{  \sqrt{(z+1)^2-r^2}    } \frac{ (z-1)^2 -r^2   }{ ((z+1)^2 -r^2)^2 }   dr \\
& \leq  &
\int_{  (z,r) \in Z_1(\varepsilon ,z) }  \frac{z^2+2z-3-r^2}{  \sqrt{(z+1)^2-r^2}    } \frac{ 1   }{  (z+1)^2 -r^2  }   dr \\
& \leq  &
\int_{  (z,r) \in Z_1(\varepsilon ,z) }  \frac{1}{  \sqrt{(z+1)^2-r^2}    }   dr\\
& \leq  &
 \frac{ \pi   }{  2  }  \qquad \mbox{\rm for all}  \qquad z \in [1,\infty)  \,.
\eeqst
Finally,
\beqst
&  &
\int_{  (z,r) \in Z_1(\varepsilon ,z) }  dr  \frac{1}{     ((z-1)^2-r^2 ) \sqrt{(z+1)^2-r^2}    } \\
&  &
\times \left| (1-z^{2}+r^2 ) 
F\Big(-\frac{1}{2},\frac{1}{2};1; \frac{ (z-1)^2 -r^2   }{ (z+1)^2 -r^2 }   \Big) 
+   2  (z -1) F\Big(\frac{1}{2},\frac{1}{2};1; \frac{ (z-1)^2 -r^2   }{ (z+1)^2 -r^2 }   \Big)    \right| \\
& \leq   &
C  \qquad \mbox{\rm for all}  \qquad z \in [1,\infty) \,.
\eeqst
In the second zone we have
\be
\label{10.4n}
\varepsilon \leq  \frac{ (z-1)^2 -r^2   }{ (z+1)^2 -r^2 } \leq 1 \quad \Longrightarrow \quad 
\frac{ 1  }{ (z-1)^2 -r^2 }  \leq  \frac{ 1   }{ \varepsilon[(z+1)^2 -r^2] }\,.
\ee
According to the formula 15.3.10 of \cite[Ch.15]{B-E}   the hypergeometric functions obey the estimates
\beq
\label{FGH+xat1}
&  &
\left| F\Big(-\frac{1}{2},\frac{1}{2};1; x   \Big) \right|  \leq C \quad \mbox{\rm and}  \quad  
\left| F\Big(\frac{1}{2},\frac{1}{2};1; x  \Big) \right|  \leq C \big(1-\ln(1-x)) \quad \mbox{\rm for all}  \quad x \in [\varepsilon ,1) \,.
\eeq
This allows to prove the estimate for the integral over the second zone
\beq 
&  &
\int_{  (z,r) \in Z_2(\varepsilon ,z) }  dr  \frac{1}{     ((z-1)^2-r^2 ) \sqrt{(z+1)^2-r^2}    } \nonumber \\
&  &
\times \left| (1-z^{2}+r^2 ) 
F\Big(-\frac{1}{2},\frac{1}{2};1; \frac{ (z-1)^2 -r^2   }{ (z+1)^2 -r^2 }   \Big) 
+   2  (z -1) F\Big(\frac{1}{2},\frac{1}{2};1; \frac{ (z-1)^2 -r^2   }{ (z+1)^2 -r^2 }   \Big)    \right|  \nonumber \\
& \leq   &
\label{10.4}
C  \qquad \mbox{\rm for all}  \qquad z \in [1,\infty) \,.
\eeq 
Indeed, for the argument of the hypergeometric functions we have
\[
\varepsilon \leq \frac{ (z-1)^2 -r^2   }{ (z+1)^2 -r^2 }  = 1- \frac{  4z  } { (z  + 1)^2 -r ^2} < 1,\quad 
\frac{  4z  } { (z  + 1)^2 -r ^2} < 1- \varepsilon \quad \mbox{\rm for all}  \quad (z,r) \in Z_2(\varepsilon ,z) \,.
\]
Hence,
\be 
\left| F\Big(\frac{1}{2},\frac{1}{2};1; \frac{ (z-1)^2 -r^2   }{ (z+1)^2 -r^2 }    \Big) \right|   
  \leq   
C \left(1-\ln\frac{  4z  } { (z  + 1)^2 -r ^2} \right) 
  \leq  
\label{10.8}
C \left(1+ \ln  z     \right) \quad \mbox{\rm for all}  \quad  (z,r) \in Z_2(\varepsilon ,z)  . 
\ee 
To prove (\ref{10.4}) we have to estimate the following two integrals
\beqst 
A_4 
& :=   &
\int_{  (z,r) \in Z_2(\varepsilon ,z) }  \frac{1}{     ((z-1)^2-r^2 ) \sqrt{(z+1)^2-r^2}    }  \left| (1-z^{2}+r^2 ) 
 \right|  dr \,, \\
A_5 
& :=   & 
\int_{  (z,r) \in Z_2(\varepsilon ,z) }   \frac{1}{     ((z-1)^2-r^2 ) \sqrt{(z+1)^2-r^2}    } 
\left| (z -1) \left(1+ \ln  z     \right)   \right| dr \,.
\eeqst 
We apply (\ref{10.4n}) to $A_4$ and obtain
\[ 
A_4 
  \leq   
C_\varepsilon \int_{  (z,r) \in Z_2(\varepsilon ,z) }    \frac{1}{    \sqrt{(z+1)^2-r^2}    }   dr 
  \leq    
C_\varepsilon \int_{  0 }^{z-1}    \frac{1}{    \sqrt{(z+1)^2-r^2}    }   dr 
  \leq   
C_\varepsilon \,, 
\] 
while
\beqst 
A_5 
& \leq    & 
C_\varepsilon  (z -1) \left(1+ \ln  z     \right) \int_{  (z,r) \in Z_2(\varepsilon ,z) }  \frac{1}{ ((z+1)^2-r^2 ) \sqrt{(z+1)^2-r^2}    }  dr \\
& \leq    & 
C_\varepsilon  (z -1) \left(1+ \ln  z     \right) \int_{  0 }^{z-1}  \frac{1}{ ((z+1)^2-r^2 )^{3/2}   }  dr \\
& \leq    & 
C_\varepsilon  (z -1)^2 \left(1+ \ln  z     \right)   \frac{1}{ \sqrt{z}(z+1)^2 }  \\
& \leq    & 
C_\varepsilon   \left(1+ \ln  z     \right)   \frac{1}{ \sqrt{z} }  \,.
\eeqst 
Thus, (\ref{10.4}) is proven.
Lemma is proven. \hfill $\square$

\section{Some Estimates of the Kernels $K_0$ and $K_1$. $L^p-L^q$ Decay Estimates for  Equation with $n=1$ and without Source Term}
\label{S11}
\setcounter{equation}{0}

\begin{theorem}
Let $u=u(x,t) $  be a solution of the Cauchy problem 
\[
u_{tt} - e^{2t}u_{xx} =0\, ,\qquad u(x,0) = \varphi_0  (x) \,, \qquad u_t(x,0) =\varphi_1  (x)\,  ,
\]
with $\varphi_0  $, $\varphi_1  \in C_0^\infty ({\mb R})$. 
If $\rho \in (1,2) $, then
\beqst
\| u(x,t) \| _{L^q({\mb R}_x)} 
& \leq &
e ^{-\frac{t}{2}} 
\| \varphi_0   (x)  \| _{L^q({\mb R}_x)}  
+  C_\rho  (e^t-1)^{  \frac{1}{\rho}  } e^{-t}  \| \varphi_0   (x)  \| _{L^p({\mb R}_x)}  \\
&  &
 +\,
C (1+  t)  (e^t  - 1)^{\frac{1 }{\rho}-1} (1  - e^{ -t}) \|  \varphi_1 (x) \| _{L^p({\mb R}_x)} 
\,,
\eeqst
for all $ t \in (0,\infty)$. Here $1 <p<\rho '$, $ \frac{1}{q} = \frac{1}{p} - \frac{1}{\rho '} $, $ \frac{1}{\rho } + \frac{1}{\rho '}= 1$.  If $\rho =1$, then 
\be 
\label{10.1LpLq}   
\|  u(x,t) \|_{ { L}^{  q} ({\mathbb R})  } 
  \le     
 C \left(  \|  \varphi_0   (x)   \|_{ { L}^{  q} ({\mathbb R})  } +  (1+  t) \|  \varphi_1   (x)   \|_{ { L}^{  q} ({\mathbb R})  } \right)
\, \quad \mbox{for all} \,\quad t \in (0,\infty) .
\ee 
\end{theorem}

For $\rho =1 $ we apply Proposition~\ref{P9.3}. 
To prove this theorem for  $\rho \not=1 $  we need some auxiliary estimates for the kernels $K_0$ and $K_1$. We start with the case of $\varphi_0  =0 $,
 where the kernel $K_1$
appears.  The application of Theorem~\ref{T1.3} and Young's inequality lead  to 
\beqst
\| u(x,t) \| _{L^q({\mb R}_x)} 
& \leq  &
2\Bigg( \int_{0}^{ e^t-1 } 
| K_1(x,t)|^\rho  dx \Bigg)^{1/\rho } \|  \varphi (x) \| _{L^p({\mb R}_x)} ,
\eeqst
where $1 <p<\rho '$, $ \frac{1}{q} = \frac{1}{p} - \frac{1}{\rho '} $, $ \frac{1}{\rho } + \frac{1}{\rho '}= 1$.  Now we have to estimate
the integral $ \left( \int_{0}^{ e^t-1 } 
| K_1(x,t)|^\rho  dx \right)^{1/\rho }$.
\begin{proposition}
We have
\beqst
\left( \int_{0}^{ e^t-1 } 
| K_1(x,t)|^\rho  dx \right)^{1/\rho } 
& \leq  &
 C (1+  t)  (e^t  - 1)^{1/\rho -1} (1  - e^{ -t}) \quad \mbox{for all} \,\quad t \in (0,\infty) \,.
\eeqst
\end{proposition}
\medskip

\noindent
{\bf Proof.} One can write
\beqst
\left( \int_{0}^{ e^t-1 } 
| K_1(x,t)|^\rho  dx \right)^{1/\rho } 
& \leq  &
\left( \int_{0}^{ e^t-1 } 
\left| \frac{1}{  \sqrt{(1+e^t)^2-x^2 } }    F\Big(\frac{1}{2},\frac{1}{2};1; 
\frac{  (e^t  - 1 )^2-x^2}
{  (e^t  + 1     )^2- x ^2} \Big)  \right|^\rho  dx \right)^{1/\rho } \,.
\eeqst
Denote $z:= e^t >1 $ and consider the first integral $
\dsp \int_{0}^{ z-1 } 
\left| \frac{1}{  \sqrt{(1+z)^2-x^2 } }    F\Big(\frac{1}{2},\frac{1}{2};1; 
\frac{  (z  - 1  )^2 - x^2}
{  (z  + 1)^2-x^2} \Big)  \right|^\rho  dx $ of the right-hand side. 
According to Lemma \ref{L10.2} we obtain that for all $z>1$ the following estimate 
\[ 
\int_{0}^{ z-1 } 
\left| \frac{1}{  \sqrt{(1+z)^2-x^2 } }    F\Big(\frac{1}{2},\frac{1}{2};1; 
\frac{  (z  - 1  )^2 - x^2}
{  (z  + 1)^2-x^2} \Big)  \right|^\rho  dx  
  \leq  
C (1+ \ln   z )^\rho (z  - 1) (z  + 1)^{-\rho } F\Big(\frac{1}{2},\frac{\rho }{2};\frac{3}{2}; \frac{  (z  - 1)^2  }{ (z  + 1)^2  } \Big)
\]
is fulfilled, provided that $1 <p<\rho '$, $ \frac{1}{q} = \frac{1}{p} - \frac{1}{\rho '} $, $ \frac{1}{\rho } + \frac{1}{\rho '}= 1$. 
In particular, if $\rho < 2$, then 
\beqst
\left( \int_{  0}^{ z  - 1} 
((z  + 1)^2  - r^2  )^{-\frac{\rho }{2}}  
F\Big(\frac{1}{2},\frac{1}{2};1; 
\frac{ (z  - 1 )^2 -  r^2   }
{ (z  + 1 )^2  - r^2} \Big)^\rho  d r \right)^{1/\rho }
& \leq &
C (1+ \ln   z )  (z  - 1)^{1/\rho } (z  + 1)^{-1 } 
\,.
\eeqst
Proposition is proven. \hfill $\square$

Thus, the theorem in the case of $\varphi_0  =0 $ is proven.
\medskip

\noindent
Now we turn to the case of $\varphi_1  =0 $,
 where the kernel $K_0$
appears.  The application of Theorem~\ref{T1.3} leads to 
\beqst
\| u(x,t) \| _{L^q({\mb R}_x)} 
& \leq  &  
e ^{-\frac{t}{2}} 
\| \varphi_0   (x)  \| _{L^q({\mb R}_x)}  
+  \, \left\| \int_{ 0}^{e^t-1}\left[ 
\varphi_0   (x - z) 
+     \varphi_0   (x  + z)  \right] K_0(z,t)\, dz  \right\| _{L^q({\mb R}_x)} \,.
\eeqst 
Similarly to the case of the second datum we arrive at
\beqst
\| u(x,t) \| _{L^q({\mb R}_x)} 
& \leq  &  
e ^{-\frac{t}{2}} 
\| \varphi_0   (x)  \| _{L^q({\mb R}_x)}  
+  \, 
\| \varphi_0   (x )  \| _{L^p({\mb R}_x)} \left( \int_{0}^{ e^t-1 } | K_0(r,t)|^\rho  dr \right)^{1/\rho }\,.
\eeqst 
The next proposition gives an estimate for  the integral $ \left( \int_{0}^{ e^t-1 } 
| K_0(r,t)|^\rho  dr \right)^{1/\rho }$.
\begin{proposition}
\label{P11.4}
Let $1 <p<\rho '$, $ \frac{1}{q} = \frac{1}{p} - \frac{1}{\rho '} $, $ \frac{1}{\rho } + \frac{1}{\rho '}= 1$, and $ \rho \in [1,2)$.  We have
\beqst
\left( \int_{0}^{ e^t-1 } 
| K_0(r,t)|^\rho  dr \right)^{1/\rho } 
& \leq  &
 C_\rho  (e^t-1)^{  \frac{1}{\rho}  }(e^t+1)^{-1 } \quad \mbox{ for all} \quad t \in (0,\infty)  \,.
\eeqst
\end{proposition}
\medskip

\noindent
{\bf Proof.} 
We turn to the integral ($z= e^t>1$)
\beqst
I_2
& := &
\left( \int_{0}^{ z-1 } 
\left|  \frac{1}{  ((z-1)^2-r^2 ) \sqrt{(z+1)^2-r^2}      }  \right|^\rho \right. \\
&  &
\left. \times  \left|  (1-z^{2}+r^2 ) 
F\Big(-\frac{1}{2},\frac{1}{2};1;  \frac{ (z-1)^2 -r^2   }{ (z+1)^2 -r^2 }   \Big) 
+   2  (z -1) F\Big(\frac{1}{2},\frac{1}{2};1; \frac{ (z-1)^2 -r^2   }{ (z+1)^2 -r^2 }  \Big)          \right|^\rho  dr \right)^{1/\rho } \,.
\eeqst
The formula (\ref{4.12})  describes the behavior of those functions at the neighbourhood of zero.
Hence, if $\varepsilon  >0$ is small, the for all $z$ and $r$ such that (\ref{9.7}) holds, 
one has (\ref{9.8}).
Consider therefore two zones, $Z_1(\varepsilon ,z) $ and $ Z_2(\varepsilon ,z)$, defined in  (\ref{9.10}) and (\ref{9.11}),
respectively.
We  split integral into two parts:
\beqst
\int_{ 0}^{e^t-1} \left|   K_0(r,t)  \right|^\rho   dr 
& = &
\int_{ (z,r) \in Z_1(\varepsilon ,z)   }   \left|   K_0(r,t)  \right|^\rho   dr 
+ \int_{ (z,r) \in Z_2(\varepsilon ,z)  }  \left|   K_0(r,t)  \right|^\rho   dr \,.
\eeqst 
In the proof of Lemma~\ref{L10.3} the   relation (\ref{9.15}) was checked in the first zone.  
If $1\leq  z \leq M$ with some constant $M$, then the argument of the hypergeometric functions is bounded,
\be
\label{argument}
\frac{ (z-1)^2 -r^2   }{ (z+1)^2 -r^2 } \leq C_M <1 \quad \mbox{for all} \quad r \in (0,z-1),
\ee
and we obtain 
\beqst
\left( \int_{ 0}^{z-1} \left|   K_0(r,t)  \right|^\rho   dr \right)^{1/\rho } 
& \leq  & 
C \left( \int_{0}^{ z-1 } 
\left|  \frac{1}{   \sqrt{(z+1)^2-r^2}      }  \right|^\rho  dr \right)^{1/\rho }  \\
& \leq  & 
C \left( (z-1)(z+1)^{-\rho }F\Big(\frac{1}{2},\frac{\rho }{2};\frac{3}{2}; \frac{ (z-1)^2     }{ (z+1)^2   }   \Big) \right)^{1/\rho } \\
& \leq  & 
C   (z-1)^{1/\rho } (z+1)^{-1 }   \,.
\eeqst
Thus, we can restrict ourselves to the  case of large $   z \geq M$ in both zones.

Consider therefore for $\rho \in (1,2)$ the integrals over the first zone
\beqst
A_6
& := &
\int_{(z,r)\in Z_1(\varepsilon ,z)} \left| \frac{1}{   \sqrt{(z+1)^2-r^2}    } \right|^\rho  dr  
 \,\, \leq    \,\,
\int_{0}^{z-1}   \left|\frac{1}{   \sqrt{(z+1)^2-r^2}    } \right|^\rho dr \\
& \leq  & 
C  (z-1)(z+1)^{-\rho }F\Big(\frac{1}{2},\frac{\rho }{2};\frac{3}{2}; \frac{ (z-1)^2     }{ (z+1)^2   }   \Big)  \\
& \leq  &
C  (z-1)(z+1)^{-\rho }
\eeqst 
and
\beqst
A_7
& := &
\int_{(z,r)\in Z_1(\varepsilon ,z)} \left| \frac{1}{    ((z-1)^2-r^2 ) \sqrt{(z+1)^2-r^2}    } 
\big(z^{2}+2z -3-r^2  \big)\left( \frac{ (z-1)^2 -r^2   }{ (z+1)^2 -r^2 }\right)^2\right|^\rho  dr \\
& \leq   &
\int_{0}^{z-1}   \left|\frac{1}{   \sqrt{(z+1)^2-r^2}    } \right|^\rho dr \\
& \leq  & 
C  (z-1)(z+1)^{-\rho }\,.
\eeqst 
In the second zone  for the argument of the hypergeometric functions we have
\[
\varepsilon \leq \frac{ (z-1)^2 -r^2   }{ (z+1)^2 -r^2 }  = 1- \frac{  4z  } { (z  + 1)^2 -r ^2} < 1,\quad 
\frac{  4z  } { (z  + 1)^2 -r ^2} < 1- \varepsilon \quad \mbox{\rm for all}  \quad (z,r) \in Z_2(\varepsilon ,z) \,,
\]
and
\[
\frac{ 1 } { (z-1)^2 -r^2   } \leq  \frac{ 1 }{ \varepsilon[(z+1)^2 -r^2] },\qquad   0 \leq r \leq z-1 \,.
\]
Hence,
\[
\left| F\Big(\frac{1}{2},\frac{1}{2};1; \frac{ (z-1)^2 -r^2   }{ (z+1)^2 -r^2 }    \Big) \right|   
  \leq  
C \left(1-\ln\frac{  4z  } { (z  + 1)^2 -r ^2} \right) \leq 
C \left(1+ \ln  z     \right) \quad \mbox{\rm for all}  \quad  (z,r) \in Z_2(\varepsilon ,z)  .
\]
We have to estimate the following two integrals
\beqst 
A_8 
& :=   &
\int_{  (z,r) \in Z_2(\varepsilon ,z) }  \left| \frac{1}{     ((z-1)^2-r^2 ) \sqrt{(z+1)^2-r^2}    }    (z^{2}-1- r^2 ) 
 \right|^\rho   dr \,, \\
A_9 
& :=   & 
\int_{  (z,r) \in Z_2(\varepsilon ,z) }   \left|  \frac{1}{     ((z-1)^2-r^2 ) \sqrt{(z+1)^2-r^2}    } 
 (z -1) \left(1+ \ln  z     \right)   \right|^\rho dr \,.
\eeqst 
We apply (\ref{10.4n}) and obtain 
\beqst 
A_8 
& \leq    &
\int_{  (z,r) \in Z_2(\varepsilon ,z) }  \left| \frac{1}{     ((z+1)^2-r^2 ) \sqrt{(z+1)^2-r^2}    }    (z^{2}-1- r^2 ) 
 \right|^\rho   dr \\
& \leq   &
\int_{0}^{z-1}   \left|\frac{1}{   \sqrt{(z+1)^2-r^2}    } \right|^\rho dr \\
& \leq  & 
C  (z-1)(z+1)^{-\rho }F\Big(\frac{1}{2},\frac{\rho }{2};\frac{3}{2}; \frac{ (z-1)^2     }{ (z+1)^2   }   \Big)  \\
& \leq  &
C  (z-1)(z+1)^{-\rho } \,,
\eeqst 
while
\beqst 
A_9 
& \leq   & 
C_\varepsilon (z -1)^\rho \left(1+ \ln  z     \right)^\rho
\int_{  (z,r) \in Z_2(\varepsilon ,z) }     ( (z+1)^2-r^2)^{-3\rho/2} dr \\
& \leq   & 
C_\varepsilon (z -1)^\rho \left(1+ \ln  z     \right)^\rho
(z -1) (z +1)^{-3\rho } F\Big(\frac{1}{2},\frac{3\rho }{2};\frac{3}{2}; \frac{ (z-1)^2     }{ (z+1)^2   }   \Big)\\
& \leq   & 
 C  (z-1)(z+1)^{-\rho }  \,.
\eeqst 
The proposition is proven. \hfill $\square $

\section{$L^p-L^q$ Decay Estimates for the Equation with  Source, $n>1$}
\label{S12}
\setcounter{equation}{0}

For the wave equation the Duhamel's principle allows to reduce the case of source term to the case of the Cauchy problem without
source term and consequently to derive the $L^p-L^q$-decay estimates for the equation. For (\ref{G-S}) the Duhamel's principle
is not applicable straightforward and we have to appeal to the representation formula of Theorem~\ref{T1.5}. In fact, one can regard that
formula as an expansion of the two-stage Duhamel's principle.
In this section we consider the Cauchy problem (\ref{1.26}) for the equation with the source term with zero initial data.   
\begin{theorem} 
\label{T11.1}
Let $u=u(x,t)$ be solution of the Cauchy problem (\ref{1.26}). Then for $n>1$ one has the following decay estimate
\beqst
&  &
\| (-\bigtriangleup )^{-s} u(x,t) \|_{ { L}^{  q} ({\mathbb R}^n)  } \\
& \le &  
C \int_{ 0}^{t} \|  f(x, b)  \|_{ { L}^{p} ({\mathbb R}^n)  } db
  \int_{ 0}^{ e^t- e^b} dr \,\, r^{ 2s-n(\frac{1}{p}-\frac{1}{q}) }  \frac{1}{\sqrt{(e^t  + e^b)^2    - r^2}}
F\left(\frac{1}{2},\frac{1}{2};1; 
\frac{  (e^t  - e^b )^2  -r^2 }
{  (e^t  + e^b )^2  -r^2  } \right) 
\eeqst
provided that $s \ge 0$, $1<p \leq 2$, $\frac{1}{p}+ \frac{1}{q}=1$, $\frac{1}{2} (n+1)\left( \frac{1}{p} - \frac{1}{q}\right) \leq 
2s \leq n \left( \frac{1}{p} - \frac{1}{q}\right) $, $-1 + n \left( \frac{1}{p} - \frac{1}{q}\right)< 2s$.
\end{theorem}
\medskip

\noindent
{\bf Proof.} According to the representation  (\ref{1.29}) and to the results of \cite{Brenner,Pecher} for the wave equation, we have
\beqst
&  &
\| (-\bigtriangleup )^{-s} u(x,t) \|_{ { L}^{  q} ({\mathbb R}^n)  } \\
& \le &  
C\int_{ 0}^{t} db
  \int_{ 0}^{ e^t- e^b}  \| (-\bigtriangleup )^{-s} v(x,r ;b)  \|_{ { L}^{  q} ({\mathbb R}^n)  }   \frac{1}{\sqrt{(e^t  + e^b)^2    - r^2}}
F\left(\frac{1}{2},\frac{1}{2};1; 
\frac{  (e^t  - e^b )^2  -r^2 }
{  (e^t  + e^b )^2  -r^2  }  \right) dr \\ 
& \le &  
C \int_{ 0}^{t} db\, \|  f(x, b)  \|_{ { L}^{p} ({\mathbb R}^n)  }   
  \int_{ 0}^{ e^t- e^b} r^{ 2s-n(\frac{1}{p}-\frac{1}{q}) }  \frac{1}{\sqrt{(e^t  + e^b)^2    - r^2}}
F\left(\frac{1}{2},\frac{1}{2};1; 
\frac{  (e^t  - e^b )^2  -r^2 }
{  (e^t  + e^b )^2  -r^2  }  \right) dr\,.
\eeqst
The theorem is proven. \hfill $\square$

We are going to transform the estimate of the last theorem to more cosy form. 
To this aim we  estimate  for $2s-n(\frac{1}{p}-\frac{1}{q})>-1 $ the  last integral of the right hand side. 
If we replace $e^t/ e^b >  1$ with $z >  1$, then the integral will be simplified.
\begin{lemma}
\label{L12.2}
Assume that $0 \ge 2s-n(\frac{1}{p}-\frac{1}{q})>-1 $. Then 
\[
  \int_{ 0}^{ z- 1}  r^{ 2s-n(\frac{1}{p}-\frac{1}{q}) }  \frac{1}{\sqrt{(z  + 1)^2    - r^2}}
F\left(\frac{1}{2},\frac{1}{2};1; 
\frac{  (z  - 1)^2  -r^2 }
{  (z  + 1)^2  -r^2  }  \right)  \, dr  
 \leq   
C     z^{-1}(z-1)^{ 1+  2s-n(\frac{1}{p}-\frac{1}{q}) } (1+ \ln z ),
\]
for all $z>1$.
\end{lemma}
\medskip

\noindent
{\bf Proof.}   If $1 <  z \leq M$ with some constant $M$, then the argument $ $ of the hypergeometric functions is bounded, see (\ref{argument}),
and
\beqst
&  &
  \int_{ 0}^{ z- 1}  r^{ 2s-n(\frac{1}{p}-\frac{1}{q}) }  \frac{1}{\sqrt{(z  + 1)^2    - r^2}}
F\left(\frac{1}{2},\frac{1}{2};1; 
\frac{  (z  - 1)^2  -r^2 }
{  (z  + 1)^2  -r^2  }  \right)  \, dr \\
 &  \leq &
C _M     (z-1)^{1+   2s-n(\frac{1}{p}-\frac{1}{q}) } , \quad \mbox{\rm for all } \quad 1 <  z \leq M\,.
\eeqst
Hence, we can restrict ourselves to the case of large $z$, that is $z\geq M$. In particular, we choose $M>6$ and split integral into two parts:
\beqst
&  &
  \int_{ 0}^{ z- 1}  r^{ 2s-n(\frac{1}{p}-\frac{1}{q}) }  \frac{1}{\sqrt{(z  + 1)^2    - r^2}}
F\left(\frac{1}{2},\frac{1}{2};1; 
\frac{  (z  - 1)^2  -r^2 }
{  (z  + 1)^2  -r^2  }  \right)  \, dr \\
 & = &
  \int_{ 0}^{ \sqrt{(z  + 1)^2-8z}}  r^{ 2s-n(\frac{1}{p}-\frac{1}{q}) }  \frac{1}{\sqrt{(z  + 1)^2    - r^2}}
F\left(\frac{1}{2},\frac{1}{2};1; 
1- \frac{  4z  }
{ (z  + 1)^2 -r ^2} \right)  \, dr \\
&     &
+   \int_{ \sqrt{(z  + 1)^2-8z}}^{ z- 1}  r^{ 2s-n(\frac{1}{p}-\frac{1}{q}) }  \frac{1}{\sqrt{(z  + 1)^2    - r^2}}
F\left(\frac{1}{2},\frac{1}{2};1; 
1- \frac{  4z  }
{ (z  + 1)^2 -r ^2}\right)  \, dr \,.
\eeqst
For the second part we have (\ref{r}) and $z\ge M >6$, 
then (\ref{9.2})
and (\ref{9.3}) 
imply 
\beqst
&  &
\int_{\sqrt{(z  + 1)^2-8z}  }^{ z-1}  r^{ 2s-n(\frac{1}{p}-\frac{1}{q}) }   \frac{1}{\sqrt{(z  + 1)^2  -r^2}}  
F\left(\frac{1}{2},\frac{1}{2};1; 1- \frac{  4z  }{ (z  + 1)^2 -r ^2} \right) \, dr  \\
&  \le &
C \int_{0}^{ z-1}  r^{ 2s-n(\frac{1}{p}-\frac{1}{q}) }   \frac{1}{\sqrt{(z  + 1)^2  -r^2}}   \, dr \\
& \le  &
C  (1+z)^{2 s - n \left(  \frac{1}{p}- \frac{1}{q} \right)}   \quad \mbox{\rm for all} \quad z\ge M >6\,.
\eeqst
For the first integral \, 
$ r \le \sqrt{(z  + 1)^2-8z} $  \, and  \, $ z\ge M >6 $ \,  imply  \,  $8z \leq (z  + 1)^2 - r^2 $. It follows  
\[
\left| F\left(\frac{1}{2},\frac{1}{2};1; 1-
\frac{  4z  }
{ (z  + 1)^2 -r ^2} \Bigg) \right| \le C \left| \ln \Bigg( \frac{  4z  }
{ (z  + 1)^2 -r ^2} \right)  \right|  \le C   (1+ \ln   z )  \,.       
\]
Then we obtain
\beqst
&  &
  \int_{ 0}^{ \sqrt{(z  + 1)^2-8z}}  r^{ 2s-n(\frac{1}{p}-\frac{1}{q}) }  \frac{1}{\sqrt{(z  + 1)^2    - r^2}}
F\left(\frac{1}{2},\frac{1}{2};1; 
1- \frac{  4z  }
{ (z  + 1)^2 -r ^2} \right)  \, dr \\
&   \leq   &
 C   (1+ \ln   z ) \int_{ 0}^{z-1}  r^{ 2s-n(\frac{1}{p}-\frac{1}{q}) }  \frac{1}{\sqrt{(z  + 1)^2    - r^2}} \, dr \\
&   \leq   &
 C   (1+ \ln   z )  (1+z)^{2 s - n \left(  \frac{1}{p}- \frac{1}{q} \right)}   \,.
\eeqst
Lemma is proven.
\hfill $\square$

\begin{corollary}
\label{C12.3}
Let $u=u(x,t)$ be solution of the Cauchy problem (\ref{1.26}). Then for $n>1$ one has the following decay estimate  
 \beq
 \label{12.1}
\| (-\bigtriangleup )^{-s} u(x,t) \|_{ { L}^{  q} ({\mathbb R}^n)  } 
& \le &  
C   e^{t \big(  2s-n(\frac{1}{p}-\frac{1}{q}) \big) } \int_{ 0}^{t} \|  f(x, b)  \|_{ { L}^{p} ({\mathbb R}^n)  }  
(1+  t -b  )\,db 
\eeq
provided that $s \ge 0$, $1<p \leq 2$, $\frac{1}{p}+ \frac{1}{q}=1$, $\frac{1}{2} (n+1)\left( \frac{1}{p} - \frac{1}{q}\right) \leq 
2s \leq n \left( \frac{1}{p} - \frac{1}{q}\right) $, $-1+ n \left( \frac{1}{p} - \frac{1}{q}\right)< 2s  $.
\end{corollary}
\medskip
 
\noindent
{\bf Proof.} Indeed, from Theorem~\ref{T11.1} we derive 
 \beqst
\| (-\bigtriangleup )^{-s} u(x,t) \|_{ { L}^{  q} ({\mathbb R}^n)  } & \le &  
C \int_{ 0}^{t} \|  f(x, b)  \|_{ { L}^{p} ({\mathbb R}^n)  }  e^{b(  2s-n(\frac{1}{p}-\frac{1}{q})) }db    \\
&  &
\times   \int_{ 0}^{ e^{t -b}- 1} dl \,\,   
\frac{l^{ 2s-n(\frac{1}{p}-\frac{1}{q}) } }{\sqrt{(e^{t-b}  + 1)^2    - r^2}}  
F\left(\frac{1}{2},\frac{1}{2};1; 
\frac{  (e^{t -b}  - 1 )^2  -l^2 }
{  (e^{t -b}  + 1 )^2  -l^2  } \right) .
\eeqst
 Next we apply Lemma~\ref{L12.2} with $z= e^{t -b}$ and arrive at (\ref{12.1}). Corollary is proven.
\hfill $\square$

\section{$L^p-L^q$ Decay Estimates for the Equation without Source, $n>1$}
\label{S13}
\setcounter{equation}{0}

The $L^p-L^q$-decay estimates for the energy of the solution of the Cauchy problem for the 
wave equation without source can be proved 
by the representation formula, $L_1-L_ \infty$ and  
$L_2-L_2$ estimates, interpolation argument. (See, e.g., \cite[Theorem~2.1]{Racke}.) 
There is also a proof of the $L^p-L^q$-decay estimates for the solution itself, that is based on the
microlocal consideration and dyadic decomposition of the phase space.  (See, e.g., \cite{Brenner,Pecher}.) 
The last one was applied in  \cite{Galstian2001,Galstian} to the equation (\ref{G-S}) and its 
result is given by (\ref{1.16}) that contains some loss of regularity.  
The application of the first approach includes the step with the Granwall inequality that brings some inaccuracy in the result.  
To avoid the  loss of regularity and obtain more sharp estimates  we 
appeal to the representation formula provided by Theorem~\ref{T1.6}.

\begin{theorem} 
\label{T13.2}
The solution  $u=u(x,t)$ of the Cauchy problem (\ref{1.30CP})  satisfies the following $L^p-L^q$ estimate
\begin{eqnarray*} 
\| (-\bigtriangleup )^{-s} u(x,t) \|_{ { L}^{  q} ({\mathbb R}^n)  } 
& \leq  &
 C(e^ t-1)^{2s-n(\frac{1}{p}-\frac{1}{q}) }  \left\{ \| \varphi_0  (x) \|_{ { L}^{p} ({\mathbb R}^n)  } 
 +  \|\varphi_1  \|_{ { L}^{p} ({\mathbb R}^n)  }(1+ t ) (1-e^{-t}) \right\}
\end{eqnarray*}
for all $t \in (0,\infty)$, provided that $s \ge 0$, $1<p \leq 2$, $\frac{1}{p}+ \frac{1}{q}=1$, $\frac{1}{2} (n+1)\left( \frac{1}{p} - \frac{1}{q}\right) \leq 
2s \leq n \left( \frac{1}{p} - \frac{1}{q}\right) $, $-1+ n \left( \frac{1}{p} - \frac{1}{q}\right)< 2s  $. 
\end{theorem}
\medskip

\noindent
{\bf Proof.} We start with the case of $\varphi _0=0$. Due  to Theorem~\ref{T1.6} for the solution $u=u (x,t)$ of the Cauchy 
problem (\ref{1.30CP}) with  $\varphi _0=0$ and to the results of \cite{Brenner,Pecher} 
we have:
\begin{eqnarray*}
\| (-\bigtriangleup )^{-s} u(x,t) \|_{ { L}^{  q} ({\mathbb R}^n)  } 
 \leq 
C \|\varphi_1  \|_{ { L}^{p} ({\mathbb R}^n)  }\int_{0}^{e^t-1}   
 r^{ 2s-n(\frac{1}{p}-\frac{1}{q}) }  \left| K_1(r,t)\right|   \, dr \,.  
\end{eqnarray*}
To continue we need the following lemma.
\begin{lemma}
\label{L13.1}
The following inequality holds   
\begin{eqnarray*}
\int_{0}^{z-1}   
 r^{ 2s-n(\frac{1}{p}-\frac{1}{q}) }  \left| K_1(r,t)\right|   \, dr  \leq C    (1+ \ln z ) z^{-1}(z-1)^{ 1+  2s-n(\frac{1}{p}-\frac{1}{q})  }
\quad \mbox{\rm for all}\quad z>1.
\end{eqnarray*}
\end{lemma}
\medskip

\noindent
{\bf Proof.} In fact, we have to estimate the integral: 
\begin{eqnarray*}
I_3
& := &
 \int_{0}^{z-1}   
r^{ 2s-n(\frac{1}{p}-\frac{1}{q}) }  \frac{1}{  \sqrt{(z+1)^2-r^2 } }    F\Big(\frac{1}{2},\frac{1}{2};1; 
\frac{  (z  - 1)^2-  r ^2}
{  (z  + 1     )^2-r^2} \Big)    \, dr    \,,
\end{eqnarray*}
where 
$ z= e^t$. The estimate for $I_3$ is given by Lemma~\ref{L12.2}. 
Thus, for the case of $\varphi _0=0$ the theorem is proven.
\medskip

Next we turn to the case of $\varphi _1=0$.    Due  to Theorem~\ref{T1.6} for the solution $u=u (x,t)$ of the Cauchy 
problem (\ref{1.30CP}) with  $\varphi _1=0$ and to the results of \cite{Brenner,Pecher} 
we have: 
\begin{eqnarray*} 
&  &
\| (-\bigtriangleup )^{-s} u(x,t) \|_{ { L}^{  q} ({\mathbb R}^n)  }  \\
& \leq & 
C e ^{-\frac{t}{2}} (e^t-1) ^{ 2s-n(\frac{1}{p}-\frac{1}{q}) } \| \varphi_0  (x) \|_{ { L}^{p} ({\mathbb R}^n)  }
+ \, C\| \varphi_0  (x) \|_{ { L}^{p} ({\mathbb R}^n)  }\int_{ 0}^{e^t-1}  r^{ 2s-n(\frac{1}{p}-\frac{1}{q}) }
|K_0(r,t)|  \,  dr.
\end{eqnarray*}
The following proposition gives the remaining   estimate for the last integral, 
$\int_{ 0}^{z-1}  r^{ 2s-n(\frac{1}{p}-\frac{1}{q}) }
|K_0(r,t)|  \,  dr$, 
and completes the proof of the theorem.
\begin{proposition}
\label{P13.4}
If  $2s-n(\frac{1}{p}-\frac{1}{q})>-1 $, then
\begin{eqnarray*} 
\int_{ 0}^{z-1}  r^{ 2s-n(\frac{1}{p}-\frac{1}{q}) }
|K_0(r,t)|  \,  dr \leq C z^{-1}(z-1)^{1+ 2s-n(\frac{1}{p}-\frac{1}{q}) }  \quad \mbox{ for all} \quad z>1.
\end{eqnarray*}
\end{proposition}
\medskip

\noindent
{\bf Proof.} We follow the arguments have been used in the proof of Proposition~\ref{P11.4}. If $1\leq  z \leq M$ with some constant $M$, then the argument  of the hypergeometric functions is bounded (\ref{argument}), and we have
\beqst
\int_{ 0}^{z-1}  r^{ 2s-n(\frac{1}{p}-\frac{1}{q}) }\left|   K_0(r,t)  \right|  dr  
& \leq &  
C  \int_{0}^{ z-1 } 
 \frac{1}{   \sqrt{(z+1)^2-r^2}      }  r^{ 2s-n(\frac{1}{p}-\frac{1}{q}) }  dr  \\
& \leq  & 
C_M  (z-1)^{1+ 2s-n(\frac{1}{p}-\frac{1}{q}) }, \qquad 1 <  z \leq M \,.
\eeqst
Thus, we can restrict ourselves to the  case of large $   z \geq M$ in both zones $Z_1(\varepsilon ,z) $ and $ Z_2(\varepsilon ,z)$, 
defined in  (\ref{9.10}) and (\ref{9.11}),
respectively. In the first zone we have (\ref{9.15}). 
Consider therefore the following inequalities,
\beqst
A_{10}
& := &
\int_{ (z,r) \in Z_1(\varepsilon ,z)  }  r^{ 2s-n(\frac{1}{p}-\frac{1}{q}) }  \frac{1}{     \sqrt{(z+1)^2-r^2}    }  \, dr \\
& \leq   &
C  \int_{0}^{ z-1 } 
 r^{ 2s-n(\frac{1}{p}-\frac{1}{q}) }   \frac{1}{   \sqrt{(z+1)^2-r^2}      } \,dr    \\
& \leq   &
C  z^{ 2s-n(\frac{1}{p}-\frac{1}{q}) }   \quad \mbox{\rm for all}  \quad z \in [1,\infty)  \,,\\
A_{11}
& := &
\int_{  (z,r) \in Z_1(\varepsilon ,z) } r^{ 2s-n(\frac{1}{p}-\frac{1}{q}) }  
 \frac{1}{   \sqrt{(z+1)^2-r^2}    } \frac{ \big|3-z^{2} -2   z  +r^2 \big|   }{ (z+1)^2 -r^2 }  
 dr \\
& \leq   &
C\int_{ 0}^{z-1}  r^{ 2s-n(\frac{1}{p}-\frac{1}{q}) }  \frac{1}{ \sqrt{(z+1)^2-r^2}    }    dr \\
& \leq   &
C  z^{ 2s-n(\frac{1}{p}-\frac{1}{q}) }   \quad \mbox{\rm for all}  \quad z \in [1,\infty) \,,
\eeqst 
and 
\beqst
A_{12}
& := &
\int_{  (z,r) \in Z_1(\varepsilon ,z) }   r^{ 2s-n(\frac{1}{p}-\frac{1}{q}) } \frac{1}{ ((z-1)^2-r^2 ) \sqrt{(z+1)^2-r^2}    } \left(\frac{ (z-1)^2 -r^2   }{ (z+1)^2 -r^2 } \right)^2  dr \\
& \leq  &
\int_{  (z,r) \in Z_1(\varepsilon ,z) }  r^{ 2s-n(\frac{1}{p}-\frac{1}{q}) }  \frac{1}{  \sqrt{(z+1)^2-r^2}    } \frac{ 1   }{  (z+1)^2 -r^2  }   dr \\
& \leq  &
\int_{  (z,r) \in Z_1(\varepsilon ,z) }  r^{ 2s-n(\frac{1}{p}-\frac{1}{q}) }  \frac{1}{  \sqrt{(z+1)^2-r^2}    } \frac{ 1   }{  4z  }   dr\\
& \leq  &
C  z^{ 2s-n(\frac{1}{p}-\frac{1}{q}) -1}   \quad \mbox{\rm for all}  \quad z \in [1,\infty) \,.
\eeqst
Finally,
\beqst
&  &
\int_{  (z,r) \in Z_1(\varepsilon ,z) }  dr\,  \frac{1}{     ((z-1)^2-r^2 ) \sqrt{(z+1)^2-r^2}    } r^{ 2s-n(\frac{1}{p}-\frac{1}{q}) }  \\
&  &
\times \left| (1-z^{2}+r^2 ) 
F\Big(-\frac{1}{2},\frac{1}{2};1; \frac{ (z-1)^2 -r^2   }{ (z+1)^2 -r^2 }   \Big) 
+   2  (z -1) F\Big(\frac{1}{2},\frac{1}{2};1; \frac{ (z-1)^2 -r^2   }{ (z+1)^2 -r^2 }   \Big)    \right| \\
& \leq   &
C z^{ 2s-n(\frac{1}{p}-\frac{1}{q})}   \quad \mbox{\rm for all}  \quad z \in [1,\infty) \,.
\eeqst
In the second zone we use (\ref{10.4n}), (\ref{FGH+xat1}),   and  (\ref{10.8}). Thus, we have to estimate the next two integrals:
\beqst 
A_{13} 
& :=   &
\int_{  (z,r) \in Z_2(\varepsilon ,z) }   r^{ 2s-n(\frac{1}{p}-\frac{1}{q}) } \frac{1}{     ((z-1)^2-r^2 ) \sqrt{(z+1)^2-r^2}    }  \left| (1-z^{2}+r^2 ) 
 \right|  dr \,, \\
A_{14} 
& :=   & 
\int_{  (z,r) \in Z_2(\varepsilon ,z) }   r^{ 2s-n(\frac{1}{p}-\frac{1}{q}) }  \frac{1}{     ((z-1)^2-r^2 ) \sqrt{(z+1)^2-r^2}    } 
\left| (z -1) \left(1+ \ln  z     \right)   \right| dr \,.
\eeqst 
We apply (\ref{10.4n}) to $A_{13}$ and obtain
\beqst 
A_{13} 
& \leq   &
C_\varepsilon \int_{  (z,r) \in Z_2(\varepsilon ,z) }   r^{ 2s-n(\frac{1}{p}-\frac{1}{q}) }   \frac{ 1   }{ [(z+1)^2 -r^2] }\frac{1}{    \sqrt{(z+1)^2-r^2}    }  \left| 
z^{2}-1- r^2 \right|  dr\\
& \leq   &
C_\varepsilon \int_{  0 }^{z-1}   r^{ 2s-n(\frac{1}{p}-\frac{1}{q}) }   \frac{1}{    \sqrt{(z+1)^2-r^2}    }   dr\\
& \leq   &
C_\varepsilon  z^{ 2s-n(\frac{1}{p}-\frac{1}{q}) }   \quad \mbox{\rm for all}  \quad z \in [1,\infty) \,, 
\eeqst 
while
\beqst 
A_{14} 
& \leq    & 
 (z -1) \left(1+ \ln  z     \right) \int_{  (z,r) \in Z_2(\varepsilon ,z) }  r^{ 2s-n(\frac{1}{p}-\frac{1}{q}) }   \frac{1}{ ((z-1)^2-r^2 ) \sqrt{(z+1)^2-r^2}    }  dr \\
& \leq    & 
C_\varepsilon  (z -1) \left(1+ \ln  z     \right) \int_{  0 }^{z-1}  r^{ 2s-n(\frac{1}{p}-\frac{1}{q}) }    \frac{1}{ ((z+1)^2-r^2 )^{3/2}   }  dr  \,.
\eeqst 
For $0 \geq a>-1 $ and $z\geq M$ the following integral can be easily estimated: 
\beqst 
&  &
\int_{  0 }^{z-1}  r^{ a}    \frac{1}{ ((z+1)^2-r^2 )^{3/2}   }  dr \\
& = &
\int_{  0 }^{z/2}  r^{ a}    \frac{1}{ ((z+1)^2-r^2 )^{3/2}   }  dr  + \int_{  z/2 }^{z-1}  r^{ a}    \frac{1}{ ((z+1)^2-r^2 )^{3/2}   }  dr  \\
& \leq  &
\frac{16}{9} \int_{  0 }^{z/2}  r^{ a } z^{  -3}    dr  + \frac{z^{ a}}{4^{ a}} \int_{  z/2 }^{z-1}     \frac{1}{ ((z+1)^2-r^2 )^{3/2}   }  dr  \\
& \leq  &
C z^{  -3+ a+1}    r^{ }     dr  + C  z^{ a -3/2}   \\
& \leq  &
C  z^{ a -3/2}  \,.
\eeqst 
Then 
$A_{14} 
  \leq      
C_\varepsilon  (z -1) \left(1+ \ln  z     \right) 
z^{ a -3/2} \leq  C_\varepsilon 
z^{ a }$.
The proposition is proven.
\hfill $\square$

\end{document}